\documentclass[13pt]{article}
\usepackage{latexsym}
\usepackage{geometry}
\usepackage{graphicx}
\usepackage{amsmath, amssymb, amsthm}
\usepackage{booktabs}
\usepackage{algorithm}
\usepackage{algorithmic}

\usepackage[utf8]{inputenc} 
\usepackage[T1]{fontenc}    

\usepackage{url}
\usepackage{natbib}

\usepackage{appendix}

\usepackage{amsfonts}
\usepackage{multirow}
\usepackage{multicol}

\usepackage{color}
\usepackage{algorithm}
\usepackage{algorithmic}

\usepackage{xcolor,colortbl}
\definecolor{LightCyan}{rgb}{0.88,1,1}

\usepackage{graphicx}
\usepackage{subfig}
\usepackage{hyperref}
\usepackage{tcolorbox}

\newtheorem{theorem}{Theorem}
\newtheorem{lemma}{Lemma}

\newtheorem{definition}{Definition}
\newtheorem{assumption}{Assumption}
\newtheorem{remark}{Remark}

\begin{document}
\title{ On Momentum-Based Gradient Methods for Bilevel Optimization with Nonconvex Lower-Level}
\author{Feihu Huang\thanks{Feihu Huang is with College of Computer Science and Technology,
Nanjing University of Aeronautics and Astronautics, Nanjing, China;
and also with MIIT Key Laboratory of Pattern Analysis and Machine Intelligence, Nanjing, China. Email: huangfeihu2018@gmail.com} }

\date{\textbf{The first version: 7 Mar. 2023}}

\maketitle

\begin{abstract}
Bilevel optimization is a popular two-level hierarchical optimization, which has been widely applied to many machine learning tasks such as hyperparameter learning, meta learning and continual learning.
Although many bilevel optimization methods
 recently have been developed, the bilevel methods are not well studied when the lower-level problem is nonconvex. To fill this gap, in the paper, we study a class of nonconvex
bilevel optimization problems, where both upper-level and lower-level problems are nonconvex, and the lower-level problem satisfies Polyak-{\L}ojasiewicz (PL) condition.
We propose an efficient momentum-based gradient bilevel method (MGBiO) to solve these deterministic problems.
Meanwhile, we propose a class of efficient momentum-based stochastic gradient bilevel methods (MSGBiO and  VR-MSGBiO) to solve these stochastic problems.
Moreover, we provide a useful convergence analysis framework for our methods. Specifically, under some mild conditions, we prove that our MGBiO method has a sample (or gradient) complexity of $O(\epsilon^{-2})$ for finding an $\epsilon$-stationary solution of the deterministic bilevel problems (i.e., $\|\nabla F(x)\|\leq \epsilon$),
which improves the existing best results by a factor of $O(\epsilon^{-1})$.
Meanwhile, we prove that our MSGBiO and VR-MSGBiO methods have sample complexities of $\tilde{O}(\epsilon^{-4})$ and $\tilde{O}(\epsilon^{-3})$, respectively, in finding an
$\epsilon$-stationary solution of the stochastic bilevel problems
(i.e., $\mathbb{E}\|\nabla F(x)\|\leq \epsilon$),
which improves the existing best results by a factor of $\tilde{O}(\epsilon^{-3})$.
Extensive experimental results on bilevel PL game and hyper-representation learning
demonstrate the efficiency of our algorithms. This paper commemorates the mathematician Boris Polyak (1935 -2023).
\end{abstract}

\section{Introduction}
Bilevel Optimization (BO)~\citep{colson2007overview}, as an effective two-level hierarchical optimization, can effectively capture the inherent nested structures of optimization problems.
Thus, it is useful in a variety of machine learning tasks such as hyperparameter learning~\citep{franceschi2018bilevel},
 meta learning~\citep{franceschi2018bilevel,ji2021bilevel}, continual learning~\citep{pham2021contextual,liu2022bome} and reinforcement learning~\citep{hong2020two}.
In the paper, we consider the following BO problem:
\begin{align} \label{eq:1}
 \min_{x \in \mathcal{X} } & \ F(x) \equiv f\big(x,y^*(x)\big) & \mbox{(Upper-Level)}  \\
 \mbox{s.t.} & \ y^*(x) \in \arg\min_{y\in \mathbb{R}^p} \ g\big(x,y\big), &  \mbox{(Lower-Level)} \nonumber
\end{align}
where $\mathcal{X}\subseteq \mathbb{R}^d$ denotes a convex set, and $F(x)\equiv f\big(x,y^*(x)\big):\mathbb{R}^d\rightarrow \mathbb{R}$ is the objective function of the upper-level, which is possibly nonconvex, and $y^*(x)$ may be a solution set when the lower level subproblem has multi local minimizers. Here $g(x,y)$ denotes the objective function of the lower-level, which is possibly nonconvex for all $x \in \mathcal{X}$ and satisfies Polyak-{\L}ojasiewicz (PL) condition~\citep{polyak1963gradient}. The PL condition
 relaxes the strong convexity and is more appealing than
convexity because some modern over-parameterized deep neural networks~\citep{frei2021proxy,song2021subquadratic}
have been shown to satisfy the PL-inequality.
In fact, Problem~(\ref{eq:1}) is associated to some machine learning applications:

\begin{table*}
  \centering
  \caption{ Comparison of sample (or gradient) \textbf{complexity }
between our algorithms and the existing bilevel algorithms in solving bilevel problems \textbf{with
lower-level non-convexity} for finding an $\epsilon$-stationary solution ($\|\nabla F(x)\|\leq \epsilon$ or its equivalent variants). Here $g(x,\cdot)$ denotes the function on the second variable $y$ with fixing variable $x$. The local PL condition implies that $g(x,\cdot)$ has multi local minimizers. \textbf{Note that} the GALET~\citep{xiao2023generalized} method simultaneously uses the PL condition, its Assumption 2 (i.e., let $\sigma_g = \inf_{x,y}\{\sigma_{\min}^{+}(\nabla^2_{yy} g(x,y))\} >0$ for all $(x,y)$) and its Assumption 1 (i.e., $\nabla^2_{yy} g(x,y)$ is Lipschitz continuous). Clearly, when Hessian matrix $\nabla^2_{yy} g(x,y)$ is singular, its Assumption 1 and Assumption 2 imply that the lower bound of the non-zero singular values $\sigma_g$ is close to zero (i.e., $\sigma_g\rightarrow 0$), under this case, the convergence results of the GALET are \textbf{meaningless}, e.g., the constant $L_w = \frac{\ell_{f,1}}{\sigma_g}+\frac{\sqrt{2}\ell_{g,2}\ell_{f,0}}{\sigma_g^2}\rightarrow + \infty$ used in its Lemmas 6 and 9. Under the other case, the PL condition, Lipschitz continuous of Hessian and its Assumption 2 (the singular values of Hessian is bounded away from 0, i.e., $\sigma_g>0$) imply that the GALET uses the strongly convex. \textbf{SC} stands for strongly convex. }
  \label{tab:1}
   \resizebox{\textwidth}{!}{
\begin{tabular}{c|c|c|c|c|c}
  \hline
  \textbf{Problem} & \textbf{Algorithm} & \textbf{Reference} & \textbf{Assumption} $g(x,\cdot)$ & \textbf{Non-Asymptotic} & \textbf{Complexity}   \\ \hline
   \multirow{6}*{Deterministic} & BVFIM/BVFSM  & \cite{liu2021avalue,liu2021towards,liu2021value}  & Nonconvex & $\times$ & $\times$ \\  \cline{2-6}
  & BGS-Opt  & \cite{arbel2022non}  &  Morse-Bott &$\times$ & $\times$ \\  \cline{2-6}
  & BOME  & \cite{liu2022bome} & PL / local-PL & $\surd$ & $O(\epsilon^{-3})$ / $O(\epsilon^{-4})$\\  \cline{2-6}
  & V-PBGD  & \cite{shen2023penalty} & PL / local-PL & $\surd$ & $O(\epsilon^{-3})$ / $O(\epsilon^{-3})$ \\  \cline{2-6}
  & GALET  & \cite{xiao2023generalized} & SC / PL & $\surd$ & $O(\epsilon^{-2})$ / \textbf{Meaningless} \\  \cline{2-6}
  & MGBiO  & Ours & PL / local-PL & $\surd$ &  {\color{red}{$O(\epsilon^{-2})$}} / {\color{red}{$O(\epsilon^{-2})$}} \\  \hline \hline
 \multirow{3}*{Stochastic} &  V-PBSGD & \cite{shen2023penalty}  & PL & $\surd$ & $\tilde{O}(\epsilon^{-6})$ \\  \cline{2-6}
   & MSGBiO  & Ours  & PL & $\surd$ & $\tilde{O}(\epsilon^{-4})$ \\  \cline{2-6}
  & VR-MSGBiO  & Ours  & PL & $\surd$ & {\color{red}{$\tilde{O}(\epsilon^{-3})$}} \\  \hline
\end{tabular}
 }
\end{table*}

\textbf{1). Bilevel Polyak-{\L}ojasiewicz Game.} In this task, its goal to find the optimal variables $x$ and $y$ by solving the following bilevel problem:
\begin{align}
 & \min_{x\in \mathbb{R}^d} \frac{1}{2}x^TPx + x^TR^1y, \\
 & \mbox{s.t.} \ \min_{y\in \mathbb{R}^d} \frac{1}{2}y^TQy + x^TR^2y, \nonumber
\end{align}
where both matrices $P\in \mathbb{R}^{d\times d}$ and $Q\in \mathbb{R}^{d\times d}$ are rank-deficient, and
$R^1, R^2 \in \mathbb{R}^{d\times d}$. Since the matrices $P$ and $Q$ are singular, the lower-level and upper-level objective functions may be not convex, while they satisfy the PL condition.

\textbf{2). Few-Shot Mate Learning}.
We consider the few-shot meta-learning problem with $k$ tasks
$\{\mathcal{T}_i,i=1,\cdots,k\}$ sampled from a task distribution, where each task
$\mathcal{T}_i$ has a loss function $\ell(x, y_i; \xi)$ over data $\xi$,
where $x \in \mathbb{R}^d$ denotes the parameters of an embedding model shared
by all tasks, and $y_i\in \mathbb{R}^p$ denotes the task-specific parameters. Our goal is to find
good parameters $x$ for all
tasks, and the optimal task-specific parameters $\{y_i\}_{i=1}^k$
for each task by minimizing its loss.

Here we take this model training as a bilevel optimization problem. In the lower-level problem, the
learner of task $\mathcal{T}_i$ is to find an optimal $y_i$
as the minimizer of its loss over a training set $\mathcal{D}^i_{tr}$.
In the upper-level problem, the meta learner to search for a good shared parameter vector
$y$ for all tasks over a test set $\mathcal{D}_{te}=\{\mathcal{D}^1_{te},\cdots,\mathcal{D}^k_{te}\}$ based on the minimizers $\{y_i\}_{i=1}^k$.
For notational simplicity, let $\tilde{y}=(y_1,\cdots,y_k)$ denotes all task-specific
parameters. Then, this bilevel problem can be represented by
\begin{align}
& \min_x \ \ell_{\mathcal{D}_{te}}(x,\tilde{y}) = \frac{1}{k}\sum_{i=1}^k \sum_{\xi\in\mathcal{D}^i_{te} }\ell_{te}(x,y^*_i;\xi) \\
& \mbox{s.t.} \ y^* = \arg\min_{\tilde{y}} \ell_{\mathcal{D}_{tr}} (x,\tilde{y}) =\frac{1}{k}\sum_{i=1}^k \sum_{\xi\in\mathcal{D}^i_{tr} }\ell_{tr}(x,y_i;\xi), \nonumber
\end{align}
where $\mathcal{D}^i_{tr}$ and $\mathcal{D}^i_{te}$ denote
the training and test datasets of task $\mathcal{T}_i$, respectively.
In the lower-level problem, we find each $y_i$ by minimizing
the task-specific loss $\ell_{tr}(x,y_i;\xi)$ for $i=1,\cdots,k$.
In practice, $y_i$ often are the parameters of \textbf{the
last layer} of a neural network and $x$ are the parameters of \textbf{the remaining layers} (e.g., some convolutional layers). Here we use the composition of
strongly convex functions with the single leaky ReLU neuron $\ell_{tr}(x,y_i;\xi) = g(\sigma(\phi(x)y_i);\xi)$ as the last layer, where $g(\cdot)$ is a strongly convex function, and $\sigma(\cdot)$ is a single leaky ReLU neuron, i.e., $\sigma(z) = \max(cz,z)$ with $c\neq 0$, and $\phi(x)$ denotes the output of the second last layer of a neural network. Following \cite{charles2018stability}, the function $\ell_{tr}(x,y_i;\xi)$ satisfies the PL condition, while
the function $\ell_{te}(x,y^*_i;\xi)$ is generally nonconvex.

The goal of Problem~(\ref{eq:1}) is to minimize an upper-level objective function, whose variables include the solution of lower-level problem.
When $g\big(x,y\big)=0$, Problem~(\ref{eq:1}) will reduce to a standard single-level constrained minimization problem.
Compared with these single-level constrained problems, the major difficulty of Problem~(\ref{eq:1})
is that the minimization of the upper and lower-level objectives are intertwined via the
mapping $y^*(x)$.
To alleviate this difficulty, recently many bilevel methods~\citep{ghadimi2018approximation,hong2020two,ji2021bilevel,chen2022single,huang2022enhanced}
for BO have been proposed by imposing the
strong convexity assumption on the lower-level objective.
The underlying idea of these methods is based on the algorithm of gradient descent
\begin{align}
 x_{t+1} = x_t -\eta\nabla F(x_t),
\end{align}
where $\eta$ is a step size, and $\nabla F(x_t)$ is a hyper-gradient of $F(x_t)$ defined as
\begin{align} \label{eq:g1}
 \nabla F(x_t) = \nabla_xf(x_t,y^*(x_t)) - \nabla^2_{xy}g(x_t,y^*(x_t))\nabla^2_{yy}g(x_t,y^*(x_t))^{-1}\nabla_yf(x_t,y^*(x_t)).
\end{align}
In fact, there are two major difficulties in estimating $\nabla F(x_t)$. The first obstacle is that we need to
estimate the optimal solution $y^*(x)$ of the lower problem for every given $x$, which updates the lower variable $y$ multiple times before updating $x$. To tackle this issue, several methods~\citep{ghadimi2018approximation,hong2020two,chen2021closing,ji2021bilevel,ji2022will} have been proposed to effectively track $y^*(x)$ without waiting for too many inner iterations before updating $x$.
The second obstacle is that we need to compute second-order
derivatives of $g(x,y)$. Many existing methods~\citep{chen2022single,dagreou2022framework} require an explicit extraction of second-order information
of $g(x,y)$ with a major focus on efficiently estimating its Jacobian and inverse Hessian.
Meanwhile, some methods~\citep{likhosherstov2020ufo,chen2022single,li2022fully,dagreou2022framework} avoid directly estimating its second-order computation and only use the first-order information of both upper and lower objectives.

The above studies including both algorithmic design and theoretical
investigation heavily rely on a restrictive strong convexity condition of lower-level problem.
In fact, many machine learning applications such as hyperparameter learning, meta learning and continual learning do not satisfy this condition. For example, when training deep neural networks, the lower-level problems of the continual learning~\citep{pham2021contextual} are nonconvex. However, the above methods cannot be easily adapted to the setting where the lower-level strong convex condition is absent.
More recently, thus, some works~\citep{liu2020generic,liu2021towards,sow2022constrained,arbel2022non,liu2022bome,
chen2023bilevel,liu2023averaged,lu2023first,shen2023penalty} have begun to study the bilevel optimization problems without lower-level strong convexity condition. These works basically can be divided into two categories: The first class mainly focuses on the lower level \textbf{convexity instead of strongly-convexity}. For example,
\cite{liu2020generic} proposed a bilevel descent aggregation (BDA) algorithmic framework for BO with the lower level convexity, and studied the asymptotic convergence analysis. Subsequently, \cite{sow2022constrained} proposed a constrained method for bilevel optimization with multiple inner minima, and the non-asymptotic convergence analysis. More recently, \cite{liu2023averaged} proposed an effective averaged method of multipliers for bilevel optimization
with lower-level convexity. \cite{chen2023bilevel} proposed the inexact gradient-free method (IGFM) to solve the bilevel problem based on zeroth order oracle, and provided the convergence properties of the IGFM in finding Goldstein stationary point. \cite{lu2023first} studied the BO when the lower-level objective is convex and nonsmooth.

The second class mainly focuses on the lower level \textbf{non-convexity instead of strongly-convexity}.
For example, \cite{liu2021towards,liu2021value,liu2021avalue,mehra2021penalty} studied the gradient-based method for bilevel optimization when the lower-level objective is nonconvex, and provided the asymptotic convergence analysis.
\cite{arbel2022non} studied the bilevel optimization with Morse-Bott lower-level objective.
\cite{liu2022bome} proposed a simple first-order method for BO problems where the lower-level problem is nonconex and satisfies PL inequality. More recently, \cite{shen2023penalty} presented a penalty-based bilevel gradient descent (PBGD) algorithm and established its
finite-time convergence for the constrained bilevel problem, where the lower-level objectives are non-convexity
and satisfy PL condition.

Although these works studied the BO without Lower-Level Strong-Convexity (LLSC) condition, the
 BO without LLSC still is not well-studied. For example, some methods~\citep{liu2021towards,arbel2022non} do not provide non-asymptotic analysis, or the existing non-asymptotic analysis shows that these methods still suffer from a high sample complexity (please see Table~\ref{tab:1}). Thus, there exists a natural question:
\begin{center}
\begin{tcolorbox}
\textbf{ Could we develop more efficient and provably bilevel methods to solve the bilevel problems with nonconvex lower-level ? }
\end{tcolorbox}
\end{center}

In the paper, we provide an affirmative answer to the above question and propose
an efficient momentum-based gradient bilevel algorithm (i.e., MGBiO) to solve deterministic Problem \eqref{eq:1}, where
the lower-level problem is nonconvex and satisfies PL condition. Recently, since
most of the existing stochastic bilevel methods ~\citep{yang2021provably,guo2021stochastic,khanduri2021near,huang2021biadam,dagreou2022framework,kwon2023fully} still rely on the lower-level strongly-convexity assumption, meanwhile, we propose a class of momentum-based stochastic gradient bilevel methods (MSGBiO and  VR-MSGBiO) to solve the stochastic bilevel problem:
\begin{align}
 \min_{x \in \mathcal{X}} & \ F(x)=\mathbb{E}_{\xi}\big[f(x,y^*(x);\xi)\big],  &
 \mbox{(Upper-Level)} \label{eq:2} \\
 \mbox{s.t.} & \ y^*(x) \in \arg\min_{y\in \mathbb{R}^p} \ \mathbb{E}_{\zeta}\big[g(x,y;\zeta)\big],   &  \mbox{(Lower-Level)} \nonumber
\end{align}
where $\xi$ and $\zeta$ are random variables, and $\mathcal{X}\subseteq \mathbb{R}^d$ denotes a convex set. Here  $F(x)=f(x,y^*(x))=\mathbb{E}_{\xi}\big[f(x,y^*(x);\xi)\big]$ is possibly nonconvex, and $g(x,y)=\mathbb{E}_{\zeta}\big[g(x,y;\zeta)\big]$ is nonconvex on any $x \in \mathcal{X}$ and satisfies PL condition.
Our main contributions are four-fold:
\begin{itemize}
\item[(1)] We propose an efficient momentum-based gradient bilevel method (i.e., MGBiO) to solve deterministic Problem (\ref{eq:1}), where the lower-level problem is nonconvex and satisfies PL condition.
\item[(2)] We propose a class of efficient momentum-based stochastic gradient bilevel methods (i.e., MSGBiO and  VR-MSGBiO) to solve the stochastic Problem (\ref{eq:2}) based on momentum and variance reduced techniques.
\item[(3)] We provide a convergence analysis framework for our algorithms. Specifically, we prove that our MGBiO method has a faster convergence rate $O(\frac{1}{\sqrt{T}})$ (i.e., $\|\nabla F(x)\|\leq O(\frac{1}{\sqrt{T}})$), which improves the best known result by a factor of $O(\frac{1}{T^{1/6}})$. Meanwhile, our MGBiO reaches a lower sample complexity of $O(\epsilon^{-2})$ in finding $\epsilon$-stationary solution of Problem~(\ref{eq:1}) ($\|\nabla F(x)\|\leq \epsilon$).
\item[(4)] We prove that our MSGBiO and VR-MSGBiO methods have sample complexities of $\tilde{O}(\epsilon^{-4})$ and $\tilde{O}(\epsilon^{-3})$, respectively, in finding an $\epsilon$-stationary solution of the stochastic Problem~(\ref{eq:2}) (i.e., $\mathbb{E}\|\nabla F(x)\|\leq \epsilon$).
\end{itemize}

\subsection*{Notations}
$\|\cdot\|$ denotes the $\ell_2$ norm for vectors and spectral norm for matrices.
$\langle x,y\rangle$ denotes the inner product of two vectors $x$ and $y$.
Given function $f(x,y)$, $f(x,\cdot)$ denotes  function \emph{w.r.t.} the second variable with fixing $x$,
and $f(\cdot,y)$ denotes function \emph{w.r.t.} the first variable
with fixing $y$.
$a_m=O(b_m)$ denotes that $a_m \leq c b_m$ for some constant $c>0$. The notation $\tilde{O}(\cdot)$ hides logarithmic terms. $\nabla_x$ denotes the partial derivative on variable $x$.
For notational simplicity, let $\nabla^2_{xy}=\nabla_x\nabla_y$ and $\nabla^2_{yy}=\nabla_y\nabla_y$. $\Pi_C[\cdot]$ denotes Euclidean projection onto the set $\{x\in \mathbb{R}^d: \|x\| \leq C\}$ for some constant $C>0$.
$\hat{\Pi}_{C}[\cdot]$ denotes a
projection on the set $\{X\in \mathbb{R}^{d\times p}: \|X\| \leq C\}$ for some constant $C>0$.
$\mathcal{S}_{[\mu,L_g]}$ denotes a
projection on the set $\{X\in \mathbb{R}^{d\times d}: \mu \leq \varrho(X) \leq L_g\}$,
where $\varrho(\cdot)$ denotes the eigenvalue function. Both $\hat{\Pi}_{C}$ and $\mathcal{S}_{[\mu,L_g]}$ can be implemented by using Singular Value Decomposition (SVD) and thresholding the singular values.

\section{ Preliminaries }
In this section, we provide some mild assumptions and useful lemmas on the above bilevel problems~(\ref{eq:1}) and~(\ref{eq:2}).

\subsection{Assumptions}
In this subsection, we give some mild assumptions on the above problems~(\ref{eq:1}) and (\ref{eq:2}).

\begin{assumption} \label{ass:1}
 The function $g(x,\cdot)$ satisfies the $\mu$-Polyak-{\L}ojasiewicz (PL) condition for some $\mu>0$ if for any given $x\in \mathcal{X}$, it holds that
\begin{align}
 \|\nabla_y g(x,y)\|^2 \geq 2\mu \big(g(x,y)-\min_y g(x,y)\big), \quad \forall y \in \mathbb{R}^p.
\end{align}
\end{assumption}

\begin{assumption} \label{ass:2}
 The function $g(x,y)$ satisfies
\begin{align}
 {\color{blue}{\varrho\big(\nabla^2_{yy}g\big(x,y^*(x)\big)\big) \in [\mu, L_g]}},
\end{align}
where $y^*(x) \in \arg\min_y g(x,y)$, and $\varrho(\cdot)$ denotes the eigenvalue (or singular-value) function and $L_g\geq\mu>0$.
\end{assumption}

\begin{assumption} \label{ass:3}
The functions $f(x,y)$ and $g(x,y)$ satisfy
\begin{itemize}
 \item[(i)] {\color{blue}{$\|\nabla_yf(x,y^*(x))\|\leq C_{fy}$ }} and {\color{blue}{$\|\nabla^2_{xy}g(x,y^*(x))\|\leq C_{gxy}$ }} at $ y^*(x) \in \arg\min_y g(x,y)$;
 \item[(ii)] The partial derivatives $\nabla_x f(x,y)$ and $\nabla_y f(x,y)$ are $L_f$-Lipschitz continuous, i.e., for $x,x_1,x_2\in \mathcal{X}$ and $y,y_1,y_2 \in \mathbb{R}^p$,
     \begin{align}
      \|\nabla_x f(x_1,y)-\nabla_x f(x_2,y)\| \leq L_f\|x_1-x_2\|, \
      \|\nabla_x f(x,y_1)-\nabla_x f(x,y_2)\| \leq L_f\|y_1-y_2\|, \nonumber \\
       \|\nabla_y f(x_1,y)-\nabla_y f(x_2,y)\| \leq L_f\|x_1-x_2\|, \
      \|\nabla_y f(x,y_1)-\nabla_y f(x,y_2)\| \leq L_f\|y_1-y_2\|; \nonumber
     \end{align}
    \item[(iii)] The partial derivatives $\nabla_x g(x,y)$ and $\nabla_y g(x,y)$ are $L_g$-Lipschitz continuous.
\end{itemize}
\end{assumption}

\begin{assumption} \label{ass:4}
 The Jacobian matrix $\nabla^2_{xy}g(x,y)$ and Hessian matrix $\nabla^2_{yy}g(x,y)$ are $L_{gxy}$-Lipschitz and $L_{gyy}$-Lipschitz continuous, respectively. E.g.,
 for all $x,x_1,x_2 \in \mathcal{X}$, $y,y_1,y_2 \in \mathbb{R}^p$, we have
 \begin{align}
 &\|\nabla^2_{xy} g(x_1,y)-\nabla^2_{xy} g(x_2,y)\| \leq L_{gxy}\|x_1-x_2\|,  \  \|\nabla^2_{xy} g(x,y_1)-\nabla^2_{xy} g(x,y_2)\| \leq L_{gxy}\|y_1-y_2\|, \nonumber \\
 &\|\nabla^2_{yy} g(x_1,y)-\nabla^2_{yy} g(x_2,y)\| \leq L_{gyy}\|x_1-x_2\|,  \  \|\nabla^2_{yy} g(x,y_1)-\nabla^2_{yy} g(x,y_2)\| \leq L_{gyy}\|y_1-y_2\|. \nonumber
 \end{align}
 \end{assumption}

\begin{assumption} \label{ass:5}
 The function $F(x)=f(x,y^*(x))$ is bounded below in $\mathcal{X}$, \emph{i.e.,} $F^* = \inf_{x\in \mathcal{X}}F(x) > -\infty$.
\end{assumption}

Assumption~\ref{ass:1} is commonly used in bilevel optimization without the lower-level strongly-convexity~\citep{liu2022bome,shen2023penalty}.
 Assumption~\ref{ass:2} ensures the non-singularity of Hessian matrix $\nabla^2_{yy}g(x,y^*(x))$.
 \textbf{Note that} our Assumption~2 requires the non-singularity of $\nabla^2_{yy}g(x,y)$ {\color{blue}{only at the minimized point}} $y=y^*(x)\in \arg\min_y g(x,y)$.
 Since $y^*(x)\in \arg\min_y g(x,y)$, we can not have negative eigenvalues at the minimizer $y^*(x)$, so Assumption~2 assumes that $\varrho\big(\nabla^2_{yy}g\big(x,y^*(x)\big)\big) \in [\mu, L_g]$ instead of $\varrho\big(\nabla^2_{yy}g\big(x,y^*(x)\big)\big) \in [-L_g,-\mu] \cup [\mu, L_g]$.
Since $\nabla^2_{yy} g(x,y)$
is a symmetric matrix, its singular values are the absolute value of eigenvalues.
Hence, we also can use $\varrho(\cdot)$ to denote the singular-value function.

\textbf{Note that} the GALET~\citep{xiao2023generalized} method simultaneously uses the PL condition, its Assumption 2 (i.e., let $\sigma_g = \inf_{x,y}\{\sigma_{\min}^{+}(\nabla^2_{yy} g(x,y))\} >0$ for all $(x,y)$) and its Assumption 1 (i.e., $\nabla^2_{yy} g(x,y)$ is Lipschitz continuous). Clearly, when Hessian matrix $\nabla^2_{yy} g(x,y)$ is singular, its Assumption 1 and Assumption 2 imply that the lower bound of the non-zero singular values $\sigma_g$ is close to zero (i.e., $\sigma_g\rightarrow 0$), under \textbf{this case}, the convergence results of the GALET are \textbf{meaningless}, e.g., the constant $L_w = \frac{\ell_{f,1}}{\sigma_g}+\frac{\sqrt{2}\ell_{g,2}\ell_{f,0}}{\sigma_g^2}\rightarrow + \infty$ used in its Lemmas 6 and 9, and $L_F = \ell_{f,0}(\ell_{f,1}+\ell_{g,2})/\sigma_g \rightarrow + \infty$ used in its Lemma 12. Under \textbf{the other case}, the PL condition, Lipschitz continuous of Hessian and its Assumption 2 (the singular values of Hessian is bounded away from 0, i.e., $\sigma_g>0$) imply that the GALET uses the strongly convex.
For example, considering a PL function: $g(x,y) = y^2+\sin^2(y)+x$ on variable $y$, where both $x$ and $y$ are scalar variables, we have $\nabla^2_{yy} g(x,y)=2+2\cos^2(y)-2\sin^2(y)$ and
$0=y^*=\arg\min_y g(x,y)$, and then we get $\nabla^2_{yy} g(x,y^*)=\nabla^2_{yy} g(x,0) =4>0$ and $\nabla^2_{yy} g(x,\pi/2)=0$. Clearly, Assumption~2 of the GALET~\citep{xiao2023generalized} can not conform to this PL function, since $\nabla^2_{yy} g(x,\pi/2)=0$. Clearly, our Assumption~\ref{ass:2} is milder than the Assumption~2 of the GALET. Meanwhile, our Assumption~\ref{ass:2} is more reasonable in practice.

 Assumption~\ref{ass:3} is commonly appeared in bilevel optimization methods \citep{ghadimi2018approximation,ji2021bilevel,liu2022bome}.
 \textbf{Note that} our Assumption~3(i) assumes that
 $\|\nabla_yf(x,y)\|$ and $\|\nabla^2_{xy}g(x,y)\|$ are bounded \textbf{only at
 the minimizer} $y^*(x) \in \arg\min_y g(x,y)$, while \citep{ghadimi2018approximation,ji2021bilevel} assume that $\|\nabla_yf(x,y)\|$ and $\|\nabla^2_{xy}g(x,y)\|$ are bounded \textbf{at
 any} $y \in \mathbb{R}^p$ and $x\in \mathcal{X}$.
 Meanwhile, the BOME~\cite{liu2022bome} uses the stricter assumption that
  $\|\nabla f(x,y)\|$, $\|\nabla g(x,y)\|$, $|f(x,y)|$ and $|g(x,y)|$
 are bounded for any $(x,y)$ in its Assumption~3.
Assumption~\ref{ass:4} is also commonly used in bilevel optimization methods \citep{ghadimi2018approximation,ji2021bilevel}.
Assumption~\ref{ass:5} ensures the feasibility of the bilevel Problem~\eqref{eq:1}.

\subsection{Useful Lemmas}
In this subsection, based on the above assumptions, we give some useful lemmas.

\begin{lemma} \label{lem:1}
Under the above Assumption \ref{ass:2}, we have, for any $x\in \mathbb{R}^d$,
\begin{align}
 \nabla F(x) =\nabla_x f(x,y^*(x)) - \nabla^2_{xy} g(x,y^*(x))\Big[\nabla^2_{yy}g(x,y^*(x))\Big]^{-1}\nabla_y f(x,y^*(x)). \nonumber
\end{align}
\end{lemma}

From the above Lemma~\ref{lem:1}, we can get the form of hyper-gradient $\nabla F(x)$ is the same as the hyper-gradient given in~(\ref{eq:g1}). Since the Hessian matrix $\nabla^2_{yy}g(x,y)$ for all $(x,y)$ may be singular in the nonconvex lower-level problems, we define a new and useful hyper-gradient estimator:
\begin{align}
\hat{\nabla} f(x,y)=\nabla_xf(x,y) - {\color{blue}{\hat{\Pi}_{C_{gxy}}\big[\nabla^2_{xy}g(x,y)\big]}} \big({\color{blue}{\mathcal{S}_{[\mu,L_g]}\big[\nabla^2_{yy}g(x,y)\big]}}\big)^{-1}
{\color{blue}{\Pi_{C_{fy}}\big[\nabla_yf(x,y)\big]}},
\end{align}
which replaces the standard hyper-gradient estimator $\breve{\nabla} f(x,y)$ used in ~\cite{ghadimi2018approximation,ji2021bilevel} for the strongly-convex lower-level problems:
\begin{align}
\breve{\nabla} f(x,y)=\nabla_xf(x,y) - \nabla^2_{xy}g(x,y) \big(\nabla^2_{yy}g(x,y)\big)^{-1}\nabla_yf(x,y).
\end{align}
\begin{lemma} \label{lem:2}
Under the above Assumptions~\ref{ass:1}-\ref{ass:3}, the functions (or mappings) $F(x)=f(x,y^*(x))$, $G(x)=g(x,y^*(x))$ and $y^*(x)\in \arg\min_{y\in \mathbb{R}^p}g(x,y)$ satisfy, for all $x_1,x_2\in \mathbb{R}^d$,
\begin{align}
 & \|y^*(x_1)-y^*(x_2)\| \leq \kappa\|x_1-x_2\|, \quad \|\nabla y^*(x_1) - \nabla y^*(x_2)\| \leq L_y\|x_1-x_2\| \nonumber \\
 & \|\nabla F(x_1) - \nabla F(x_2)\|\leq L_F\|x_1-x_2\|, \quad \|\nabla G(x_1) - \nabla G(x_2)\|\leq L_G\|x_1-x_2\| \nonumber
\end{align}
where $\kappa=C_{gxy}/\mu$, $L_y=\big( \frac{C_{gxy}L_{gyy}}{\mu^2} +  \frac{L_{gxy}}{\mu} \big) (1+ \frac{C_{gxy}}{\mu})$,
$L_F=\Big(L_f + L_f\kappa + C_{fy}\big( \frac{C_{gxy}L_{gyy}}{\mu^2} +  \frac{L_{gxy}}{\mu} \big)\Big)(1+\kappa)$ and $L_G=\Big(L_g + L_g\kappa + C_{gy}\big( \frac{C_{gxy}L_{gyy}}{\mu^2} +  \frac{L_{gxy}}{\mu} \big)\Big)(1+\kappa)$.
\end{lemma}

\begin{lemma} \label{lem:3}
Let $\hat{\nabla} f(x,y)=\nabla_xf(x,y) - {\color{blue}{\hat{\Pi}_{C_{gxy}}\big[\nabla^2_{xy}g(x,y)\big] }} \big({\color{blue}{\mathcal{S}_{[\mu,L_g]}\big[\nabla^2_{yy}g(x,y)\big]}}\big)^{-1}{\color{blue}{\Pi_{C_{fy}}\big[\nabla_yf(x,y)\big]}}$ and
$\nabla F(x_t) = \nabla f(x_t,y^*(x_t))$,
 we have
 \begin{align}
 \|\hat{\nabla} f(x,y)-\nabla F(x)\|^2 \leq \hat{L}^2\|y^*(x)-y\|^2\leq \frac{2\hat{L}^2}{\mu}\big(g(x,y)-\min_y g(x,y)\big),
\end{align}
where $\hat{L}^2 = 4\big(L^2_f+ \frac{L^2_{gxy}C^2_{fy}}{\mu^2} + \frac{L^2_{gyy} C^2_{gxy}C^2_{fy}}{\mu^4} +
 \frac{L^2_fC^2_{gxy}}{\mu^2}\big)$.
\end{lemma}

\begin{algorithm}[tb]
\caption{ MGBiO Algorithm }
\label{alg:1}
\begin{algorithmic}[1] 
\STATE {\bfseries Input:} $T$, parameters $\{\gamma, \lambda, \eta_t\}$
and initial input $x_1 \in \mathcal{X}$ and $y_1 \in \mathbb{R}^p$; \\
\FOR{$t = 1, 2, \ldots, T$}
\STATE $v_t = \nabla_y g(x_t,y_t)$, $u_t = \nabla_x f(x_t,y_t)$, $h_t =\Pi_{C_{fy}}\big[\nabla_yf(x_t,y_t)\big]$, $G_t = \hat{\Pi}_{C_{gxy}}\big[\nabla^2_{xy}g(x_t,y_t)\big]$;
\STATE  {\color{blue}{ $H_t = \mathcal{S}_{[\mu,L_g]}\big[\nabla_{yy}^2 g(x_t,y_t)\big] = U_t\Theta_tU_t^T$ }}, where $\theta_{t,i}\in [\mu,L_g]$ for all $i=1,\cdots,p$;
\STATE $w_t=u_t - G_t(H_t)^{-1}h_t= \nabla_x f(x_t,y_t) - G_t\Big(\sum_{i=1}^p \big(U^T_{t,i}h_t \big)/\theta_{t,i}U_{t,i}\Big)$;
\STATE $\tilde{x}_{t+1} = \arg\min_{x\in \mathcal{X}}\big\{ \langle w_t, x\rangle + \frac{1}{2\gamma}\|x-x_t\|^2\big\}$ and
 $x_{t+1} = x_t+\eta_t(\tilde{x}_{t+1}-x_t)$;
\STATE $\tilde{y}_{t+1} = y_t-\lambda v_t$ and
 $y_{t+1} = y_t+\eta_t(\tilde{y}_{t+1}-y_t)$;
\ENDFOR
\STATE {\bfseries Output:} Chosen uniformly random from $\{x_t\}_{t=1}^{T}$.
\end{algorithmic}
\end{algorithm}

\section{Momentum-Based Gradient Bilevel Method}
In the section, we propose an efficient momentum-based gradient bilevel (MGBiO) method to solve the deterministic bilevel Problem~(\ref{eq:1}).

\subsection{Method}
The detailed procedure of our MGBiO method is provided in Algorithm~\ref{alg:1}.
At the line 5 of Algorithm~\ref{alg:1}, we estimate the gradient $\nabla F(x_t)$ based on the above Lemma~\ref{lem:1}:
\begin{align}
w_t = \hat{\nabla}f(x_t,y_t) =\nabla_xf(x_t,y_t) + \hat{\Pi}_{C_{gxy}}\big[\nabla^2_{xy}g(x_t,y_t)\big] \big(\mathcal{S}_{[\mu,L_g]}\big[\nabla^2_{yy}g(x_t,y_t)\big]\big)^{-1}\Pi_{C_{fy}}\big[\nabla_yf(x_t,y_t)\big].
\end{align}
Since the projection operator $\mathcal{S}_{[\mu,L_g]}$ is implemented by using Singular Value Decomposition (SVD) and thresholding the singular values, we have
\begin{align}
 H_t = \mathcal{S}_{[\mu,L_g]}\big[\nabla_{yy}^2 g(x_t,y_t)\big] = U_t\Theta_tU_t^T,
\end{align}
where $U_t$ is a real orthogonal matrix, i.e, $U_t^TU_t=U_tU_t^T=I_p$, and $\Theta_t=\mbox{diag}(\theta_t) \in \mathbb{R}^{p\times p}$ is a diagonal matrix with $\theta_{t,i} \in [\mu,L_g]$ for all $i=1,2,\cdots, p$. Clearly, we can easily obtain
\begin{align}
 H_t^{-1} = U_t\Theta_t^{-1}U_t^T, \quad \Theta_t^{-1} = \mbox{diag}(1/\theta_t),
\end{align}
where $1/\theta_t = (1/\theta_{t,1},\cdots,1/\theta_{t,p})$.
Thus, we can get
\begin{align} \label{eq:3}
 w_t=u_t - G_t(H_t)^{-1}h_t & = \nabla_x f(x_t,y_t) - G_tU_t\Theta_t^{-1}U_t^Th_t \nonumber \\
 & = \nabla_x f(x_t,y_t) - G_t \Big(\sum_{i=1}^p \big(U^T_{t,i}h_t \big)/\theta_{t,i}U_{t,i}\Big),
\end{align}
where $U_{t,i}$ is the $i$-th column of matrix $U_t = [U_{t,1},\cdots, U_{t,p}]$.
Based on the above equality~(\ref{eq:3}), calculations of the gradient $w_t$ mainly depends on the above SVD.
Recently, the existing many methods~\citep{tzeng2013split} can compute SVD of the large-scale matrices.
Thus, our algorithm can also solve the large-scale bilevel problems in high-dimension setting.

At the line 6 of Algorithm~\ref{alg:1}, when $\mathcal{X}\subset \mathbb{R}^d$, we have
$$\tilde{x}_{t+1}=\mathbb{P}_{\mathcal{X}}(x_t-\gamma\hat{\nabla}f(x_t,y_t))=\arg\min_{x\in \mathcal{X}}
\big\{ \langle \hat{\nabla}f(x_t,y_t), x-x_t\rangle + \frac{1}{2\gamma}\|x-x_t\|^2 \big\},$$
and then we can define the gradient mapping
$\mathcal{G}(x_t,\hat{\nabla}f(x_t,y_t),\gamma) = \frac{1}{\gamma}\big(x_t-\mathbb{P}_{\mathcal{X}}(x_t-\gamma\hat{\nabla}f(x_t,y_t))\big)$.
Meanwhile, we can also define a gradient mapping $\mathcal{G}(x_t,\nabla F(x_t),\gamma) = \frac{1}{\gamma}\big(x_t-\mathbb{P}_{\mathcal{X}}(x_t-\gamma\nabla F(x_t))\big)$ as a convergence metric also
used in~\citep{shen2023penalty}, where
$$\mathbb{P}_{\mathcal{X}}(x_t-\gamma\nabla F(x_t)) =\arg\min_{x\in \mathcal{X}}
\big\{ \langle \nabla F(x_t), x-x_t\rangle + \frac{1}{2\gamma}\|x-x_t\|^2 \big\}.$$

At the line 6 of Algorithm~\ref{alg:1}, we have
$\tilde{x}_{t+1}=x_t -\gamma w_t$ when $\mathcal{X}=\mathbb{R}^d$.
Meanwhile, we also use the
momentum iteration to get $x_{t+1}=x_t+\eta_t(\tilde{x}_{t+1}-x_t)$, where $\eta_t \in (0,1]$.
Since $\eta_t \in (0,1]$ and $x_t,\tilde{x}_{t+1}\in \mathcal{X}$, we have $x_{t+1}=(1-\eta_t)x_t+\eta_t\tilde{x}_{t+1}\in \mathcal{X}$.
At the line 7 of Algorithm~\ref{alg:1}, we also simultaneously use the gradient descent and momentum iterations to update the inner variable $y$.

Although our Algorithm~\ref{alg:1} requires to compute Hessian matrix, its inverse and
its projecting matrix over spectral set, \textbf{our algorithm can avoid computing inverse of Hessian matrix}. In fact, many bilevel optimization methods~\citep{ji2021bilevel,xiao2023generalized} still compute Hessian matrix $\nabla^2_{yy}g(x,y)$.
\textbf{Note that} the GALET~\citep{xiao2023generalized} method not only calculates Hessian matrices, but also computes the product of Hessian matrices multiple times when updating variable $x$ at each iteration.
From the following experimental results, our Algorithm~\ref{alg:1} outperforms the GALET~\citep{xiao2023generalized} method.

\subsection{Convergence Results}
In this subsection,
we provide the convergence results of our MGBiO method under the above mild assumptions.

\subsubsection{Convergence with Unimodal $g(x,y)$}
We first present the convergence properties of our MGBiO method when $g(x,\cdot)$ has \textbf{a unique minimizer} $y^*(x)=\arg\min_y f(x,y)$ and satisfies the PL condition for all $x$. We first give a useful lemma.

\begin{lemma} \label{lem:4}
Suppose the sequence $\{x_t,y_t\}_{t=1}^T$ be generated from Algorithm \ref{alg:1}.
Under the above Assumptions~\ref{ass:1}-\ref{ass:2}, given $\gamma\leq \frac{\lambda\mu}{8L_G}$ and $0<\lambda \leq \frac{1}{2L_g\eta_t}$ for all $t\geq 1$, we have
\begin{align}
g(x_{t+1},y_{t+1}) - G(x_{t+1})
& \leq (1-\frac{\eta_t\lambda\mu}{2}) \big(g(x_t,y_t) -G(x_t)\big) + \frac{\eta_t}{8\gamma}\|\tilde{x}_{t+1}-x_t\|^2  -\frac{\eta_t}{4\lambda}\|\tilde{y}_{t+1}-y_t\|^2 \nonumber \\
& \quad + \eta_t\lambda\|\nabla_y g(x_t,y_t)-v_t\|^2,
\end{align}
where $G(x_t)=g(x_t,y^*(x_t))$ with $y^*(x_t) \in \arg\min_{y}g(x_t,y)$ for all $t\geq 1$.
\end{lemma}

The above Lemma~\ref{lem:4} shows the properties of the residuals $g(x_t,y_t)-G(x_t)\geq 0$ for all $t\geq 1$. In fact, this lemma can also be used in the following Algorithms~\ref{alg:2} and \ref{alg:3}.

In convergence analysis,
we measure the convergence of MGBiO by the following \emph{Lyapunov} function (i.e., potential function),
 for any $t\geq 1$
\begin{align}
 \Omega_t = F(x_t) + g(x_t,y_t)-G(x_t),
\end{align}
where $G(x_t)=g(x_t,y^*(x_t))=\min_y g(x,y)$.

\begin{theorem}  \label{th:1}
 Under the above Assumptions (\ref{ass:1}-\ref{ass:5}), in the Algorithm \ref{alg:1}, let $\eta_t=\eta$ for all $t\geq 0$, $0< \gamma \leq \min\big(\frac{1}{2L_F\eta},\frac{\lambda\mu^2}{16\hat{L}^2} \big)$, $0<\lambda\leq \frac{1}{2L_g\eta}$. When $\mathcal{X}\subseteq\mathbb{R}^d$, we can get
 \begin{align}
 \frac{1}{T}\sum_{t=1}^T\|\mathcal{G}(x_t,\nabla F(x_t),\gamma)\|
 & \leq \frac{1}{T} \sum_{t=1}^T \big[ \|w_t-\nabla F(x_t)\| + \frac{1}{\gamma}\|\tilde{x}_{t+1}-x_t\|  \big] \leq \frac{4\sqrt{R}}{\sqrt{3T\gamma\eta}},
\end{align}
when $\mathcal{X}=\mathbb{R}^d$, we can get
 \begin{align}
 \frac{1}{T}\sum_{t=1}^T\|\nabla F(x_t)\|
  & \leq \frac{1}{T} \sum_{t=1}^T \big[ \|w_t-\nabla F(x_t)\| + \frac{1}{\gamma}\|\tilde{x}_{t+1}-x_t\|  \big] \leq \frac{4\sqrt{R}}{\sqrt{3T\gamma\eta}},
\end{align}
where $R= F(x_1) -F^* +g(x_1,y_1)-G(x_1)$.
\end{theorem}

\begin{remark}
Without loss of generality, let $\eta_t=\eta \in (0,1)$, i.e., $\eta=O(1)$ and $\gamma=\min\big(\frac{1}{2L_F\eta},\frac{\lambda\mu^2}{16\hat{L}^2} \big)=O(1)$.
Meanwhile, we have $R= F(x_1) -F^* +g(x_1,y_1)-G(x_1)=O(1)$. Then we have \\ $\frac{1}{T}\sum_{t=1}^T\|\mathcal{G}(x_t,\nabla F(x_t),\gamma)\| \leq O(\frac{1}{\sqrt{T}}) \leq \epsilon$ or $\frac{1}{T}\sum_{t=1}^T\|\nabla F(x_t)\| \leq O(\frac{1}{\sqrt{T}}) \leq \epsilon$. Thus we can get a lower sample (gradient) complexity of $1*T=O(\epsilon^{-2})$ in finding $\epsilon$-stationary solution of Problem~(\ref{eq:1}).
\end{remark}

\subsubsection{Convergence with multimodal $g(x,y)$}
Here, we give the convergence properties of our MGBiO method when $g(x,\cdot)$ has \textbf{multi local minimizers } $y^\diamond(x,y)\in \arg\min_y f(x,y)$ and satisfies the local PL condition for all $x$.
As~\cite{liu2022bome}, we first define the attraction points:

\begin{definition}
Given any $(x,y)$, if the sequence $\{y_t\}_{t=0}^{\infty}$ generated by gradient descent $y_t=y_{t-1}-\lambda \nabla_y g(x,y_{t-1})$ starting from $y_0=y$ converges to $y^\diamond(x,y)$, we say that $y^\diamond(x,y)$ is
the attraction point of $(x,y)$ with step size $\lambda>0$.
\end{definition}

Meanwhile, an attraction basin be formed by the same attraction point in set of $(x,y)$. The following analysis needs to assume the PL condition within the individual attraction basins.

\begin{assumption} \label{ass:1g}
(\textbf{Local PL Condition in Attraction Basins})
Assume that for any $(x,y)$, $y^\diamond(x,y)$ exists.
 $g(x,\cdot)$ satisfies the local PL condition in attraction basins, 
i.e., for any $(x,y)$, there exists a constant $\mu>0$ such that
\begin{align}
 \|\nabla_y g(x,y)\|^2 \geq 2\mu \big(g(x,y)- G^\diamond(x)\big),
\end{align}
where $G^\diamond(x)=g(x,y^\diamond(x,y))$.
\end{assumption}

\begin{assumption} \label{ass:2g}
 The function $g\big(x,y^\diamond(x,y)\big)$ satisfies
\begin{align}
 \varrho\big(\nabla^2_{yy}g\big(x,y^\diamond(x,y)\big)\big) \in [\mu, L_g],
\end{align}
where $y^\diamond(x,y)$ is the attraction point of $(x,y)$, and $\varrho(\cdot)$ denotes the eigenvalue (or singular-value) function  and $L_g\geq\mu>0$.
\end{assumption}

\begin{assumption} \label{ass:3g}
The functions $f(x,y)$ and $g(x,y)$ satisfy
\begin{itemize}
 \item[(i)] {\color{blue}{$\|\nabla_yf(x,y^\diamond(x,y))\|\leq C_{fy}$ }} and {\color{blue}{$\|\nabla^2_{xy}g(x,y^\diamond(x,y))\|\leq C_{gxy}$ }} at $ y^\diamond(x,y) \in \arg\min_y g(x,y)$;
 \item[(ii)] The partial derivatives $\nabla_x f(x,y)$ and $\nabla_y f(x,y)$ are $L_f$-Lipschitz continuous;
 \item[(iii)] The partial derivatives $\nabla_x g(x,y)$ and $\nabla_y g(x,y)$ are $L_g$-Lipschitz continuous.
\end{itemize}
\end{assumption}

\begin{assumption} \label{ass:5g}
 The function $F^\diamond(x)=f(x,y^\diamond(x,y))$ is bounded below in $\mathcal{X}$, \emph{i.e.,} $F^\diamond = \inf_{x\in \mathcal{X}}F^\diamond(x) > -\infty$.
\end{assumption}

\begin{theorem}  \label{th:1g}
 Under the above Assumptions (\ref{ass:4},\ref{ass:1g},\ref{ass:2g},\ref{ass:3g},\ref{ass:5g}), in the Algorithm \ref{alg:1}, let $\eta_t=\eta$ for all $t\geq 0$, $0< \gamma \leq \min\big(\frac{1}{2L_F\eta},\frac{\lambda\mu^2}{16\hat{L}^2} \big)$, $0<\lambda\leq \frac{1}{2L_g\eta}$. When $\mathcal{X}\subseteq\mathbb{R}^d$, we can get
 \begin{align}
 \frac{1}{T}\sum_{t=1}^T\|\mathcal{G}(x_t,\nabla F^\diamond(x_t),\gamma)\|
 & \leq \frac{1}{T} \sum_{t=1}^T \big[ \|w_t-\nabla F^\diamond(x_t)\| + \frac{1}{\gamma}\|\tilde{x}_{t+1}-x_t\|  \big] \leq \frac{4\sqrt{R}}{\sqrt{3T\gamma\eta}},
\end{align}
when $\mathcal{X}=\mathbb{R}^d$, we can get
 \begin{align}
 \frac{1}{T}\sum_{t=1}^T\|\nabla F^\diamond(x_t)\|
  & \leq \frac{1}{T} \sum_{t=1}^T \big[ \|w_t-\nabla F^\diamond(x_t)\| + \frac{1}{\gamma}\|\tilde{x}_{t+1}-x_t\|  \big] \leq \frac{4\sqrt{R}}{\sqrt{3T\gamma\eta}},
\end{align}
where $R= F(x_1) -F^\diamond +g(x_1,y_1)-G^\diamond(x_1)$.
\end{theorem}

\begin{remark}
The proof of Theorem~\ref{th:1g} \textbf{can follow} the proof of Theorem~\ref{th:1}.
Meanwhile, we can also get a lower sample (gradient) complexity of $1*T=O(\epsilon^{-2})$ in finding $\epsilon$-stationary solution of Problem~(\ref{eq:1}) \textbf{under local PL condition}.
\end{remark}

\begin{algorithm}[tb]
\caption{ MSGBiO Algorithm }
\label{alg:2}
\begin{algorithmic}[1] 
\STATE {\bfseries Input:} $T$, parameters $\{\gamma, \lambda, \eta_t, \alpha_t, \hat{\alpha}_t, \tilde{\alpha}_t, \beta_t, \hat{\beta}_t\}$
and initial input $x_1 \in \mathcal{X}$ and $y_1 \in \mathbb{R}^p$; \\
\STATE {\bfseries initialize:} Draw two independent samples $\xi_1$ and $\zeta_1$,
and then compute $u_1 = \nabla_x f(x_1,y_1;\xi_1)$, $h_1 =\Pi_{C_{fy}}\big[ \nabla_y f(x_1,y_1;\xi_1) \big]$, $v_1 = \nabla_y g(x_1,y_1;\zeta_1)$, $G_1 = \hat{\Pi}_{C_{gxy}}\big[\nabla_{xy}^2 g(x_1,y_1;\zeta_1)\big]$, {\color{blue}{$H_1 = \mathcal{S}_{[\mu,L_g]}\big[\nabla_{yy}^2 g(x_1,y_1;\zeta_1)\big]$}} and $w_1 = u_1 - G_1(H_1)^{-1}h_1$;  \\
\FOR{$t = 1, 2, \ldots, T$}
\STATE $\tilde{x}_{t+1} = \arg\min_{x\in \mathcal{X}}\big\{ \langle w_t, x\rangle + \frac{1}{2\gamma}\|x-x_t\|^2\big\}$ and
 $x_{t+1} = x_t+\eta_t(\tilde{x}_{t+1}-x_t)$;
\STATE $\tilde{y}_{t+1} = y_t-\lambda v_t$ and
 $y_{t+1} = y_t+\eta_t(\tilde{y}_{t+1}-y_t)$;
\STATE Draw two independent samples $\xi_{t+1}$ and $\zeta_{t+1}$;
\STATE $u_{t+1} = \beta_{t+1}\nabla_x f(x_{t+1},y_{t+1};\xi_{t+1}) + (1-\beta_{t+1})u_t $;
\STATE $h_{t+1} = \Pi_{C_{fy}}\big[\hat{\beta}_{t+1}\nabla_y f(x_{t+1},y_{t+1};\xi_{t+1}) + (1-\hat{\beta}_{t+1})h_t \big]$;
\STATE $v_{t+1} = \alpha_{t+1}\nabla_y g(x_{t+1},y_{t+1};\zeta_{t+1}) + (1-\alpha_{t+1})v_t $;
\STATE $G_{t+1} = \hat{\Pi}_{C_{gxy}}\big[ \hat{\alpha}_{t+1}\nabla^2_{xy} g(x_{t+1},y_{t+1};\zeta_{t+1}) + (1-\hat{\alpha}_{t+1})G_t \big]= W_{t+1}\Sigma_{t+1}V_{t+1}^T$;
\STATE {\color{blue}{$H_{t+1} = \mathcal{S}_{[\mu,L_g]}\big[ \tilde{\alpha}_{t+1}\nabla^2_{yy} g(x_{t+1},y_{t+1};\zeta_{t+1}) + (1-\tilde{\alpha}_{t+1})H_t\big]= U_{t+1}\Theta_{t+1}U_{t+1}^T$ }} ;
\STATE $w_{t+1}=u_{t+1} - G_{t+1}(H_{t+1})^{-1}h_{t+1}= u_{t+1} - W_{t+1}\Sigma_{t+1}V_{t+1}^T\Big(\sum_{i=1}^p \big(U^T_{t+1,i}h_{t+1} \big)/\theta_{t+1,i}U_{t+1,i}\Big)$;
\ENDFOR
\STATE {\bfseries Output:} Chosen uniformly random from $\{x_t\}_{t=1}^{T}$.
\end{algorithmic}
\end{algorithm}

\begin{algorithm}[tb]
\caption{ VR-MSGBiO Algorithm }
\label{alg:3}
\begin{algorithmic}[1] 
\STATE {\bfseries Input:} $T$, parameters $\{\gamma, \lambda, \eta_t, \alpha_t, \hat{\alpha}_t, \tilde{\alpha}_t, \beta_t, \hat{\beta}_t\}$
and initial input $x_1 \in \mathcal{X}$ and $y_1 \in \mathbb{R}^p$; \\
\STATE {\bfseries initialize:} Draw two independent samples $\xi_1$ and $\zeta_1$,
and then compute $u_1 = \nabla_x f(x_1,y_1;\xi_1)$, $h_1 =\Pi_{C_{fy}}\big[ \nabla_y f(x_1,y_1;\xi_1) \big]$, $v_1 = \nabla_y g(x_1,y_1;\zeta_1)$, $G_1 = \hat{\Pi}_{C_{gxy}}\big[\nabla_{xy}^2 g(x_1,y_1;\zeta_1)\big]$, ${\color{blue}{H_1 = \mathcal{S}_{[\mu,L_g]}\big[\nabla_{yy}^2 g(x_1,y_1;\zeta_1)\big]}}$ and $w_1 = u_1 - G_1(H_1)^{-1}h_1$;  \\
\FOR{$t = 1, 2, \ldots, T$}
\STATE $\tilde{x}_{t+1} = \arg\min_{x\in \mathcal{X}}\big\{ \langle w_t, x\rangle + \frac{1}{2\gamma}\|x-x_t\|^2\big\}$ and
 $x_{t+1} = x_t+\eta_t(\tilde{x}_{t+1}-x_t)$;
\STATE $\tilde{y}_{t+1} = y_t-\lambda v_t$ and
 $y_{t+1} = y_t+\eta_t(\tilde{y}_{t+1}-y_t)$;
\STATE Draw two independent samples $\xi_{t+1}$ and $\zeta_{t+1}$;
\STATE $u_{t+1} = \nabla_x f(x_{t+1},y_{t+1};\xi_{t+1}) + (1-\beta_{t+1})(u_t-\nabla_x f(x_{t},y_{t};\xi_{t+1})) $;
\STATE $h_{t+1} = \Pi_{C_{fy}}\big[\nabla_y f(x_{t+1},y_{t+1};\xi_{t+1}) + (1-\hat{\beta}_{t+1})(h_t-\nabla_y f(x_{t},y_{t};\xi_{t+1})) \big]$;
\STATE $v_{t+1} = \nabla_y g(x_{t+1},y_{t+1};\zeta_{t+1}) + (1-\alpha_{t+1})(v_t-\nabla_y g(x_{t},y_{t};\zeta_{t+1})) $;
\STATE $G_{t+1} = \hat{\Pi}_{C_{gxy}}\big[ \nabla^2_{xy} g(x_{t+1},y_{t+1};\zeta_{t+1}) + (1-\hat{\alpha}_{t+1})(G_t-\nabla^2_{xy} g(x_{t},y_{t};\zeta_{t+1}) )\big] = W_{t+1}\Sigma_{t+1}V_{t+1}^T$;
\STATE ${\color{blue}{H_{t+1} = \mathcal{S}_{[\mu,L_g]}\big[ \nabla^2_{yy} g(x_{t+1},y_{t+1};\zeta_{t+1}) + (1-\tilde{\alpha}_{t+1})(H_t-\nabla^2_{yy} g(x_t,y_t;\zeta_{t+1}))\big]= U_{t+1}\Theta_{t+1}U_{t+1}^T}}$;
\STATE $w_{t+1}=u_{t+1} - G_{t+1}(H_{t+1})^{-1}h_{t+1}= u_{t+1} - W_{t+1}\Sigma_{t+1}V_{t+1}^T\Big(\sum_{i=1}^p \big(U^T_{t+1,i}h_{t+1} \big)/\theta_{t+1,i}U_{t+1,i}\Big)$;
 \ENDFOR
\STATE {\bfseries Output:} Chosen uniformly random from $\{x_t\}_{t=1}^{T}$.
\end{algorithmic}
\end{algorithm}

\section{Momentum-Based Stochastic Gradient Bilevel Methods}
In this section, we present a class of efficient momentum-based stochastic gradient bilevel methods (i.e., MSGBiO and  VR-MSGBiO) to solve the stochastic bilevel problem (\ref{eq:2}) based on momentum and variance reduced techniques.

\subsection{Methods}
In the subsection, we first propose a momentum-based stochastic gradient bilevel method (MSGBiO) method by simultaneously using momentum technique to estimate gradient and update variables. The detailed procedure of our MSGBiO is provided in Algorithm~\ref{alg:2}.

At the lines 4-5 of Algorithm~\ref{alg:2}, we use the momentum iteration to update variables
$x$ and $y$, i.e., $x_{t+1} = x_t+\eta_t(\tilde{x}_{t+1}-x_t)$ and $y_{t+1} = y_t+\eta_t(\tilde{y}_{t+1}-y_t)$.
Note that we use the same tuning parameter $\eta_t$ in updating variables
$x$ and $y$. At the lines 7-11 of Algorithm~\ref{alg:2}, we use the standard momentum technique
to estimate the stochastic
gradients or Jacobian matrices. Here we use the projected operators, e.g., $\mathcal{S}_{[\mu,L_g]}[\cdot]$ is a
projection on the set $\{X\in \mathbb{R}^{d\times d}: \mu \leq \varrho(X) \leq L_g\}$,
where $\varrho(\cdot)$ denotes the eigenvalue function.

Since both $\hat{\Pi}_{C}$ and $\mathcal{S}_{[\mu,L_g]}$ can be implemented by using Singular Value Decomposition (SVD) and thresholding the singular values, we have for all $t\geq 1$
\begin{align}
 H_t = U_t\Theta_tU_t^T, \quad G_t = W_t\Sigma_tV_t^T,
\end{align}
where $U_t$ is a real orthogonal matrix, i.e, $U_t^TU_t=U_tU_t^T=I_p$, and $\Theta_t=\mbox{diag}(\theta_t) \in \mathbb{R}^{p\times p}$ is a diagonal matrix with $\theta_{t,i} \in [\mu,L_g]$ for all $i=1,2,\cdots, p$, and
$W_t$ and $V_t$ are orthogonal matrices, i.e., $W_tW_t^T=I_d$ and $V_tV_t^T = I_p$, and $\Sigma_t \in \mathbb{R}^{d\times p}$ is a rectangular diagonal matrix with $\max_{1\leq i \leq \max(d,p)} \Sigma_{t,i,i} \leq C_{gxy}$.
Then we can get for all $t\geq 1$
\begin{align}
 w_t=u_t - G_t(H_t)^{-1}h_t = u_t - W_t\Sigma_tV_t^T\Big(\sum_{i=1}^p \big(U^T_{t,i}h_t \big)/\theta_{t,i}U_{t,i}\Big),
\end{align}
where $U_{t,i}$ is the $i$-th column of matrix $U_t = [U_{t,1},\cdots, U_{t,p}]$.

Next, we propose a variance-reduced stochastic gradient bilevel method (VR-MSGBiO) method by using the variance-reduced momentum-based technique of STORM~\citep{cutkosky2019momentum}/ProxHSGD~\citep{tran2022hybrid} to estimate gradients
and Jacobian matrix. Algorithm~\ref{alg:3} shows the detailed procedure of our VR-MSGBiO algorithm.
Specifically, at the lines 7-11 of Algorithm~\ref{alg:3}, we use the variance-reduced momentum technique to estimate the gradients or Jacobian matrices. For example, we use the gradient estimator $u_{t+1}$ to estimate gradient $\nabla_x f(x_{t+1},y_{t+1})$. At the line 7 of Algorithm~\ref{alg:2}, we use the \textbf{basic momentum-based} gradient estimator to estimate $\nabla_x f(x_{t+1},y_{t+1})$:
\begin{align}
u_{t+1} = \beta_{t+1}\nabla_x f(x_{t+1},y_{t+1};\xi_{t+1}) + (1-\beta_{t+1})u_t, \label{eq:26a}
\end{align}
while at the line 7 of Algorithm~\ref{alg:3}, we use the \textbf{variance-reduced momentum-based} gradient estimator:
\begin{align}
 u_{t+1} & = \nabla_x f(x_{t+1},y_{t+1};\xi_{t+1}) + (1-\beta_{t+1})(u_t-\nabla_x f(x_{t},y_{t};\xi_{t+1})) \nonumber \\
 & = \beta_{t+1}\nabla_x f(x_{t+1},y_{t+1};\xi_{t+1}) + (1-\beta_{t+1})(u_t + {\color{blue}{\nabla_x f(x_{t+1},y_{t+1};\xi_{t+1})-\nabla_x f(x_{t},y_{t};\xi_{t+1})}} ). \label{eq:27a}
\end{align}
Compared the above equalities~(\ref{eq:26a}) with (\ref{eq:27a}), the additional term $\nabla_x f(x_{t+1},y_{t+1};\xi_{t+1})-\nabla_x f(x_{t},y_{t};\xi_{t+1})$ is a variance reduced term used in SARAH~\citep{nguyen2017sarah} and
SPIDER~\citep{fang2018spider,wang2019spiderboost}.

\subsection{Convergence Results}
In this subsection,
we provide the convergence results of our MSGBiO and VR-MSGBiO methods.
We first provide some mild assumptions for Problem~(\ref{eq:2}).

\begin{assumption} \label{ass:6}
The functions $f(x,y)=\mathbb{E}_{\xi}[f(x,y;\xi)]$, $g(x,y)=\mathbb{E}_{\zeta}[g(x,y;\zeta)]$, $f(x,y;\xi)$ and $g(x,y;\zeta)$ satisfy
\begin{itemize}
 \item[(i)] {\color{blue}{$\|\nabla_yf(x,y^*(x))\|\leq C_{fy}$}} and {\color{blue}{$\|\nabla^2_{xy}g(x,y^*(x))\|\leq C_{gxy}$}} at $y^*(x) \in \arg\min_y g(x,y)$; 
 \item[(ii)] The partial derivatives $\nabla_x f(x,y;\xi)$ and $\nabla_y f(x,y;\xi)$ are $L_f$-Lipschitz continuous, i.e., for $x,x_1,x_2\in \mathcal{X}$ and $y,y_1,y_2 \in \mathbb{R}^p$,
     \begin{align}
      \|\nabla_x f(x_1,y;\xi)-\nabla_x f(x_2,y;\xi)\| \leq L_f\|x_1-x_2\|, \
      \|\nabla_x f(x,y_1;\xi)-\nabla_x f(x,y_2;\xi)\| \leq L_f\|y_1-y_2\|, \nonumber \\
      \|\nabla_y f(x_1,y;\xi)-\nabla_y f(x_2,y;\xi)\| \leq L_f\|x_1-x_2\|, \
      \|\nabla_y f(x,y_1;\xi)-\nabla_y f(x,y_2;\xi)\| \leq L_f\|y_1-y_2\|; \nonumber
     \end{align}
  \item[(iii)] The partial derivatives $\nabla_xg(x,y;\zeta)$ and $\nabla_yg(x,y;\zeta)$ are $L_g$-Lipschitz continuous.
\end{itemize}
\end{assumption}

\begin{assumption} \label{ass:7}
 The partial derivatives $\nabla^2_{xy}g(x,y;\zeta)$ and $\nabla^2_{yy}g(x,y;\zeta)$ are $L_{gxy}$-Lipschitz and $L_{gyy}$-Lipschitz, e.g.,
 for all $x,x_1,x_2 \in \mathcal{X}$ and $y,y_1,y_2 \in \mathbb{R}^p$
 \begin{align}
 \|\nabla^2_{xy} g(x_1,y;\zeta)-\nabla^2_{xy} g(x_2,y;\zeta)\| \leq L_{gxy}\|x_1-x_2\|,  \  \|\nabla^2_{xy} g(x,y_1;\zeta)-\nabla^2_{xy} g(x,y_2;\zeta)\| \leq L_{gxy}\|y_1-y_2\|. \nonumber
 \end{align}
 \end{assumption}

\begin{assumption} \label{ass:8}
Stochastic partial derivatives $\nabla_x f(x,y;\xi)$, $\nabla_y f(x,y;\xi)$, $\nabla_x g(x,y;\zeta)$, $\nabla_y g(x,y;\zeta)$, $\nabla^2_{xy} g(x,y;\zeta)$ and $\nabla^2_{yy} g(x,y;\zeta)$ are unbiased
with bounded variance, e.g.,
\begin{align}
\mathbb{E}[\nabla_x f(x,y;\xi)] = \nabla_x f(x,y), \ \mathbb{E}\|\nabla_x f(x,y;\xi) - \nabla_x f(x,y) \|^2 \leq \sigma^2. \nonumber
\end{align}
\end{assumption}

Assumptions~\ref{ass:6}-\ref{ass:8} are commonly used in the stochastic bilevel optimization problems \citep{ghadimi2018approximation,hong2020two,yang2021provably,khanduri2021near,huang2022fast}.
Based on Assumptions~\ref{ass:6} and~\ref{ass:7}, we have $\|\nabla_x f(x_1,y)-\nabla_x f(x_2,y)\|=\|\mathbb{E}[\nabla_x f(x_1,y;\xi)-\nabla_x f(x_2,y;\xi)]\|
\leq \mathbb{E}\|\nabla_x f(x_1,y;\xi)-\nabla_x f(x_2,y;\xi)\| \leq L_f\|x_1-x_2\|$, i.e., we can also obtain $\nabla_x f(x,y)$ is $L_f$-Lipschitz continuous,
which is similar for $\nabla_y f(x,y)$, $\nabla_y g(x,y)$, $\nabla^2_{xy}g(x,y)$ and $\nabla^2_{yy}g(x,y)$.
Clearly, Assumptions~\ref{ass:3} and~\ref{ass:4} is milder than Assumptions~\ref{ass:6} and~\ref{ass:7}.
Next, we provide a useful lemma:
\begin{lemma} \label{lem:5}
When the gradient estimator $w_t$ generated from Algorithm~\ref{alg:2} or~\ref{alg:3}, for all $t\geq 1$,
 we have
\begin{align}
 \|w_t-\nabla F(x_t)\|^2
 & \leq 8\|u_t - \nabla_xf(x_t,y_t)\|^2 + \frac{8C^2_{fy}}{\mu^2}\|G_t- {\color{blue}{\hat{\Pi}_{C_{gxy}}\big[\nabla^2_{xy}g(x_t,y_t)\big]}}\|^2   \nonumber \\
 & \quad + 8\kappa^2 \|h_t -{\color{blue}{\Pi_{C_{fy}}\big[\nabla_yf(x_t,y_t)\big]}} \|^2 + \frac{8\kappa^2 C^2_{fy}}{\mu^2}\|H_t-
 {\color{blue}{\mathcal{S}_{[\mu,L_g]}\big[\nabla^2_{yy} g(x_t,y_t)\big]}}\|^2  \nonumber \\
 & \quad + \frac{4\hat{L}^2}{\mu}\big(g(x_t,y_t)-G(x_t)\big),
\end{align}
where $\kappa=\frac{C_{gxy}}{\mu}$.
\end{lemma}

\subsubsection{Convergence Results of MSGBiO}

In the convergence analysis,
we measure the convergence of MSGBiO by the following \emph{Lyapunov} function (i.e., potential function),
for any $t\geq 1$
\begin{align}
 \Phi_t & = \mathbb{E}\Big [F(x_t) + g(x_t,y_t)-G(x_t) + \gamma \big(\|\nabla_x f(x_t,y_t)-u_t\|^2 + \|{\color{blue}{\Pi_{C_{fy}}\big[\nabla_yf(x_t,y_t)\big]}}-h_t\|^2 \nonumber \\
 & \qquad + \|{\color{blue}{\hat{\Pi}_{C_{gxy}}\big[\nabla^2_{xy}g(x_t,y_t)\big]}} - G_t\|^2 + \|{\color{blue}{\mathcal{S}_{[\mu,L_g]}\big[\nabla^2_{yy} g(x_t,y_t)\big]}} - H_t\|^2\big) + \lambda\|\nabla_y g(x_t,y_t) - v_t\|^2 \Big]. \nonumber
\end{align}

\begin{theorem} \label{th:2}
 Under the above Assumptions (\ref{ass:1}-\ref{ass:5}, \ref{ass:8}), in the Algorithm \ref{alg:2}, let $\eta_t=\frac{k}{(m+t)^{1/2}}$ for all $t\geq 0$, $\beta_{t+1}=c_1\eta_t$, $\hat{\beta}_{t+1}=c_2\eta_t$, $\alpha_{t+1}=c_3\eta_t$, $\hat{\alpha}_{t+1}=c_4\eta_t$, $\tilde{\alpha}_{t+1}=c_5\eta_t$, $10 \leq c_1 \leq \frac{m^{1/2}}{k}$, $10\kappa^2 \leq c_2 \leq \frac{m^{1/2}}{k}$, $1 \leq c_3 \leq \frac{m^{1/2}}{k}$, $\frac{10C^2_{fy}}{\mu^2} \leq c_4 \leq \frac{m^{1/2}}{k}$, $\frac{10C^2_{fy}\kappa^2}{\mu^2} \leq c_5 \leq \frac{m^{1/2}}{k}$, $m\geq \max\big(k^2, (c_1k)^2,(c_2k)^2, (c_3k)^2, (c_4k)^2, (c_5k)^2\big)$, $k>0$, $0< \gamma \leq \min\big(\frac{m^{1/2}}{2L_Fk},\frac{\lambda\mu^2}{16\hat{L}^2},\frac{\sqrt{5}}{4\breve{L}},\frac{1}{32L^2_g\lambda},\frac{5}{8\breve{L}^2\lambda}\big)$ and $0<\lambda\leq \min\big(\frac{m^{1/2}}{2L_gk},\frac{1}{4L_g}\big)$. When $\mathcal{X}\subseteq \mathbb{R}^d$, we can get
\begin{align}
  \frac{1}{T}\sum_{t=1}^T\mathbb{E}\|\mathcal{G}(x_t,\nabla F(x_t),\gamma)\|
 & \leq \frac{1}{T} \sum_{t=1}^T \mathbb{E} \big[ \|w_t-\nabla F(x_t)\| + \frac{1}{\gamma}\|\tilde{x}_{t+1}-x_t\|  \big]   \leq \frac{\sqrt{2M}m^{1/4}}{\sqrt{T}} + \frac{\sqrt{2M}}{T^{1/4}};
\end{align}
When $\mathcal{X}= \mathbb{R}^d$, we can get
\begin{align}
  \frac{1}{T}\sum_{t=1}^T\mathbb{E}\|\nabla F(x_t)\|
 & \leq \frac{1}{T} \sum_{t=1}^T \mathbb{E} \big[ \|w_t-\nabla F(x_t)\| + \frac{1}{\gamma}\|\tilde{x}_{t+1}-x_t\|  \big]   \leq \frac{\sqrt{2M}m^{1/4}}{\sqrt{T}} + \frac{\sqrt{2M}}{T^{1/4}},
\end{align}
where $M= \frac{4(F(x_1)- F^*+g(x_1,y_1)-G(x_1))}{k\gamma} + \frac{16\sigma^2}{k} + \frac{4\lambda\sigma^2}{\gamma k} + \frac{16m\sigma^2\ln(m+T)}{k}+ \frac{4m\lambda\sigma^2\ln(m+T)}{k\gamma}$
and $\breve{L}^2=L^2_f+\frac{L^2_f}{\kappa^2}+\frac{\mu^2L^2_{gxy}}{C^2_{fy}}+\frac{\mu^2L^2_{gyy}}{C^2_{fy}\kappa^2}$.
\end{theorem}

\begin{remark}
Without loss of generality, let $k=O(1)$, $m=O(1)$, $c_1=O(1)$, $c_2=O(1)$, $c_3=O(1)$, $c_4=O(1)$ and $c_5=O(1)$.
 Then we have $M=\tilde{O}(1)$. Thus we have $\frac{1}{T}\sum_{t=1}^T\|\mathcal{G}(x_t,\nabla F(x_t),\gamma)\| \leq \tilde{O}(\frac{1}{T^{1/4}}) \leq \epsilon$ or $\frac{1}{T}\sum_{t=1}^T\|\nabla F(x_t)\| \leq \tilde{O}(\frac{1}{T^{1/4}}) \leq \epsilon$. Thus we can get a lower sample (gradient) complexity of $1*T=\tilde{O}(\epsilon^{-4})$ in finding $\epsilon$-stationary solution of Problem~(\ref{eq:2}).
\end{remark}

\subsubsection{ Convergence Results of VR-MSGBiO}
In the convergence analysis,
we measure the convergence of VR-MSGBiO by the following \emph{Lyapunov} function (i.e., potential function),
for any $t\geq 1$
 \begin{align}
 \Gamma_t & = \mathbb{E}\Big [F(x_t) + g(x_t,y_t)-G(x_t) + \frac{\gamma}{\eta_{t-1}} \big(\|\nabla_x f(x_t,y_t)-u_t\|^2 + \|{\color{blue}{\Pi_{C_{fy}}\big[\nabla_yf(x_t,y_t)\big]}}- h_t\|^2 \nonumber \\
 & \qquad  + \|{\color{blue}{ \hat{\Pi}_{C_{gxy}}\big[\nabla^2_{xy}g(x_t,y_t)\big] }} - G_t\|^2 + \|{\color{blue}{\mathcal{S}_{[\mu,L_g]}\big[\nabla^2_{yy} g(x_t,y_t)\big]}} - H_t\|^2\big) + \frac{\lambda}{\eta_{t-1}}\|\nabla_y g(x_t,y_t) - v_t\|^2 \Big]. \nonumber
 \end{align}

\begin{theorem}  \label{th:3}
 Under the above Assumptions (\ref{ass:1}, \ref{ass:2}, \ref{ass:5}, \ref{ass:6}-\ref{ass:8}), in the Algorithm \ref{alg:3}, let $\eta_t=\frac{k}{(m+t)^{1/3}}$ for all $t\geq 0$, $\beta_{t+1}=c_1\eta_t^2$, $\hat{\beta}_{t+1}=c_2\eta_t^2$, $\alpha_{t+1}=c_3\eta_t^2$, $\hat{\alpha}_{t+1}=c_4\eta_t^2$, $\tilde{\alpha}_{t+1}=c_5\eta_t^2$, $m\geq \max\big(2,k^3, (c_1k)^3,(c_2k)^3, (c_3k)^3, (c_4k)^3, (c_5k)^3\big)$, $k>0$, $c_1 \geq \frac{2}{3k^3} + 10$, $c_2 \geq \frac{2}{3k^3} + 10\kappa^2$, $c_3 \geq \frac{2}{3k^3} + 1$, $c_4 \geq \frac{2}{3k^3} + \frac{10C^2_{fy}}{\mu^2}$, $c_5 \geq \frac{2}{3k^3} + +  \frac{10C^2_{fy}\kappa^2}{\mu^2}$, $0< \gamma \leq \min\big(\frac{m^{1/3}}{2L_Fk},\frac{\lambda\mu^2}{16\hat{L}^2},\frac{1}{8\check{L}},\frac{1}{64L^2_g\lambda},
 \frac{1}{32\check{L}^2\lambda},\frac{\lambda\mu}{8L_G}\big)$, $0<\lambda\leq \min\big(\frac{1}{4\sqrt{2}L_g},\frac{m^{1/3}}{2L_gk}\big)$. When $\mathcal{X}\subseteq\mathbb{R}^d$,
 we can get
 \begin{align}
 \frac{1}{T}\sum_{t=1}^T\mathbb{E}\|\mathcal{G}(x_t,\nabla F(x_t),\gamma)\|
 \leq \frac{1}{T} \sum_{t=1}^T \mathbb{E} \big[ \|w_t-\nabla F(x_t)\| + \frac{1}{\gamma}\|\tilde{x}_{t+1}-x_t\|  \big]  \leq \frac{\sqrt{2\breve{M}}m^{1/6}}{\sqrt{T}} + \frac{\sqrt{2\breve{M}}}{T^{1/3}};
\end{align}
When $\mathcal{X}=\mathbb{R}^d$, we can get
 \begin{align}
 \frac{1}{T}\sum_{t=1}^T\mathbb{E}\|\nabla F(x_t)\|
  \leq \frac{1}{T} \sum_{t=1}^T \mathbb{E} \big[ \|w_t-\nabla F(x_t)\| + \frac{1}{\gamma}\|\tilde{x}_{t+1}-x_t\|  \big]   \leq \frac{\sqrt{2\breve{M}}m^{1/6}}{\sqrt{T}} + \frac{\sqrt{2\breve{M}}}{T^{1/3}},
\end{align}
where $\breve{M}= \frac{4(F(x_1)- F^*+g(x_1,y_1)-G(x_1))}{k\gamma} + \frac{16\sigma^2m^{1/3}}{k^2} + \frac{4\lambda\sigma^2m^{1/3}}{\gamma k^2} + \big( 2k^2\hat{c}^2\sigma^2 + \frac{2k^2c^2_3\lambda\sigma^2}{\gamma}\big)\ln(m+T)$,
 $\check{L}^2=2L^2_f+L^2_{gxy}+L^2_{gyy}$ and $\hat{c}^2=c^2_1+c^2_2+c^2_4+c^2_5$.
\end{theorem}

\begin{remark}
Without loss of generality, let $k=O(1)$, $m=O(1)$, $c_1=O(1)$, $c_2=O(1)$, $c_3=O(1)$, $c_4=O(1)$ and $c_5=O(1)$.
 Then we have $\breve{M}=\tilde{O}(1)$. Thus we have $\frac{1}{T}\sum_{t=1}^T\|\mathcal{G}(x_t,\nabla F(x_t),\gamma)\| \leq \tilde{O}(\frac{1}{T^{1/3}}) \leq \epsilon$ or $\frac{1}{T}\sum_{t=1}^T\|\nabla F(x_t)\| \leq \tilde{O}(\frac{1}{T^{1/3}}) \leq \epsilon$. Thus we can get a lower sample (gradient) complexity of $1*T=\tilde{O}(\epsilon^{-3})$ in finding $\epsilon$-stationary solution of Problem~(\ref{eq:2}).
\end{remark}

\begin{figure}[ht]
\centering
 \subfloat{\includegraphics[width=0.48\textwidth]{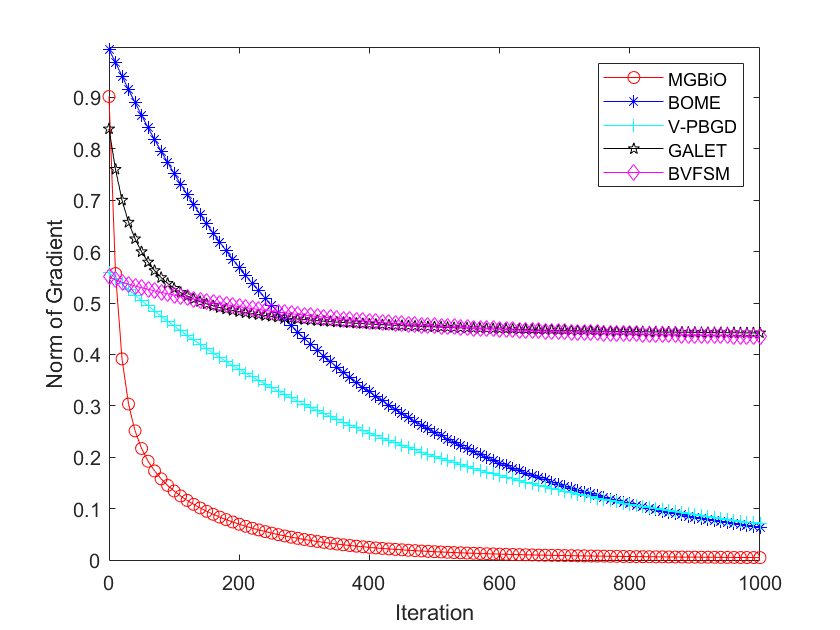}}
  \hfill
 \subfloat{\includegraphics[width=0.48\textwidth]{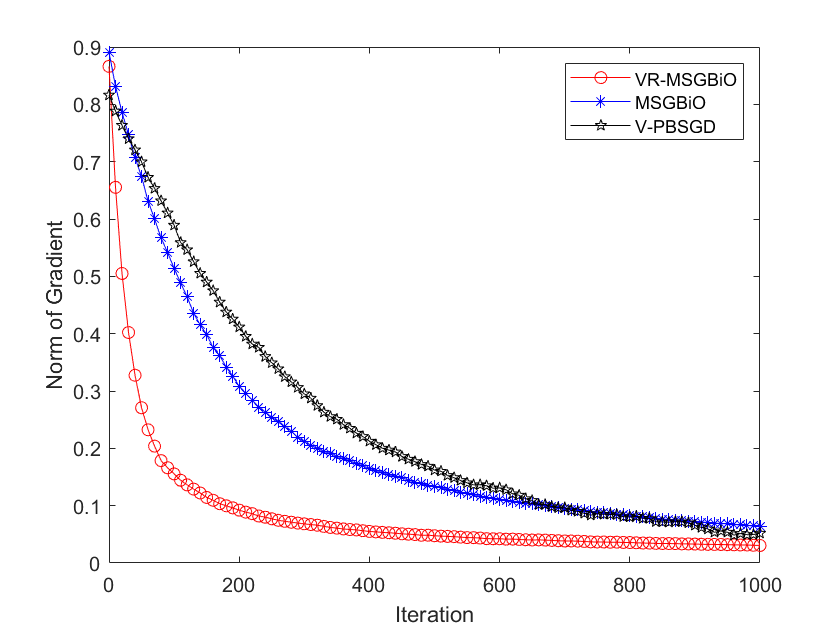}}
  \hfill
\caption{ Norm of gradient against number of iteration: \textbf{Left}: deterministic algorithms, \textbf{Right}: stochastic algorithms.}
\label{fig:1}
\end{figure}

\section{Numerical Experiments}
In the section, we conduct numerical experiments on bilevel PL game and
hyper-representation learning
to demonstrate the efficiency of our algorithms. In the experiments, we compare our algorithms
with the baselines in the above Table~\ref{tab:1}. Meanwhile, we also add two baselines:
BVFSM~\citep{liu2021avalue} and GALET~\citep{xiao2023generalized}.


\begin{figure}[ht]
\centering
 \subfloat{\includegraphics[width=0.33\textwidth]{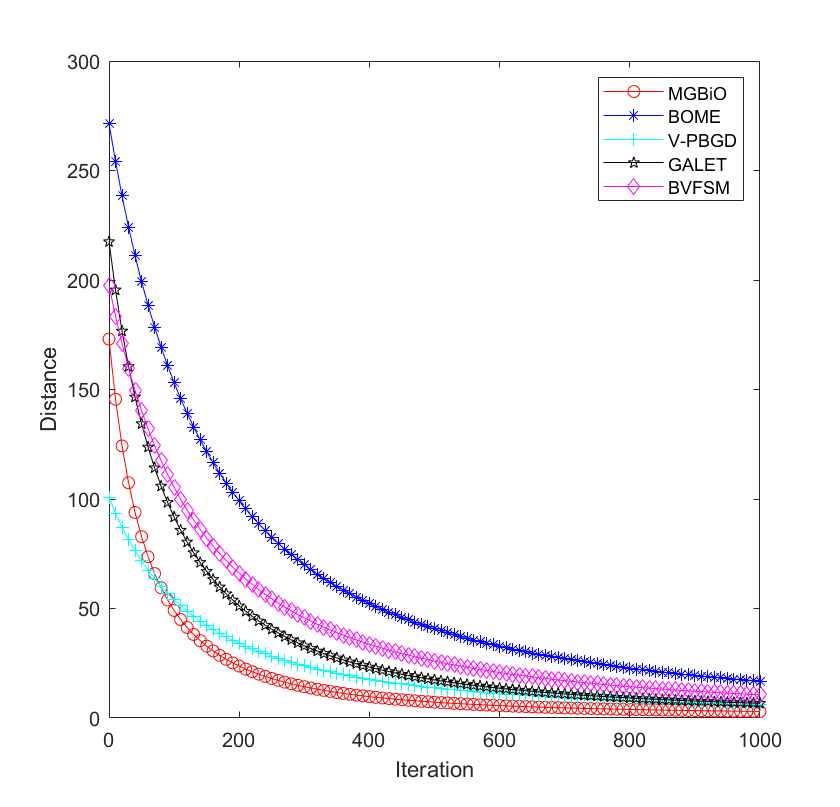}}
  \hfill
 \subfloat{\includegraphics[width=0.33\textwidth]{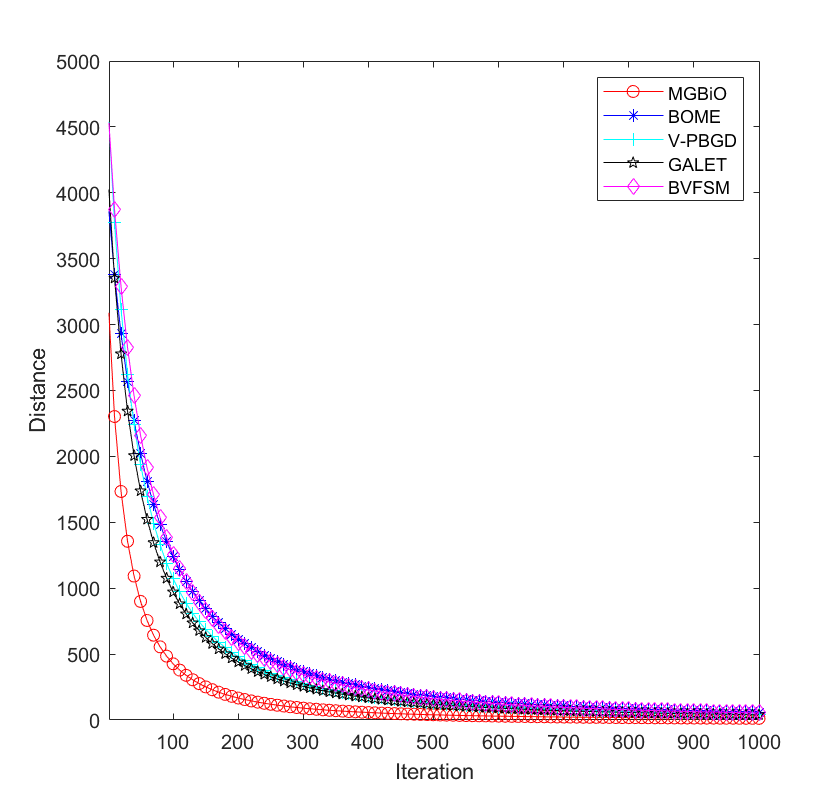}}
  \hfill
 \subfloat{\includegraphics[width=0.33\textwidth]{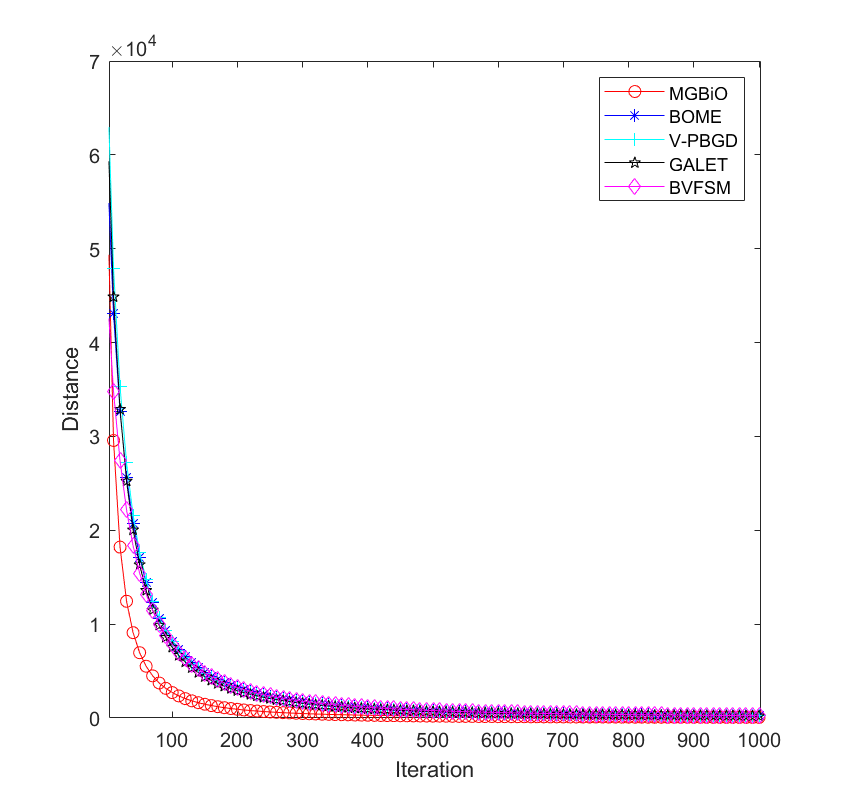}}
  \hfill
\caption{ The ratio of distance $\|UU^T-H^*\|^2_F/\|H^*\|^2_F$ of \textbf{deterministic} algorithms under the case of $d=50$ (Left), $d=100$ (Middle) and $d=200$ (Right).}
\label{fig:2}
\end{figure}

\begin{figure}[ht]
\centering
 \subfloat{\includegraphics[width=0.33\textwidth]{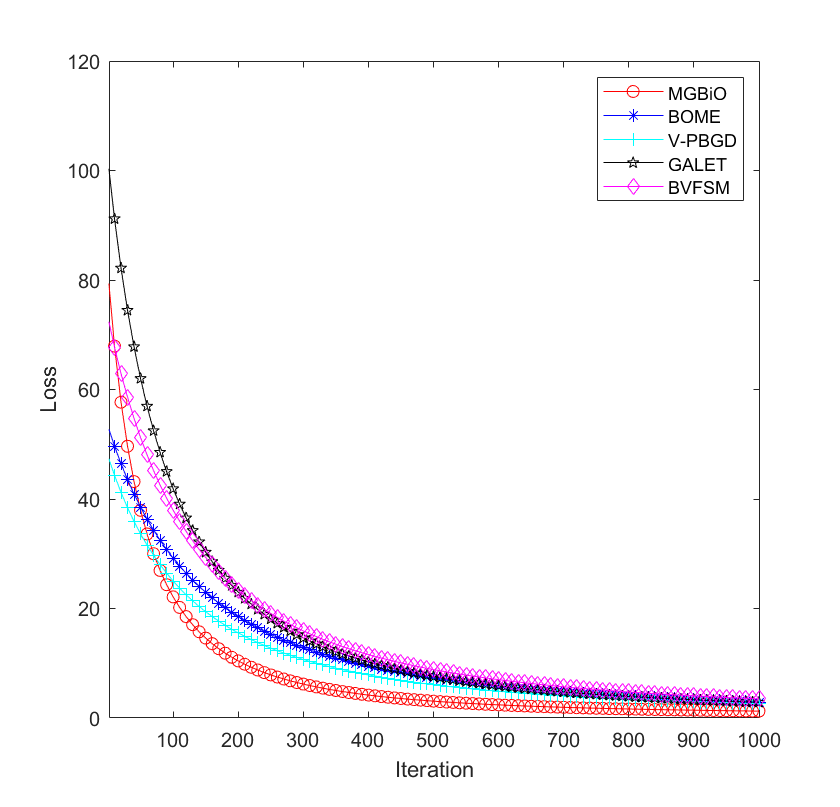}}
  \hfill
 \subfloat{\includegraphics[width=0.33\textwidth]{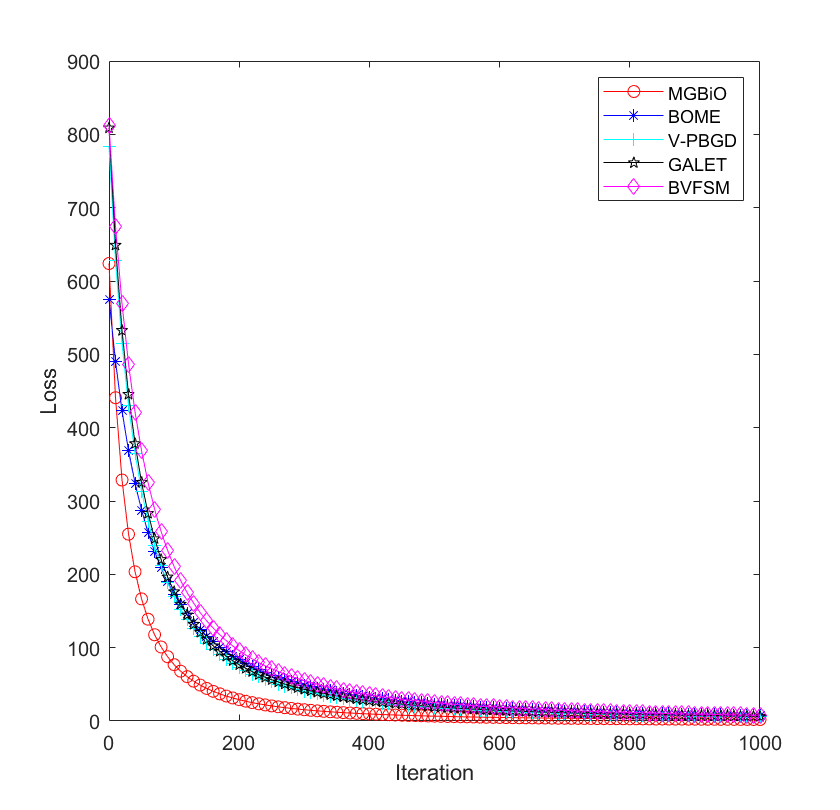}}
  \hfill
  \subfloat{\includegraphics[width=0.33\textwidth]{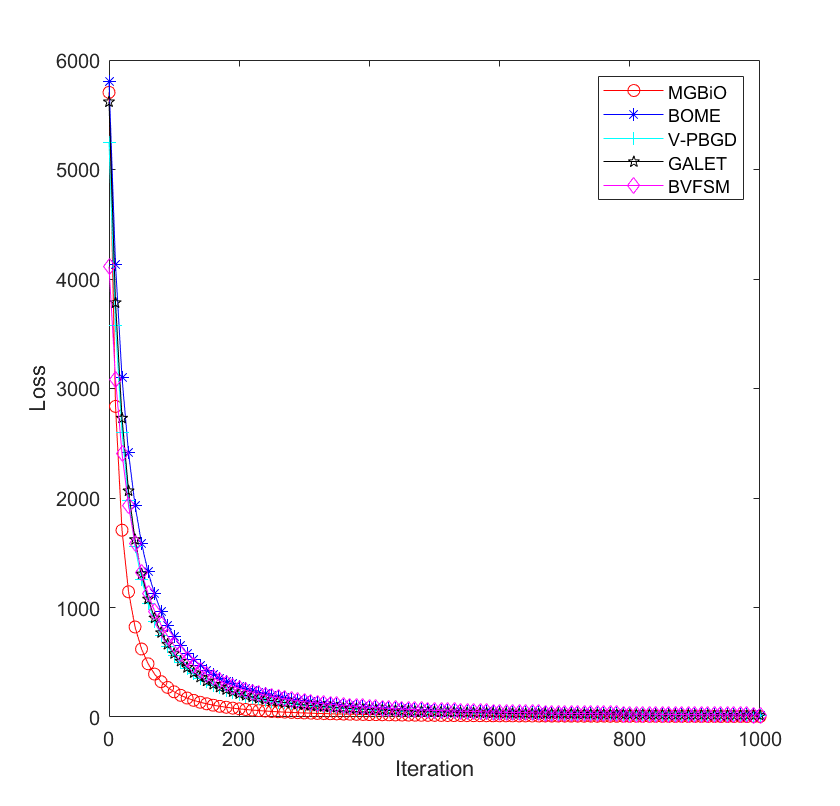}}
  \hfill
\caption{The loss of the \textbf{deterministic} algorithms under the case of $d=50$ (Left), $d=100$ (Middle) and $d=200$ (Right).}
\label{fig:3}
\end{figure}

\subsection{ Bilevel Polyak-{\L}ojasiewicz Game}
In the subsection, we consider the following bilevel Polyak-{\L}ojasiewicz Game:
\begin{align}
 & \min_{x\in \mathbb{R}^d} \ \frac{1}{2}x^TPx + x^TR^1y,   \\
 & \mbox{s.t.} \ \min_{y\in \mathbb{R}^d} \frac{1}{2}y^TQy + x^TR^2y, \nonumber
\end{align}
where $P=\frac{1}{n}\sum_{i=1}^n p_i(p_i)^T$, $Q=\frac{1}{n}\sum_{i=1}^n q_i(q_i)^T$,
$R^1= \frac{1}{n}\sum_{i=1}^n 0.01r^1_i(r^1_i)^T$ and $R^2= \frac{1}{n}\sum_{i=1}^n 0.01r^2_i(r^2_i)^T$.
Here the samples $\{p_i\}_{i=1}^n$,
$\{q_i\}_{i=1}^n$, $\{r^1_i\}_{i=1}^n$ and $\{r^2_i\}_{i=1}^n$ are independently drawn from Gaussian distributions $\mathcal{N}(0,\Sigma_{P})$, $\mathcal{N}(0,\Sigma_{Q})$, $\mathcal{N}(0,\Sigma_{R^1})$ and $\mathcal{N}(0,\Sigma_{R^2})$, respectively.
Meanwhile, we set $\Sigma_{P}=U^1D^1(U^1)^T$, where $U^1\in \mathbb{R}^{d\times l} \ (l<d)$
is column orthogonal, and $D^1\in \mathbb{R}^{l\times l}$ is diagonal and its diagonal elements are distributed uniformly in the interval $[\mu,L]$ with $0<\mu<L$. Let $\Sigma_{Q}=U^2D^2(U^2)^T$, where $U^2\in \mathbb{R}^{d\times l}$
is column orthogonal, and $D^2\in \mathbb{R}^{l\times l}$ is diagonal and its diagonal elements are distributed uniformly in the interval $[\mu,L]$ with $0<\mu<L$.
We also set $\Sigma_{R^1}=0.001V^1(V^1)^T$ and $\Sigma_{R^2}=0.001V^2(V^2)^T$, where each element of $V^1, V^2\in \mathbb{R}^{d\times d}$ is independently sampled from normal distribution $\mathcal{N}(0,1)$. Since
the covariance matrices $\Sigma_P$ and $\Sigma_Q$ are rank-deficient,
it is ensured that both $P$ and $Q$ are singular.
Hence the lower-level and upper-level objective functions may be not convex, while they satisfy the PL condition. In this experiment, we set $d=50$, $l=48$ and $n=2500$. For fair comparison, we set the basic learning rate as
$0.01$ for all algorithms. In the stochastic algorithms, we set mini-batch size as 10. From Figure~\ref{fig:1}, our deterministic MGBiO algorithm outperforms the baselines in Table~\ref{tab:1}. Meanwhile, our stochastic VR-MSGBiO algorithm outperforms the other algorithms.

\begin{figure}[ht]
\centering
 \subfloat{\includegraphics[width=0.33\textwidth]{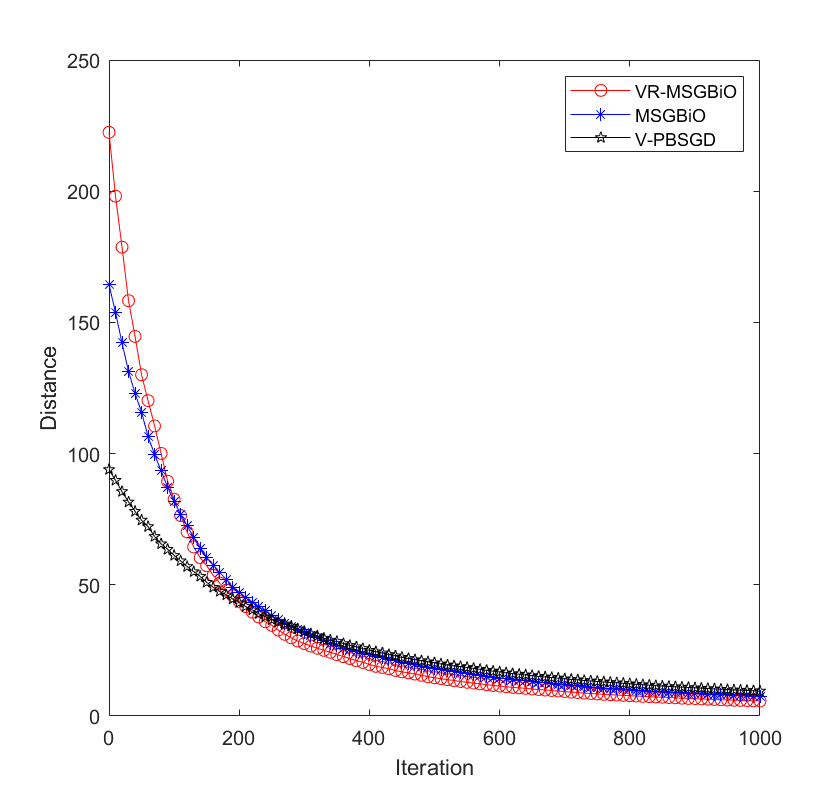}}
  \hfill
 \subfloat{\includegraphics[width=0.33\textwidth]{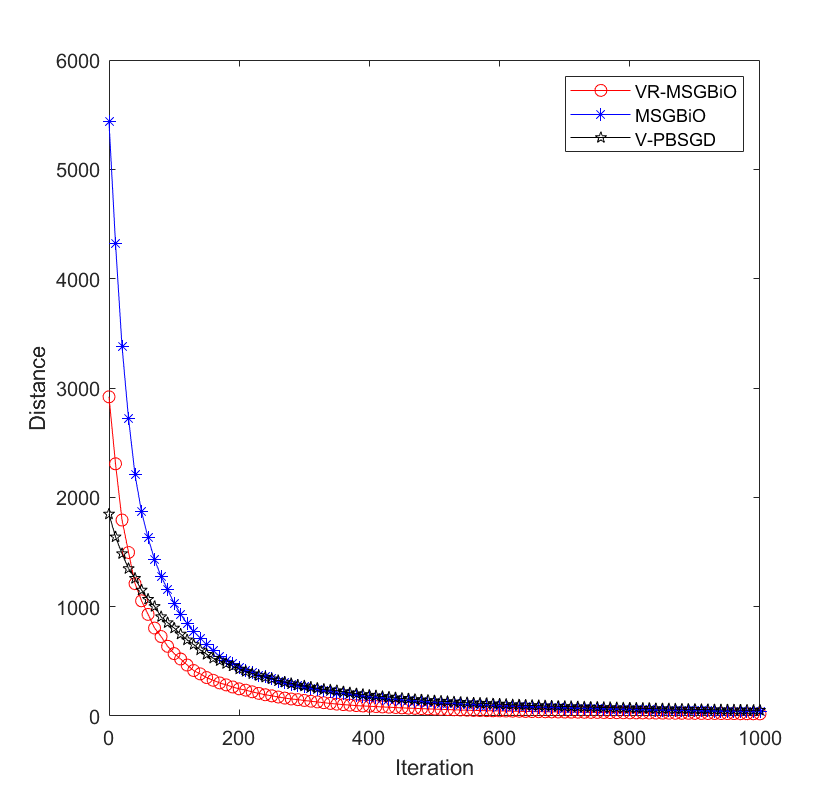}}
  \hfill
  \subfloat{\includegraphics[width=0.33\textwidth]{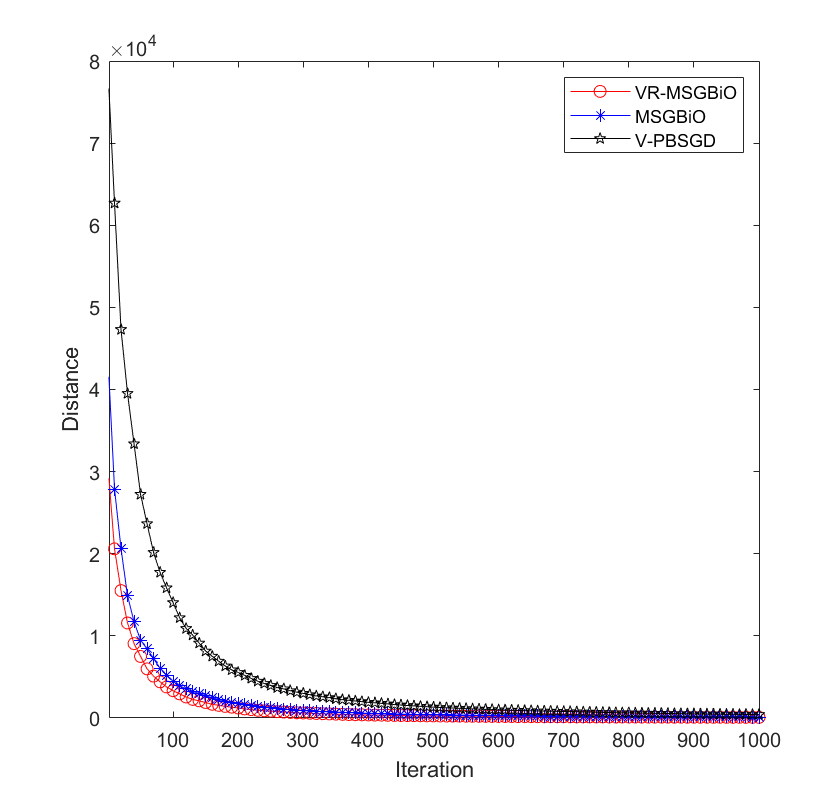}}
   \hfill
\caption{The ratio of distance $\|UU^T-H^*\|^2_F/\|H^*\|^2_F$ of the \textbf{stochastic} algorithms under the case of $d=50$ (Left), $d=100$ (Middle) and $d=200$ (Right).}
\label{fig:4}
\end{figure}

\begin{figure}[ht]
\centering
 \subfloat{\includegraphics[width=0.33\textwidth]{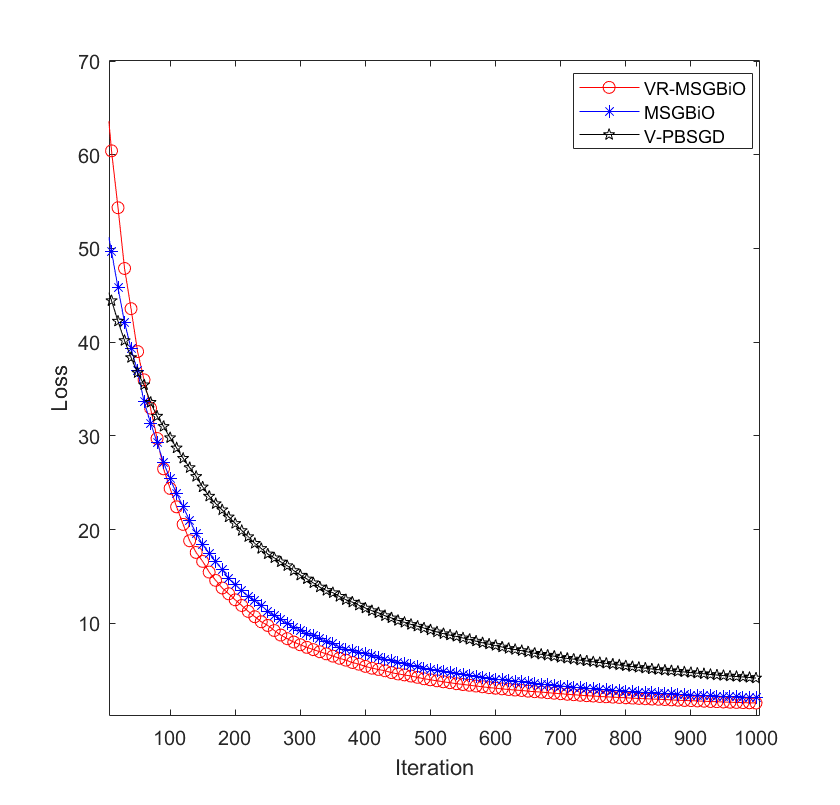}}
  \hfill
 \subfloat{\includegraphics[width=0.33\textwidth]{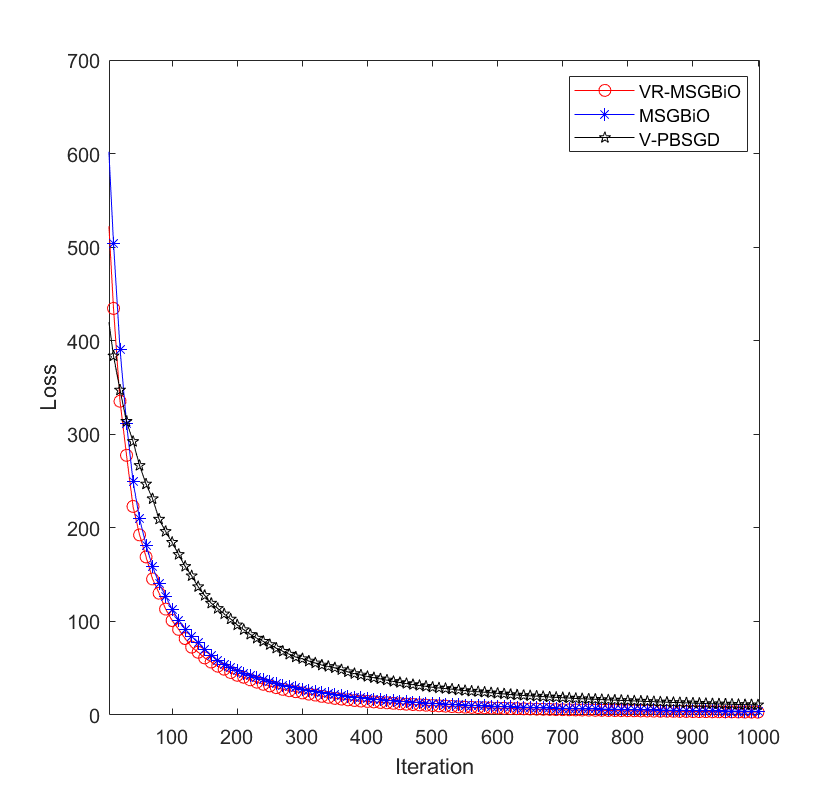}}
  \hfill
  \subfloat{\includegraphics[width=0.33\textwidth]{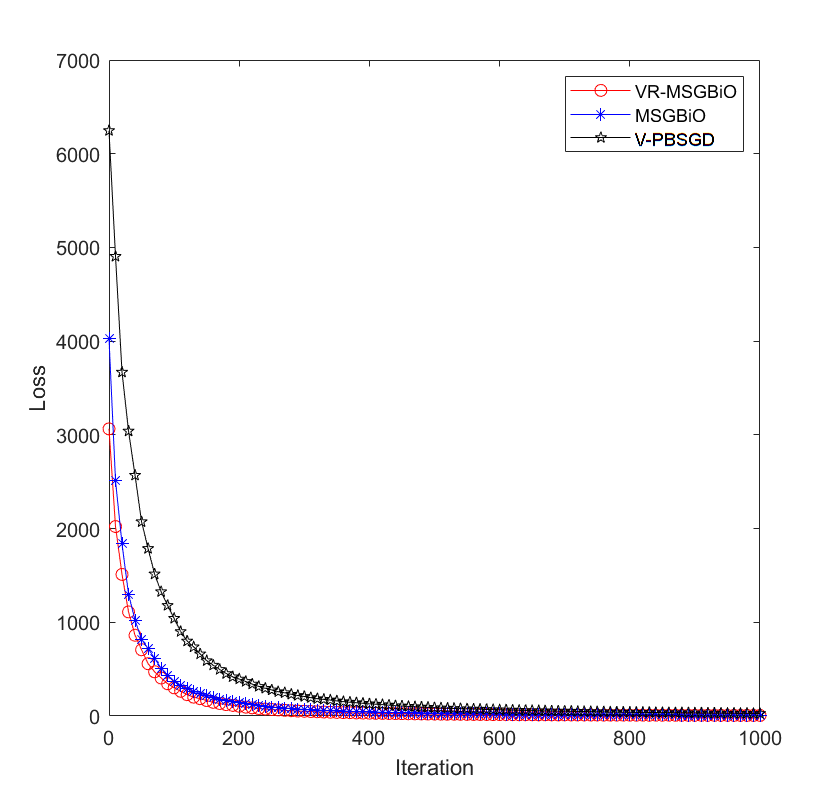}}
  \hfill
\caption{The loss of the \textbf{stochastic} algorithms under the case of $d=50$ (Left), $d=100$ (Middle) and $d=200$ (Right).}
\label{fig:5}
\end{figure}

\subsection{ Hyper-Representation Learning}
In the subsection, we conduct the hyper-representation learning task to verify the efficiency of our methods.  Learning hyper-representation is one of key points of meta learning, which extract better feature representations to be applied to many different tasks. Here we specifically consider the hyper-representation learning in matrix sensing task. Given $n$ sensing matrices $\{C_i\in \mathbb{R}^{d\times d}\}_{i=1}^n$ with $n$ observations $e_i=\langle C_i, H^*\rangle=\mbox{trace}(C_i^TH^*)$, where $H^*=U^*(U^*)^T$ is a low-rank symmetric matrix with
$U^*\in \mathbb{R}^{d\times r}$. The goal of matrix sensing task is to find the matrix $U^*$, which can be
represented the following problem:
\begin{align}
 \min_{U\in \mathbb{R}^{d\times r}} \frac{1}{n}\sum_{i=1}^n \ell_i(U)= \frac{1}{2}\big(\langle C_i, UU^T\rangle-e_i\big)^2.
\end{align}
Then we consider the hyper-representation learning in matrix sensing task, which be rewritten the following
bilevel optimization problem:
\begin{align}
& \min_{x\in \mathbb{R}^{d\times r-1}}\frac{1}{|D_v|}\sum_{i\in D_v} \ell_i(x,y^*(x)),  \\
 & \mbox{s.t.} \  y^*(x)\in \arg\min_{y\in \mathbb{R}^{d}}\frac{1}{|D_t|}\sum_{i\in D_t} \ell_i(x,y), \nonumber
\end{align}
where $U=[y;x]\in \mathbb{R}^{d\times r}$ is a concatenation of $x$ and $y$. Here we define variable $x$ to be the first $r-1$ columns of $U$ and variable $y$ to be the last column. Meanwhile, $D_t$ denotes the training dataset, and $D_v$ denotes the validation dataset.
The ground truth low-rank matrix $H^*$ is generated by $H^*=U^*(U^*)^T$, where each entry of $U^*$ is drawn from
normal distribution $\mathbb{N}(0,\frac{1}{d})$ independently. We randomly generate $n=20d$ samples of sensing
matrices $\{C_i\}_{i=1}^n$ from standard normal distribution, and then compute the corresponding no-noise
labels $e_i = \langle C_i, H^*\rangle$. We split all samples into two dataset: a train dataset
$D_t$ with 40\% data and a validation dataset $D_v$ with 60\% data.

In the experiment, for fair comparison, we set the basic learning rate as
$0.008$ for all algorithms. In the stochastic algorithms, we set mini-batch size as 10.
Let $\ell(U) = \frac{1}{2n}\sum_{i=1}^n\big(\langle C_i, UU^T\rangle-e_i\big)^2$ denote the loss.
Figures~\ref{fig:2}-\ref{fig:3} show that our deterministic MGBiO algorithm outperforms the baselines in Table~\ref{tab:1}. Figures~\ref{fig:4}-\ref{fig:5} show that our stochastic VR-MSGBiO algorithm has a better performance than the other algorithms.


\section{Conclusions }
In the paper, we studied a class of non-convex bilevel optimization problems, where both upper-level and lower-level subproblems are nonconvex, and the lower-level subproblem satisfies PL condition. To solve these deterministic bilevel problems, we proposed an efficient momentum-based gradient bilevel (MGBiO) method, which reaches a lower gradient complexity of $O(\epsilon^{-2})$ in finding an $\epsilon$-stationary solution. Meanwhile, to solve these stochastic bilevel problems, we presented a class of efficient momentum-based stochastic gradient bilevel methods (i.e., MSGBiO, VR-MSGBiO). In particular, our VR-MSGBiO obtains an near-optimal gradient complexity $\tilde{O}(\epsilon^{-3})$ in finding an $\epsilon$-stationary solution. Experimental results on bilevel PL game and hyper-representation learning demonstrate the efficiency of our algorithms.

\small

\bibliographystyle{plainnat}

\bibliography{NPLBiO}

\newpage

\appendix

\section{Appendix}
In this section, we provide the detailed convergence analysis of our algorithms.
We first review some useful lemmas.

\begin{lemma} \label{lem:A1}
(\cite{nesterov2018lectures})
Assume that $f(x)$ is a differentiable convex function and $\mathcal{X}$ is a convex set.
 $x^* \in \mathcal{X}$ is the solution of the
constrained problem $\min_{x\in \mathcal{X}}f(x)$, if
\begin{align}
 \langle \nabla f(x^*), x-x^*\rangle \geq 0, \ \forall x\in \mathcal{X}.
\end{align}
\end{lemma}

\begin{lemma} \label{lem:A2}
(\cite{karimi2016linear})
 The function $f(x): \mathbb{R}^d\rightarrow \mathbb{R}$ is $L$-smooth and satisfies PL condition with constant $\mu$, then it also satisfies error bound (EB) condition with $\mu$, i.e., for all $x \in \mathbb{R}^d$
\begin{align}
 \|\nabla f(x)\| \geq \mu\|x^*-x\|,
\end{align}
where $x^* \in \arg\min_{x} f(x)$. It also satisfies quadratic growth (QG) condition with $\mu$, i.e.,
\begin{align}
 f(x)-\min_x f(x) \geq \frac{\mu}{2}\|x^*-x\|^2.
\end{align}
\end{lemma}

\subsection{Convergence Analysis of MGBiO}
In this subsection, we detail the convergence analysis of our MGBiO algorithm.
We give some useful lemmas.

\begin{lemma} \label{lem:B1}
(Restatement of Lemma 1)
Under the above Assumption \ref{ass:2}, we have, for any $x\in \mathbb{R}^d$,
\begin{align}
 \nabla F(x) =\nabla_x f(x,y^*(x)) - \nabla^2_{xy} g(x,y^*(x))\Big[\nabla^2_{yy}g(x,y^*(x))\Big]^{-1}\nabla_y f(x,y^*(x)). \nonumber
\end{align}
\end{lemma}

\begin{proof}
Since $F(x)=f(x,y^*(x))$, we have
\begin{align}
\nabla F(x)=\nabla_x f(x,y^*(x)) + \nabla y^*(x)\nabla_y f(x,y^*(x)).
\end{align}
Meanwhile, according to the optimal condition of the Lower-Level problem in Problem \eqref{eq:1}, we have
\begin{align} \label{eq:B1}
 \nabla_y g(x,y^*(x))=0,
\end{align}
then further differentiating the above equality \eqref{eq:B1} on the variable $x$, we can obtain
\begin{align}
 \nabla_x\nabla_y g(x,y^*(x)) + \nabla y^*(x)\nabla_y\nabla_y g(x,y^*(x))=0.
\end{align}
According to Assumption \ref{ass:1}, the matrix  $\nabla_y\nabla_y g(x,y^*(x))$ is reversible. Thus, we have
\begin{align}
\nabla y^*(x) = - \nabla_x\nabla_y g(x,y^*(x))\big(\nabla_y\nabla_y g(x,y^*(x))\big)^{-1}.
\end{align}
Then we have
\begin{align}
\nabla F(x) & =\nabla_x f(x,y^*(x)) - \nabla_x\nabla_y g(x,y^*(x)) \big(\nabla_y\nabla_y g(x,y^*(x))\big)^{-1}\nabla_y f(x,y^*(x)) \nonumber \\
& =\nabla_x f(x,y^*(x)) - \nabla^2_{xy} g(x,y^*(x))\Big[\nabla^2_{yy}g(x,y^*(x))\Big]^{-1}\nabla_y f(x,y^*(x)).
\end{align}

\end{proof}

\begin{lemma} \label{lem:C1}
(Restatement of Lemma 2)
Under the above Assumptions (\ref{ass:1}-\ref{ass:3}), the functions (or mappings) $F(x)=f(x,y^*(x))$, $G(x)=g(x,y^*(x))$ and $y^*(x)\in \arg\min_{y\in \mathbb{R}^p}g(x,y)$ satisfy, for all $x_1,x_2\in \mathbb{R}^d$,
\begin{align}
 & \|y^*(x_1)-y^*(x_2)\| \leq \kappa\|x_1-x_2\|, \quad \|\nabla y^*(x_1) - \nabla y^*(x_2)\| \leq L_y\|x_1-x_2\| \nonumber \\
 & \|\nabla F(x_1) - \nabla F(x_2)\|\leq L_F\|x_1-x_2\|, \quad \|\nabla G(x_1) - \nabla G(x_2)\|\leq L_G\|x_1-x_2\| \nonumber
\end{align}
where $\kappa=C_{gxy}/\mu$, $L_y=\big( \frac{C_{gxy}L_{gyy}}{\mu^2} +  \frac{L_{gxy}}{\mu} \big) (1+ \frac{C_{gxy}}{\mu})$,
$L_F=\Big(L_f + L_f\kappa + C_{fy}\big( \frac{C_{gxy}L_{gyy}}{\mu^2} +  \frac{L_{gxy}}{\mu} \big)\Big)(1+\kappa)$ and $L_G=\Big(L_g + L_g\kappa + C_{gy}\big( \frac{C_{gxy}L_{gyy}}{\mu^2} +  \frac{L_{gxy}}{\mu} \big)\Big)(1+\kappa)$.
\end{lemma}

\begin{proof}
From the above lemma \ref{lem:B1}, we have
\begin{align}
\nabla y^*(x) = - \nabla^2_{xy}g(x,y^*(x))\Big(\nabla^2_{yy} g(x,y^*(x))\Big)^{-1}.
\end{align}
According to the above Assumptions \ref{ass:2} and \ref{ass:3}, we have for any $x\in \mathbb{R}^d$,
\begin{align}
\|\nabla y^*(x)\| & = \|- \nabla^2_{xy} g(x,y^*(x))\Big(\nabla^2_{yy} g(x,y^*(x))\Big)^{-1}\| \nonumber \\
& \leq \| \nabla^2_{xy} g(x,y^*(x))\|\|\big(\nabla^2_{yy} g(x,y^*(x))\big)^{-1}\| \leq \frac{C_{gxy}}{\mu}.
\end{align}
By using the Mean Value Theorem, we have for any $x_1,x_2\in \mathbb{R}^d$
\begin{align} \label{eq:CB1}
\|y^*(x_1)-y^*(x_1)\| = \|\nabla y^*(z)(x_1-x_2)\| \leq \frac{C_{gxy}}{\mu}\|x_1-x_2\|,
\end{align}
where $z\in (x_1,x_2)$.

For any $x_1,x_2\in \mathbb{R}^d$, we have
\begin{align} \label{eq:CB2}
 & \|\nabla y^*(x_1) - \nabla y^*(x_2)\| \nonumber \\
 & = \big\|\nabla^2_{xy} g(x_1,y^*(x_1))\Big(\nabla^2_{yy} g(x_1,y^*(x_1))\Big)^{-1}  -  \nabla^2_{xy} g(x_2,y^*(x_2))\Big(\nabla^2_{yy} g(x_2,y^*(x_2))\Big)^{-1} \big\| \nonumber \\
 & \leq \|\nabla^2_{xy} g(x_1,y^*(x_1))\Big(\nabla^2_{yy} g(x_1,y^*(x_1))\Big)^{-1} - \nabla^2_{xy} g(x_1,y^*(x_1))\Big(\nabla^2_{yy} g(x_2,y^*(x_2))\Big)^{-1}\| \nonumber \\
 & \quad + \|\nabla^2_{xy} g(x_1,y^*(x_1))\Big(\nabla^2_{yy} g(x_2,y^*(x_2))\Big)^{-1} -  \nabla^2_{xy} g(x_2,y^*(x_2))\Big(\nabla^2_{yy} g(x_2,y^*(x_2))\Big)^{-1}\| \nonumber \\
 & \leq C_{gxy}\sum_{m=1}\|\Big(\nabla^2_{yy} g(x_1,y^*(x_1))\Big)^{-1} - \Big(\nabla^2_{yy} g(x_2,y^*(x_2))\Big)^{-1}\| \nonumber \\
 & \quad + \frac{1}{\mu}\sum_{m=1}\|\nabla^2_{xy} g(x_1,y^*(x_1))-  \nabla^2_{xy} g(x_2,y^*(x_2))\| \nonumber \\
 & \leq C_{gxy}\sum_{m=1}\|\Big(\nabla^2_{yy} g(x_2,y^*(x_2))\Big)^{-1}\Big(\nabla^2_{yy} g(x_2,y^*(x_2)) - \nabla^2_{yy} g(x_1,y^*(x_1))\Big) \nonumber \\
 & \quad \cdot \Big(\nabla^2_{yy} g(x_1,y^*(x_1))\Big)^{-1}\|
 + \frac{1}{\mu}\sum_{m=1}\Big( L_{gxy}\|x_1-x_2\| + L_{gxy}\|y^*(x_1)-y^*(x_2)\|\Big) \nonumber \\
 & \leq \frac{C_{gxy}}{\mu^2}\sum_{m=1}\|\big(\nabla^2_{yy} g(x_2,y^*(x_2)) -
 \nabla^2_{yy} g(x_1,y^*(x_1))\big)\|
 + \frac{1}{\mu}\sum_{m=1}\Big( L_{gxy}\|x_1-x_2\| + L_{gxy}\|y^*(x_1)-y^*(x_2)\|\Big) \nonumber \\
 & \leq \frac{C_{gxy}}{\mu^2}\sum_{m=1}\big(L_{gyy}\|x_1-x_2\| + L_{gyy}\|y^*(x_1)-y^*(x_2)\|\big) + \frac{1}{\mu}\sum_{m=1}\Big( L_{gxy}\|x_1-x_2\| + L_{gxy}\|y^*(x_1)-y^*(x_2)\|\Big) \nonumber \\
 & \leq \big( \frac{C_{gxy}L_{gyy}}{\mu^2} +  \frac{L_{gxy}}{\mu} \big) (1+ \frac{C_{gxy}}{\mu})\|x_1-x_2\|,
\end{align}
where the last inequality holds by the above inequality (\ref{eq:CB1}).

From the above lemma \ref{lem:1}, we have
\begin{align}
 \nabla F(x) = \nabla_x f(x,y^*(x)) + \nabla y^*(x)^T\nabla_y f(x,y^*(x)),
\end{align}
then for any $x_1,x_2\in \mathbb{R}^d$, we have
\begin{align}
 & \|\nabla F(x_1) -\nabla F(x_2)\| \nonumber \\
 & = \|\nabla_x f(x_1,y^*(x_1)) + \nabla y^*(x_1)^T\Big(\nabla_y f(x_1,y^*(x_1))\Big) \nonumber \\
 & \qquad - \nabla_x f(x_1,y^*(x_2)) - \nabla y^*(x_2)^T\Big(\nabla_y f(x_2,y^*(x_2))\Big) \| \nonumber \\
 & \leq L_f\big(\|x_1-x_2\|+\|y^*(x_1)-y^*(x_2)\|\big) + \|\nabla y^*(x_1)^T\Big(\nabla_y f(x_1,y^*(x_1))\Big) \nonumber \\
 & \quad - \nabla y^*(x_2)^T\Big(\nabla_y f(x_1,y^*(x_1))\Big)  + \nabla y^*(x_2)^T\Big(\nabla_y f(x_1,y^*(x_1))\Big)  \nonumber \\
 & \quad - \nabla y^*(x_2)^T\Big(\nabla_y f(x_2,y^*(x_2))\Big)\| \nonumber \\
 & \leq L_f\big(\|x_1-x_2\|+\|y^*(x_1)-y^*(x_2)\|\big) + C_{fy}\|\nabla y^*(x_1)-\nabla y^*(x_2)\| \nonumber \\
 & \quad + \frac{C_{gxy}}{\mu}\|\nabla_y f(x_1,y^*(x_1)) - \nabla_y f(x_2,y^*(x_2))\| \nonumber \\
 & \leq L_f(1+\frac{C_{gxy}}{\mu})\|x_1-x_2\| + C_{fy}\|\nabla y^*(x_1)-\nabla y^*(x_2)\| + \frac{C_{gxy}L_f}{\mu}(1+\frac{C_{gxy}}{\mu})\|x_1-x_2\| \nonumber \\
 & \leq \Big(L_f + \frac{C_{gxy}L_f}{\mu} + C_{fy}\big( \frac{C_{gxy}L_{gyy}}{\mu^2} +  \frac{L_{gxy}}{\mu} \big)\Big)(1+\frac{C_{gxy}}{\mu})\|x_1-x_2\|,
\end{align}
where the last inequality holds by (\ref{eq:CB2}).

Since $\nabla G(x) =  \nabla_x g(x,y^*(x)) + \nabla y^*(x)^T\nabla_y g(x,y^*(x))$, we have for any $x_1,x_2\in \mathbb{R}^d$
\begin{align}
& \|\nabla G(x_1)-\nabla G(x_2)\| \nonumber \\
& = \|\nabla_x g(x_1,y^*(x_1)) + \nabla y^*(x_1)^T\Big(\nabla_y g(x_1,y^*(x_1))\Big) \nonumber \\
 & \qquad - \nabla_x g(x_1,y^*(x_2)) - \nabla y^*(x_2)^T\Big(\nabla_y g(x_2,y^*(x_2))\Big) \| \nonumber \\
 & \leq L_g\big(\|x_1-x_2\|+\|y^*(x_1)-y^*(x_2)\|\big) + \|\nabla y^*(x_1)^T\Big(\nabla_y g(x_1,y^*(x_1))\Big) \nonumber \\
 & \quad - \nabla y^*(x_2)^T\Big(\nabla_y g(x_1,y^*(x_1))\Big)  + \nabla y^*(x_2)^T\Big(\nabla_y g(x_1,y^*(x_1))\Big)  \nonumber \\
 & \quad - \nabla y^*(x_2)^T\Big(\nabla_y g(x_2,y^*(x_2))\Big)\| \nonumber \\
 & \leq L_g\big(\|x_1-x_2\|+\|y^*(x_1)-y^*(x_2)\|\big) + C_{gy}\|\nabla y^*(x_1)-\nabla y^*(x_2)\| \nonumber \\
 & \quad + \frac{C_{gxy}}{\mu}\|\nabla_y g(x_1,y^*(x_1)) - \nabla_y g(x_2,y^*(x_2))\| \nonumber \\
 & \leq L_g(1+\frac{C_{gxy}}{\mu})\|x_1-x_2\| + C_{gy}\|\nabla y^*(x_1)-\nabla y^*(x_2)\| + \frac{C_{gxy}L_g}{\mu}(1+\frac{C_{gxy}}{\mu})\|x_1-x_2\| \nonumber \\
 & \leq \Big(L_g + \frac{C_{gxy}L_g}{\mu} + C_{gy}\big( \frac{C_{gxy}L_{gyy}}{\mu^2} +  \frac{L_{gxy}}{\mu} \big)\Big)(1+\frac{C_{gxy}}{\mu})\|x_1-x_2\|.
\end{align}

\end{proof}

\begin{lemma} \label{lem:D1}
(Restatement of Lemma 3)
Let $\hat{\nabla} f(x,y)=\nabla_xf(x,y) -\hat{\Pi}_{C_{gxy}}\big[\nabla^2_{xy}g(x,y)\big] \big(\mathcal{S}_{[\mu,L_g]}\big[\nabla^2_{yy}g(x,y)\big]\big)^{-1}\Pi_{C_{fy}}\big[\nabla_yf(x,y)\big]$
and $\nabla F(x) = \nabla f(x,y^*(x))$,
 we have
 \begin{align}
 \|\hat{\nabla} f(x,y)-\nabla F(x)\|^2 \leq \hat{L}^2\|y^*(x)-y\|^2\leq \frac{2\hat{L}^2}{\mu}\big(g(x,y)-\min_y g(x,y)\big),
\end{align}
where $\hat{L}^2 = 4\big(L^2_f+ \frac{L^2_{gxy}C^2_{fy}}{\mu^2} + \frac{L^2_{gyy} C^2_{gxy}C^2_{fy}}{\mu^4} +
 \frac{L^2_fC^2_{gxy}}{\mu^2}\big)$.
\end{lemma}

\begin{proof}
  Since $\nabla F(x)=\nabla f(x,y^*(x))$, we have
\begin{align}
 & \|\nabla f(x,y^*(x))-\hat{\nabla} f(x,y)\|^2 \nonumber \\
 & = \|\nabla_xf(x,y^*(x)) - \nabla^2_{xy}g(x,y^*(x)) \big(\nabla^2_{yy}g(x,y^*(x))\big)^{-1}\nabla_yf(x,y^*(x)) \nonumber \\
 & \quad  - \nabla_xf(x,y) + \hat{\Pi}_{C_{gxy}}\big[\nabla^2_{xy}g(x,y)\big] \big(\mathcal{S}_{[\mu,L_g]}\big[\nabla^2_{yy}g(x,y)\big]\big)^{-1}\Pi_{C_{fy}}\big[\nabla_yf(x,y)\big]\|^2 \nonumber \\
 & = \|\nabla_xf(x,y^*(x)) - \nabla_xf(x,y) - \nabla^2_{xy}g(x,y^*(x)) \big(\nabla^2_{yy}g(x,y^*(x))\big)^{-1}\nabla_yf(x,y^*(x)) \nonumber \\
 & \quad  + \hat{\Pi}_{C_{gxy}}\big[\nabla^2_{xy}g(x,y)\big]  \big(\nabla^2_{yy}g(x,y^*(x))\big)^{-1}\nabla_yf(x,y^*(x)) - \hat{\Pi}_{C_{gxy}}\big[\nabla^2_{xy}g(x,y)\big]  \big(\nabla^2_{yy}g(x,y^*(x))\big)^{-1}\nabla_yf(x,y^*(x)) \nonumber \\
 & \quad + \hat{\Pi}_{C_{gxy}}\big[\nabla^2_{xy}g(x,y)\big]  \big(\mathcal{S}_{[\mu,L_g]}\big[\nabla^2_{yy}g(x,y)\big]\big)^{-1}\nabla_yf(x,y^*(x)) -  \hat{\Pi}_{C_{gxy}}\big[\nabla^2_{xy}g(x,y)\big]  \big(\mathcal{S}_{[\mu,L_g]}\big[\nabla^2_{yy}g(x,y)\big]\big)^{-1}\nabla_yf(x,y^*(x)) \nonumber \\
 & \quad + \hat{\Pi}_{C_{gxy}}\big[\nabla^2_{xy}g(x,y)\big]  \big(\mathcal{S}_{[\mu,L_g]}\big[\nabla^2_{yy}g(x,y)\big]\big)^{-1}\Pi_{C_{fy}}\big[\nabla_yf(x,y)\big]\|^2 \nonumber \\
 & \leq 4\|\nabla_xf(x,y^*(x)) - \nabla_xf(x,y)\|^2 + \frac{4C^2_{fy}}{\mu^2}\|\nabla^2_{xy}g(x,y^*(x))-\hat{\Pi}_{C_{gxy}}\big[\nabla^2_{xy}g(x,y)\big]\|^2
 \nonumber \\
 & \quad + \frac{4C^2_{gxy}C^2_{fy}}{\mu^4}\|\nabla^2_{yy}g(x,y^*(x))-\mathcal{S}_{[\mu,L_g]}\big[\nabla^2_{yy}g(x,y)\big]\|^2 +
 \frac{4C^2_{gxy}}{\mu^2}\|\nabla_yf(x,y^*(x))-\Pi_{C_{fy}}\big[\nabla_yf(x,y)\big]\|^2 \nonumber \\
 & \mathop{=}^{(i)} 4\|\nabla_xf(x,y^*(x)) - \nabla_xf(x,y)\|^2 + \frac{4C^2_{fy}}{\mu^2}\|\hat{\Pi}_{C_{gxy}}\big[\nabla^2_{xy}g(x,y^*(x))\big]- \hat{\Pi}_{C_{gxy}}\big[\nabla^2_{xy}g(x,y)\big]\|^2
 \nonumber \\
 & \quad + \frac{4C^2_{gxy}C^2_{fy}}{\mu^4}\|\mathcal{S}_{[\mu,L_g]}\big[\nabla^2_{yy}g(x,y^*(x))\big]-\mathcal{S}_{[\mu,L_g]}\big[\nabla^2_{yy}g(x,y)\big]\|^2 +
 \frac{4C^2_{gxy}}{\mu^2}\|\Pi_{C_{fy}}\big[\nabla_yf(x,y^*(x))\big]-\Pi_{C_{fy}}\big[\nabla_yf(x,y)\big]\|^2 \nonumber \\
 & \leq 4\|\nabla_xf(x,y^*(x)) - \nabla_xf(x,y)\|^2 + \frac{4C^2_{fy}}{\mu^2}\|\nabla^2_{xy}g(x,y^*(x))-\nabla^2_{xy}g(x,y)\|^2
 \nonumber \\
 & \quad + \frac{4C^2_{gxy}C^2_{fy}}{\mu^4}\|\nabla^2_{yy}g(x,y^*(x))-\nabla^2_{yy}g(x,y)\|^2 +
 \frac{4C^2_{gxy}}{\mu^2}\|\nabla_yf(x,y^*(x))-\nabla_yf(x,y)\|^2 \nonumber \\
 & \leq 4\big(L^2_f+ \frac{L^2_{gxy}C^2_{fy}}{\mu^2} + \frac{L^2_{gyy} C^2_{gxy}C^2_{fy}}{\mu^4} +
 \frac{L^2_fC^2_{gxy}}{\mu^2}\big)\|y^*(x)-y\|^2 \nonumber \\
 & = \hat{L}^2\|y^*(x)-y\|^2 \leq \frac{2\hat{L}^2}{\mu}\big(g(x,y)-\min_y g(x,y)\big),
\end{align}
where the above equality (i) due to Assumption~\ref{ass:2}, i.e., $\mathcal{S}_{[\mu,L_g]}\big[\nabla^2_{yy}g(x,y^*(x))\big]=\nabla^2_{yy}g(x,y^*(x))$ and Assumption~\ref{ass:3}, i.e., $\hat{\Pi}_{C_{gxy}}\big[\nabla^2_{xy}g(x,y^*(x))\big]=\nabla^2_{xy}g(x,y^*(x))$, $\Pi_{C_{fy}}\big[\nabla_yf(x,y^*(x))\big]=\nabla_yf(x,y^*(x))$, and the second last inequality is due to Assumptions~\ref{ass:2}-\ref{ass:4};
the last inequality holds by Lemma~\ref{lem:A2}.

\end{proof}

\begin{lemma} \label{lem:E1}
(Restatement of Lemma 4)
Suppose the sequence $\{x_t,y_t\}_{t=1}^T$ be generated from Algorithm \ref{alg:1} or \ref{alg:2} or \ref{alg:3}.
Under the above Assumptions~\ref{ass:1}-\ref{ass:2}, given $\gamma\leq \frac{\lambda\mu}{8L_G}$ and $0<\lambda \leq \frac{1}{2L_g\eta_t}$ for all $t\geq 1$, we have
\begin{align}
g(x_{t+1},y_{t+1}) - G(x_{t+1})
& \leq (1-\frac{\eta_t\lambda\mu}{2}) \big(g(x_t,y_t) -G(x_t)\big) + \frac{\eta_t}{8\gamma}\|\tilde{x}_{t+1}-x_t\|^2  -\frac{\eta_t}{4\lambda}\|\tilde{y}_{t+1}-y_t\|^2 \nonumber \\
& \quad + \eta_t\lambda\|\nabla_y g(x_t,y_t)-v_t\|^2,
\end{align}
where $G(x_t)=g(x_t,y^*(x_t))$ with $y^*(x_t) \in \arg\min_{y}g(x_t,y)$ for all $t\geq 1$.
\end{lemma}

\begin{proof}
By using Assumption~\ref{lem:3}, i.e., $L_g$-smoothness of $g(x,\cdot)$, such that
\begin{align}
    g(x_{t+1},y_{t+1}) \leq g(x_{t+1},y_t) + \langle \nabla_y g(x_{t+1},y_t), y_{t+1}-y_t \rangle + \frac{L_g}{2}\|y_{t+1}-y_t\|^2 ,
\end{align}
then we have
\begin{align} \label{eq:E1}
    & g(x_{t+1},y_{t+1}) \leq g(x_{t+1},y_t) + \eta_t\langle \nabla_y g(x_{t+1},y_t), \tilde{y}_{t+1}-y_t \rangle + \frac{L_g\eta^2_t}{2}\|\tilde{y}_{t+1}-y_t\|^2.
\end{align}

Next, we bound the inner product in \eqref{eq:E1},
\begin{align} \label{eq:E2}
    & \eta_t\langle \nabla_y g(x_{t+1},y_t), \tilde{y}_{t+1}-y_t \rangle
     = -\eta_t\lambda\langle \nabla_y g(x_{t+1},y_t), v_t \rangle \nonumber \\
    & = -\frac{\eta_t\lambda}{2}\Big( \|\nabla_y g(x_{t+1},y_t)\|^2 + \|v_t\|^2 - \|\nabla_y g(x_{t+1},y_t)-\nabla_y g(x_t,y_t) + \nabla_y g(x_t,y_t)-v_t\|^2 \Big) \nonumber \\
    & \leq -\frac{\eta_t\lambda}{2} \|\nabla_y g(x_{t+1},y_t)\|^2 -\frac{\eta_t}{2\lambda} \|\tilde{y}_{t+1}-y_t\|^2 + \eta_t\lambda L^2_g \|x_{t+1}-x_t\|^2 + \eta_t\lambda\|\nabla_y g(x_t,y_t)-v_t\|^2 \nonumber \\
    & \leq -\eta_t\lambda\mu\big(g(x_{t+1},y_t)-G(x_{t+1})\big)-\frac{\eta_t}{2\lambda} \|\tilde{y}_{t+1}-y_t\|^2 + \eta_t\lambda L^2_g \|x_{t+1}-x_t\|^2 + \eta_t\lambda\|\nabla_y g(x_t,y_t)-v_t\|^2,
\end{align}
where the last inequality is due to the quadratic growth condition of $\mu$-PL functions, i.e.,
\begin{align}
    \|\nabla_y g(x_{t+1},y_t)\|^2 \geq 2\mu\big( g(x_{t+1},y_t)-\min_{y'}g(x_{t+1},y')\big) = 2\mu\big( g(x_{t+1},y_t) - G(x_{t+1})\big).
\end{align}
Substituting \eqref{eq:E2} in \eqref{eq:E1}, we have
\begin{align} \label{eq:E3}
    g(x_{t+1},y_{t+1})
    & \leq g(x_{t+1},y_t)-\eta_t\lambda\mu\big(g(x_{t+1},y_t)- G(x_{t+1})\big)-\frac{\eta_t}{2\lambda} \|\tilde{y}_{t+1}-y_t\|^2 + \eta_t\lambda L^2_g \|x_{t+1}-x_t\|^2 \nonumber \\
    & \quad + \eta_t\lambda\|\nabla_y g(x_t,y_t)-v_t\|^2 + \frac{L_g\eta^2_t}{2}\|\tilde{y}_{t+1}-y_t\|^2,
\end{align}
then rearranging the terms, we can obtain
\begin{align} \label{eq:E4}
    g(x_{t+1},y_{t+1})- G(x_{t+1})
    & \leq (1-\eta_t\lambda\mu)\big(g(x_{t+1},y_t)- G(x_{t+1})\big)-\frac{\eta_t}{2\lambda} \|\tilde{y}_{t+1}-y_t\|^2 + \eta_t\lambda L^2_g \|x_{t+1}-x_t\|^2 \nonumber \\
    & \quad + \eta_t\lambda\|\nabla_y g(x_t,y_t)-v_t\|^2 + \frac{L_g\eta^2_t}{2}\|\tilde{y}_{t+1}-y_t\|^2.
\end{align}

Next, using $L_g$-smoothness of function $f(\cdot,y_t)$, such that
\begin{align}
     g(x_{t+1},y_t) \leq g(x_t,y_t) + \langle \nabla_x g(x_t,y_t), x_{t+1}-x_t \rangle +  \frac{L_g}{2}\|x_{t+1}-x_t\|^2 ,
\end{align}
then we have
\begin{align}
    &g(x_{t+1},y_t) - g(x_t,y_t) \nonumber \\
    & \leq \langle \nabla_x g(x_t,y_t), x_{t+1}-x_t \rangle + \frac{L_g}{2}\|x_{t+1}-x_t\|^2 \nonumber \\
    & = \eta_t\langle \nabla_x g(x_t,y_t) - \nabla G(x_t), \tilde{x}_{t+1}-x_t \rangle + \eta_t\langle \nabla G(x_t), \tilde{x}_{t+1}-x_t \rangle + \frac{L_g\eta^2_t}{2}\|\tilde{x}_{t+1}-x_t\|^2 \nonumber \\
    & \leq \frac{\eta_t}{8\gamma}\|\tilde{x}_{t+1}-x_t\|^2 + 2\eta_t\gamma\|\nabla_x g(x_t,y_t) - \nabla G(x_t)\|^2  + \eta_t\langle \nabla G(x_t), \tilde{x}_{t+1}-x_t \rangle + \frac{L_g\eta^2_t}{2}\|\tilde{x}_{t+1}-x_t\|^2 \nonumber \\
    & \leq \frac{\eta_t}{8\gamma}\|\tilde{x}_{t+1}-x_t\|^2 + 2L^2_g\eta_t\gamma \|y_t - y^*(x_t)\|^2 + G(x_{t+1}) - G(x_t) \nonumber \\
    & \quad + \frac{\eta^2_tL_G}{2}\|\tilde{x}_{t+1}-x_t\|^2 + \frac{\eta^2_tL_g}{2}\|\tilde{x}_{t+1}-x_t\|^2 \nonumber \\
    & \leq \frac{4L^2_g\eta_t\gamma}{\mu} \big(g(x_t,y_t) - G(x_t)\big) + G(x_{t+1})- G(x_t)+ (\frac{\eta_t }{8\gamma}+\eta^2_tL_G)\|\tilde{x}_{t+1}-x_t\|^2,
\end{align}
where the second last inequality is due to
$L_G$-smoothness of function $G(x)$, and the last inequality holds by Lemma~\ref{lem:A2} and $L_g\leq L_G$.
Then we have
\begin{align} \label{eq:E5}
    g(x_{t+1},y_t) -G(x_{t+1}) & = g(x_{t+1},y_t)- g(x_t,y_t) + g(x_t,y_t)- G(x_t) + G(x_t) -G(x_{t+1})\nonumber \\
    & \leq (1+\frac{4L^2_g\eta_t\gamma}{\mu}) \big(g(x_t,y_t) -G(x_t)\big) + (\frac{\eta_t }{8\gamma}+\eta^2_tL_G)\|\tilde{x}_{t+1}-x_t\|^2.
\end{align}

Substituting \eqref{eq:E5} in \eqref{eq:E4}, we get
\begin{align}
    & g(x_{t+1},y_{t+1})- G(x_{t+1})\nonumber \\
    & \leq (1-\eta_t\lambda\mu)(1+\frac{4L^2_g\eta_t\gamma}{\mu}) \big(g(x_t,y_t) -G(x_t)\big) + \eta_t(\frac{1}{8\gamma}+\eta_tL_G)(1-\eta_t\lambda\mu)\|\tilde{x}_{t+1}-x_t\|^2 \nonumber \\
    & \quad -\frac{\eta_t}{2\lambda} \|\tilde{y}_{t+1}-y_t\|^2 + \eta_t\lambda L^2_g \|x_{t+1}-x_t\|^2 + \eta_t\lambda\|\nabla_y g(x_t,y_t)-v_t\|^2 + \frac{L_g\eta^2_t}{2}\|\tilde{y}_{t+1}-y_t\|^2 \nonumber \\
    & = (1-\eta_t\lambda\mu)(1+\frac{4L^2_g\eta_t\gamma}{\mu}) \big(g(x_t,y_t) -G(x_t)\big) + \eta_t\big(\frac{1}{8\gamma}+\eta_tL_G-\frac{\eta_t\lambda\mu}{8\gamma}-\eta^2_tL_G\lambda\mu+\eta^2_tL^2_g\lambda\big)\|\tilde{x}_{t+1}-x_t\|^2 \nonumber \\
    & \quad -\frac{\eta_t}{2}\big(\frac{1}{\lambda}-L_g\eta_t\big) \|\tilde{y}_{t+1}-y_t\|^2 + \eta_t\lambda\|\nabla_y g(x_t,y_t)-v_t\|^2 \nonumber \\
    & \leq (1-\frac{\eta_t\lambda\mu}{2}) \big(g(x_t,y_t) -G(x_t)\big) + \frac{\eta_t}{8\gamma}\|\tilde{x}_{t+1}-x_t\|^2  -\frac{\eta_t}{4\lambda}\|\tilde{y}_{t+1}-y_t\|^2 + \eta_t\lambda\|\nabla_y g(x_t,y_t)-v_t\|^2,
\end{align}
where the last inequality holds by $\gamma\leq \frac{\lambda\mu}{8L_G}$, $L_G\geq L_g(1+\kappa)^2$ and $\lambda\leq \frac{1}{2L_g\eta_t}$ for all $t\geq 1$, i.e.,
\begin{align}
   & \gamma\leq \frac{\lambda\mu}{8L_G} \Rightarrow \lambda \geq \frac{8L_G\gamma}{\mu} \geq   \frac{8L_g}{\mu}(1+\kappa)^2\gamma\geq 8\kappa^2\gamma \Rightarrow  \frac{\eta_t\lambda\mu}{2} \geq \frac{4L^2_f\eta_t\gamma}{\mu} \nonumber \\
   & L_G\geq L_g(1+\kappa)^2 \geq L_g(1+\kappa) \Rightarrow \eta^2_tL_G\lambda\mu \geq \eta^2_tL^2_f\lambda \nonumber \\
   &\lambda \leq \frac{1}{2\eta_tL_g}  \Rightarrow \frac{1}{2\lambda} \geq \eta_t L_g, \ \forall t\geq 1.
\end{align}

\end{proof}

\begin{theorem}  \label{th:A1}
(Restatement of Theorem 1)
 Under the above Assumptions (\ref{ass:1}-\ref{ass:5}), in the Algorithm \ref{alg:1}, let $\eta_t=\eta$ for all $t\geq 0$, $0< \gamma \leq \min\big(\frac{1}{2L_F\eta},\frac{\lambda\mu^2}{16\hat{L}^2} \big)$, $0<\lambda\leq \frac{1}{2L_g\eta}$. When $\mathcal{X}\subseteq \mathbb{R}^d$, we can get
 \begin{align}
 \frac{1}{T}\sum_{t=1}^T\|\mathcal{G}(x_t,\nabla F(x_t),\gamma)\|
 & \leq \frac{1}{T} \sum_{t=1}^T \big[ \|w_t-\nabla F(x_t)\| + \frac{1}{\gamma}\|\tilde{x}_{t+1}-x_t\|  \big] \leq \frac{4\sqrt{R}}{\sqrt{3T\gamma\eta}},
\end{align}
when $\mathcal{X}=\mathbb{R}^d$, we can get
 \begin{align}
 \frac{1}{T}\sum_{t=1}^T\|\nabla F(x_t)\|
  & \leq \frac{1}{T} \sum_{t=1}^T \big[ \|w_t-\nabla F(x_t)\| + \frac{1}{\gamma}\|\tilde{x}_{t+1}-x_t\|  \big] \leq \frac{4\sqrt{R}}{\sqrt{3T\gamma\eta}},
\end{align}
where $R= F(x_1) -F^* +g(x_1,y_1)-G(x_1)$.
\end{theorem}

\begin{proof}
According the line 6 of Algorithm~\ref{alg:1}, we have
\begin{align} \label{eq:AG1}
\tilde{x}_{t+1} = \arg\min_{x\in \mathcal{X}}
\Big\{ \langle w_t, x-x_t\rangle + \frac{1}{2\gamma}\|x-x_t\|^2 \Big\}.
\end{align}
By using the optimal condition of the above subproblem~(\ref{eq:AG1}), we have
\begin{align}
 \langle w_t + \frac{1}{\gamma}(\tilde{x}_{t+1}-x_t), x_t-\tilde{x}_{t+1}\rangle \geq 0,
\end{align}
and then we have
\begin{align} \label{eq:AG2}
 \langle w_t,\tilde{x}_{t+1} -x_t\rangle \leq -\frac{1}{\gamma}\|\tilde{x}_{t+1}-x_t\|^2.
\end{align}
According to Lemma~\ref{lem:C1}, i.e., the function $F(x)$ has $L_F$-Lipschitz continuous gradient,
we have
\begin{align} \label{eq:AG3}
  F(x_{t+1}) & \leq F(x_t) + \langle \nabla F(x_t), x_{t+1}-x_t\rangle + \frac{L_F}{2}\|x_{t+1}-x_t\|^2 \nonumber \\
  & = F(x_t) + \eta_t \langle \nabla F(x_t) -w_t +w_t, \tilde{x}_{t+1}-x_t\rangle + \frac{\eta^2_tL_F}{2}\|\tilde{x}_{t+1}-x_t\|^2 \nonumber \\
  & \leq F(x_t) + \eta_t \langle \nabla F(x_t) -w_t , \tilde{x}_{t+1}-x_t\rangle - \frac{\eta_t}{\gamma}\|\tilde{x}_{t+1}-x_t\|^2 + \frac{\eta^2_tL_F}{2}\|\tilde{x}_{t+1}-x_t\|^2 \nonumber \\
  & \leq F(x_t) + \eta_t\gamma\|\nabla F(x_t) -w_t\|^2 + \frac{\eta_t}{4\gamma}\|\tilde{x}_{t+1}-x_t\|^2 - \frac{\eta_t}{\gamma}\|\tilde{x}_{t+1}-x_t\|^2 + \frac{\eta^2_tL_F}{2}\|\tilde{x}_{t+1}-x_t\|^2 \nonumber \\
  & \leq F(x_t) + \eta_t\gamma\|\nabla F(x_t) -w_t\|^2 - \frac{\eta_t}{2\gamma}\|\tilde{x}_{t+1}-x_t\|^2,
\end{align}
where the second inequality holds by the above inequality~(\ref{eq:AG2}), and the last inequality holds by $\gamma\leq \frac{1}{2L_F\eta_t}$.

Since $w_t=\hat{\nabla}f(x_t,y_t)$ in Algorithm~\ref{alg:1}, we have
\begin{align}\label{eq:AG4}
 \|w_t - \nabla F(x_t)\|^2
  = \| \hat{\nabla} f(x_t,y_t) - \nabla F(x_t) \|^2 \leq \frac{2\hat{L}^2}{\mu}\big(g(x_t,y_t)-G(x_t)\big).
\end{align}
Plugging the above inequalities~(\ref{eq:AG4}) into (\ref{eq:AG3}), we have
\begin{align} \label{eq:AG5}
  F(x_{t+1})  \leq F(x_t) +\frac{2\eta_t\gamma\hat{L}^2}{\mu}\big(g(x_t,y_t)-G(x_t)\big) - \frac{\eta_t}{2\gamma}\|\tilde{x}_{t+1}-x_t\|^2,
\end{align}
According to Lemma~\ref{lem:E1}, we have
\begin{align} \label{eq:AG6}
& g(x_{t+1},y_{t+1})-G(x_{t+1})-\big(g(x_t,y_t)-G(x_t)\big) \nonumber \\
& \leq -\frac{\eta_t\lambda\mu}{2} \big(g(x_t,y_t) -G(x_t)\big) + \frac{\eta_t}{8\gamma}\|\tilde{x}_{t+1}-x_t\|^2   -\frac{\eta_t}{4\lambda}\|\tilde{y}_{t+1}-y_t\|^2+ \eta_t\lambda\|\nabla_y g(x_t,y_t)-v_t\|^2 \nonumber \\
& \leq  -\frac{\eta_t\lambda\mu}{2} \big(g(x_t,y_t) -G(x_t)\big) + \frac{\eta_t}{8\gamma}\|\tilde{x}_{t+1}-x_t\|^2,
\end{align}
where the last inequality is due to $\frac{\eta_t}{4\lambda}>0$ and $v_t=\nabla_y g(x_t,y_t)$.

Next, we define a useful Lyapunov function (i.e. potential function), for any $t\geq 1$
\begin{align}
 \Omega_t = F(x_t) + g(x_t,y_t)-G(x_t).
\end{align}

By using the above inequalities~(\ref{eq:AG5}) and~(\ref{eq:AG6}), we have
\begin{align}
 \Omega_{t+1} - \Omega_t & = F(x_{t+1}) -F(x_t) + g(x_{t+1},y_{t+1})-G(x_{t+1})-\big(g(x_t,y_t)-G(x_t)\big)
 \nonumber \\
 & \leq \frac{2\eta_t\gamma\hat{L}^2}{\mu}\big(g(x_t,y_t)-G(x_t)\big) - \frac{\eta_t}{2\gamma}\|\tilde{x}_{t+1}-x_t\|^2 -\frac{\eta_t\lambda\mu}{2} \big(g(x_t,y_t) -G(x_t)\big) + \frac{\eta_t}{8\gamma}\|\tilde{x}_{t+1}-x_t\|^2   \nonumber \\
& \leq - \frac{3\eta_t}{8\gamma}\|\tilde{x}_{t+1}-x_t\|^2 -\frac{3\eta_t\lambda\mu}{8} \big(g(x_t,y_t) -G(x_t)\big),
\end{align}
where the last inequality is due to $\gamma\leq \frac{\lambda\mu^2}{16\hat{L}^2}$ and $g(x_t,y_t) -G(x_t)\geq 0$.
Let $\eta_t=\eta$ for all $t\geq1$, then we can obtain
\begin{align} \label{eq:AG7}
\frac{1}{\gamma}\|\tilde{x}_{t+1}-x_t\|^2 + \lambda\mu \big(g(x_t,y_t) -G(x_t)\big) \leq \frac{8(\Omega_t  -\Omega_{t+1})}{3\eta}
\end{align}

Since $\gamma\leq \frac{\lambda\mu^2}{16\hat{L}^2}$ and $g(x_t,y_t) -G(x_t)\geq 0$, we have
\begin{align}
&\frac{1}{\gamma}\|\tilde{x}_{t+1}-x_t\|^2 +\frac{ 16\hat{L}^2\gamma}{\mu} \big(g(x_t,y_t) -G(x_t)\big)
\nonumber \\
&  \leq
\frac{1}{\gamma}\|\tilde{x}_{t+1}-x_t\|^2 + \lambda\mu \big(g(x_t,y_t) -G(x_t)\big) \leq \frac{8(\Omega_t  -\Omega_{t+1})}{3\eta},
\end{align}
then we can get
\begin{align}
&\frac{1}{\gamma^2}\|\tilde{x}_{t+1}-x_t\|^2 +\frac{ 16\hat{L}^2}{\mu} \big(g(x_t,y_t) -G(x_t)\big)  \leq \frac{8(\Omega_t  -\Omega_{t+1})}{3\gamma\eta}.
\end{align}

According to Lemma~\ref{lem:D1}, we can obtain
\begin{align} \label{eq:AG8}
\frac{1}{\gamma^2}\|\tilde{x}_{t+1}-x_t\|^2 + 8\|\nabla F(x_t) - \hat{\nabla}f(x_t,y_t)\|^2 & \leq  \frac{1}{\gamma^2}\|\tilde{x}_{t+1}-x_t\|^2 + 8\hat{L}^2\|y_t-y^*(x_t)\|^2\nonumber \\
& \leq \frac{1}{\gamma^2}\|\tilde{x}_{t+1}-x_t\|^2 +\frac{ 16\hat{L}^2}{\mu} \big(g(x_t,y_t) -G(x_t)\big)  \nonumber \\
& \leq \frac{8(\Omega_t  -\Omega_{t+1})}{3\gamma\eta}.
\end{align}

Summing the above inequality~(\ref{eq:AG8}) from $t=1$ to $T$, we have
\begin{align}
 \sum_{t=1}^T\big( \frac{1}{\gamma^2}\|\tilde{x}_{t+1}-x_t\|^2 + \|\nabla F(x_t) - \hat{\nabla}f(x_t,y_t)\|^2 \big) & \leq  \sum_{t=1}^T\big( \frac{1}{\gamma^2}\|\tilde{x}_{t+1}-x_t\|^2 + 8\|\nabla F(x_t) - \hat{\nabla}f(x_t,y_t)\|^2 \big) \nonumber \\
& \leq \frac{8(\Omega_1 -\Omega_T)}{3\gamma\eta},
\end{align}
then we can get
\begin{align}
 & \frac{1}{T}\sum_{t=1}^T\big( \frac{1}{\gamma}\|\tilde{x}_{t+1}-x_t\| + \|\nabla F(x_t) - \hat{\nabla}f(x_t,y_t)\| \big) \leq \frac{8(\Omega_1  -\Omega_T)}{3T\gamma\eta} \nonumber \\
 & \leq \frac{8(F(x_1) -F^* +g(x_1,y_1)-G(x_1))}{3T\gamma\eta}.
\end{align}
According to Jensen's inequality, we have
\begin{align} \label{eq:AG9}
 & \frac{1}{T}\sum_{t=1}^T\big( \frac{1}{\gamma}\|\tilde{x}_{t+1}-x_t\| + \|\nabla F(x_t) - \hat{\nabla}f(x_t,y_t)\| \big) \nonumber \\
 & \leq \Big(\frac{2}{T}\sum_{t=1}^T\big( \frac{1}{\gamma^2}\|\tilde{x}_{t+1}-x_t\|^2 + \|\nabla F(x_t) - \hat{\nabla}f(x_t,y_t)\|^2 \big) \Big)^{1/2} \nonumber \\
 & \leq \frac{4\sqrt{F(x_1) -F^* +g(x_1,y_1)-G(x_1)}}{\sqrt{3T\gamma\eta}}.
\end{align}

When $\mathcal{X}\subset \mathbb{R}^d$ and $\tilde{x}_{t+1}=\mathbb{P}_{\mathcal{X}}(x_t-\gamma\hat{\nabla}f(x_t,y_t))=\arg\min_{x\in \mathcal{X}}
\big\{ \langle \hat{\nabla}f(x_t,y_t), x-x_t\rangle + \frac{1}{2\gamma}\|x-x_t\|^2 \big\}$, we can define the gradient mapping
$\mathcal{G}(x_t,\hat{\nabla}f(x_t,y_t),\gamma) = \frac{1}{\gamma}\big(x_t-\mathbb{P}_{\mathcal{X}}(x_t-\gamma\hat{\nabla}f(x_t,y_t))\big)$.
Meanwhile, we can also define a gradient mapping $\mathcal{G}(x_t,\nabla F(x_t),\gamma) = \frac{1}{\gamma}\big(x_t-\mathbb{P}_{\mathcal{X}}(x_t-\gamma\nabla F(x_t))\big)$.
Then we have
\begin{align} \label{eq:AG10}
 \|\mathcal{G}(x_t,\nabla F(x_t),\gamma)\| & \leq \|\mathcal{G}(x_t,\hat{\nabla}f(x_t,y_t),\gamma)\| + \|\mathcal{G}(x_t,\nabla F(x_t),\gamma)-\mathcal{G}(x_t,\hat{\nabla}f(x_t,y_t),\gamma)\| \nonumber \\
 & = \frac{1}{\gamma} \|\tilde{x}_{t+1}-x_t\| + \frac{1}{\gamma}\|\mathbb{P}_{\mathcal{X}}(x_t-\gamma\hat{\nabla}f(x_t,y_t))-\mathbb{P}_{\mathcal{X}}(x_t-\gamma\nabla F(x_t))\| \nonumber \\
 & \leq \frac{1}{\gamma} \|\tilde{x}_{t+1}-x_t\| + \|\hat{\nabla}f(x_t,y_t)-\nabla F(x_t)\|.
\end{align}
Putting the above inequalities~(\ref{eq:AG10}) into~(\ref{eq:AG9}), we can obtain
 \begin{align}
 \frac{1}{T}\sum_{t=1}^T\|\mathcal{G}(x_t,\nabla F(x_t),\gamma)\|
 & \leq \frac{1}{T}\sum_{t=1}^T\Big(\frac{1}{\gamma} \|\tilde{x}_{t+1}-x_t\| + \|\hat{\nabla}f(x_t,y_t)-\nabla F(x_t)\| \Big) \nonumber \\
 & \leq \frac{4\sqrt{F(x_1) -F^* +g(x_1,y_1)-G(x_1)\big) }}{\sqrt{3T\gamma\eta}}.
\end{align}

When $\mathcal{X}= \mathbb{R}^d$ and $\tilde{x}_{t+1}=x_t-\gamma\hat{\nabla}f(x_t,y_t)=\arg\min_{x\in \mathbb{R}^d}
\big\{ \langle \hat{\nabla}f(x_t,y_t), x-x_t\rangle + \frac{1}{2\gamma}\|x-x_t\|^2 \big\}$, we have
$\hat{\nabla}f(x_t,y_t) = \frac{1}{\gamma}(x_t - \tilde{x}_{t+1})$, then we can get
\begin{align} \label{eq:AG11}
 \|\nabla F(x_t)\| \leq \|\hat{\nabla}f(x_t,y_t)\| + \|\hat{\nabla}f(x_t,y_t)-\nabla F(x_t)\| =
 \frac{1}{\gamma}\|x_t - \tilde{x}_{t+1}\| + \|\hat{\nabla}f(x_t,y_t)-\nabla F(x_t)\|.
\end{align}
Putting the above inequalities~(\ref{eq:AG11}) into~(\ref{eq:AG9}), we can get
 \begin{align}
 \frac{1}{T}\sum_{t=1}^T\|\nabla F(x_t)\|
 & \leq \frac{1}{T}\sum_{t=1}^T\Big(\frac{1}{\gamma} \|\tilde{x}_{t+1}-x_t\| + \|\hat{\nabla}f(x_t,y_t)-\nabla F(x_t)\| \Big) \nonumber \\
 & \leq \frac{4\sqrt{F(x_1) -F^* +g(x_1,y_1)-G(x_1)}}{\sqrt{3T\gamma\eta}}.
\end{align}

\end{proof}

\subsection{Convergence Analysis of MSGBiO Algorithm}
In this subsection, we provide the detailed convergence analysis of MSGBiO algorithm.

\begin{lemma} \label{lem:B2}
(Restatement of Lemma 5)
When the gradient estimator $w_t$ generated from Algorithm~\ref{alg:2} or \ref{alg:3}, for all $t\geq 1$,
 we have
\begin{align}
 \|w_t-\nabla F(x_t)\|^2
 & \leq 8\|u_t - \nabla_xf(x_t,y_t)\|^2 + \frac{8C^2_{fy}}{\mu^2}\|G_t-\hat{\Pi}_{C_{gxy}}\big[\nabla^2_{xy}g(x_t,y_t)\big]\|^2 \nonumber \\
 & \quad + \frac{8\kappa^2 C^2_{fy}}{\mu^2}\|H_t-
 \mathcal{S}_{[\mu,L_g]}\big[\nabla^2_{yy} g(x_t,y_t)\big]\|^2
 + 8\kappa^2 \|h_t -\Pi_{C_{fy}}\big[\nabla_yf(x_t,y_t)\big]\|^2 \nonumber \\
 & \quad  + \frac{4\hat{L}^2}{\mu}\big(g(x_t,y_t)-G(x_t)\big),
\end{align}
where $\kappa=\frac{C_{gxy}}{\mu}$.
\end{lemma}

\begin{proof}
According to Lemmas~\ref{lem:A2} and \ref{lem:D1}, we have
 \begin{align}
 \|\hat{\nabla} f(x_t,y_t)- \nabla F(x_t)\|^2 \leq \hat{L}^2\|y^*(x_t)-y_t\|^2\leq \frac{2\hat{L}^2}{\mu}\big(g(x_t,y_t)-G(x_t)\big).
\end{align}
Then we have
\begin{align} \label{eq:BB1}
 \|w_t-\nabla F(x_t)\|^2 & = \|w_t-\hat{\nabla} f(x_t,y_t) + \hat{\nabla} f(x_t,y_t)-\nabla F(x_t)\|^2 \nonumber \\
 & \leq 2\|w_t-\hat{\nabla} f(x_t,y_t)\|^2 + 2\|\hat{\nabla} f(x_t,y_t)-\nabla F(x_t)\|^2 \nonumber \\
 & \leq 2\|w_t-\hat{\nabla} f(x_t,y_t)\|^2 + \frac{4\hat{L}^2}{\mu}\big(g(x_t,y_t)-G(x_t)\big).
\end{align}

Next we consider the term $\|w_t - \hat{\nabla} f(x_t,y_t)\|^2$ when $w_t = u_t - G_t(H_t)^{-1}h_t$ generated from Algorithms~\ref{alg:2}
or~\ref{alg:3},
\begin{align} \label{eq:BB2}
 & \|w_t - \hat{\nabla} f(x_t,y_t)\|^2 \nonumber \\
 & = \|u_t - G_t(H_t)^{-1}h_t - \nabla_xf(x_t,y_t) + \hat{\Pi}_{C_{gxy}}\big[\nabla^2_{xy}g(x_t,y_t)\big] \big(\mathcal{S}_{[\mu,L_g]}\big[\nabla^2_{yy}g(x_t,y_t)\big]\big)^{-1}\Pi_{C_{fy}}\big[\nabla_yf(x_t,y_t)\big]\|^2 \nonumber \\
 & = \|u_t - \nabla_xf(x_t,y_t) - G_t(H_t)^{-1}h_t + \hat{\Pi}_{C_{gxy}}\big[\nabla^2_{xy}g(x_t,y_t)\big] \big(H_t\big)^{-1}h_t - \hat{\Pi}_{C_{gxy}}\big[\nabla^2_{xy}g(x_t,y_t)\big] \big(H_t\big)^{-1}h_t \nonumber \\
 & \quad + \hat{\Pi}_{C_{gxy}}\big[\nabla^2_{xy}g(x_t,y_t)\big] \big( \mathcal{S}_{[\mu,L_g]}\big[\nabla^2_{yy} g(x_t,y_t)\big]\big)^{-1}h_t - \hat{\Pi}_{C_{gxy}}\big[\nabla^2_{xy}g(x_t,y_t)\big] \big( \mathcal{S}_{[\mu,L_g]}\big[\nabla^2_{yy} g(x_t,y_t)\big]\big)^{-1}h_t \nonumber \\
 & \quad + \hat{\Pi}_{C_{gxy}}\big[\nabla^2_{xy}g(x_t,y_t)\big] \big(\mathcal{S}_{[\mu,L_g]}\big[\nabla^2_{yy}g(x_t,y_t)\big]\big)^{-1}\Pi_{C_{fy}}\big[\nabla_yf(x_t,y_t)\big] \|^2 \nonumber \\
 & \leq 4\|u_t - \nabla_xf(x_t,y_t)\|^2 + \frac{4C^2_{fy}}{\mu^2}\|G_t-\hat{\Pi}_{C_{gxy}}\big[\nabla^2_{xy}g(x_t,y_t)\big] \|^2 + \frac{4C^2_{gxy}C^2_{fy}}{\mu^4}\|H_t-
 \mathcal{S}_{[\mu,L_g]}\big[\nabla^2_{yy} g(x_t,y_t)\big]\|^2 \nonumber \\
 & \quad + \frac{4C^2_{gxy}}{\mu^2} \|h_t - \Pi_{C_{fy}}\big[\nabla_yf(x_t,y_t)\big] \|^2,
\end{align}
where the last inequality holds by Assumptions~\ref{ass:1}-\ref{ass:4},
and the projection operators $\hat{\Pi}_{C_{fy}}[\cdot]$
and $\mathcal{S}_{[\mu,L_g]}[\cdot]$ in Algorithms~\ref{alg:2}
and \ref{alg:3}.

By combining the above inequalities \eqref{eq:BB1} with \eqref{eq:BB2}, we have
\begin{align}
 & \|w_t-\nabla F(x_t)\|^2 \nonumber \\
 & = \|w_t-\hat{\nabla} f(x_t,y_t) + \hat{\nabla} f(x_t,y_t)-\nabla F(x_t)\|^2 \nonumber \\
 & \leq 2\|w_t-\hat{\nabla} f(x_t,y_t)\|^2 + 2\|\hat{\nabla} f(x_t,y_t)-\nabla F(x_t)\|^2 \nonumber \\
 & \leq 8\|u_t - \nabla_xf(x_t,y_t)\|^2 + \frac{8C^2_{fy}}{\mu^2}\|G_t-\hat{\Pi}_{C_{gxy}}\big[\nabla^2_{xy}g(x_t,y_t)\big] \|^2 + \frac{8C^2_{gxy}C^2_{fy}}{\mu^4}\|H_t-
 \mathcal{S}_{[\mu,L_g]}\big[\nabla^2_{yy} g(x_t,y_t)\big]\|^2 \nonumber \\
 & \quad + \frac{8C^2_{gxy}}{\mu^2} \|h_t -\Pi_{C_{fy}}\big[\nabla_yf(x_t,y_t)\big]\|^2 + \frac{4\hat{L}^2}{\mu}\big(g(x_t,y_t)-G(x_t)\big).
\end{align}

\end{proof}

\begin{lemma} \label{lem:C2}
 Assume that the stochastic partial derivatives $u_{t+1}$, $h_{t+1}$, $v_{t+1}$, $G_{t+1}$ and $H_{t+1}$ be generated from Algorithm \ref{alg:2}, we have
 \begin{align}
 \mathbb{E}\|\nabla_x f(x_{t+1},y_{t+1}) - u_{t+1}\|^2
 & \leq (1-\beta_{t+1})\mathbb{E} \|\nabla_x f(x_t,y_t) - u_t\|^2 + \beta_{t+1}^2\sigma^2  \\
 & \quad + 2L_f^2\eta_t^2/\beta_{t+1}\big(\mathbb{E}\|\tilde{x}_{t+1} - x_t\|^2 + \mathbb{E}\|\tilde{y}_{t+1} - y_t\|^2 \big), \nonumber
 \end{align}
 \begin{align}
 \mathbb{E}\|\Pi_{C_{fy}}\big[\nabla_yf(x_{t+1},y_{t+1})\big]- h_{t+1}\|^2
 & \leq (1-\hat{\beta}_{t+1}) \mathbb{E} \|\Pi_{C_{fy}}\big[\nabla_yf(x_t,y_t)\big] - h_t\|^2 + \hat{\beta}_{t+1}^2\sigma^2 \\
 &\quad + 2L_f^2\eta_t^2/\hat{\beta}_{t+1}\big(\mathbb{E}\|\tilde{x}_{t+1} - x_t\|^2 + \mathbb{E}\|\tilde{y}_{t+1} - y_t\|^2 \big), \nonumber
 \end{align}
  \begin{align}
 \mathbb{E}\|\nabla_y g(x_{t+1},y_{t+1}) - v_{t+1}\|^2
 & \leq (1-\alpha_{t+1}) \mathbb{E} \|\nabla_y g(x_t,y_t) - v_t\|^2 + \alpha_{t+1}^2\sigma^2  \\
 & \quad + 2L_g^2\eta_t^2/\alpha_{t+1}\big(\mathbb{E}\|\tilde{x}_{t+1} - x_t\|^2 + \mathbb{E}\|\tilde{y}_{t+1} - y_t\|^2 \big), \nonumber
 \end{align}
 \begin{align}
 \mathbb{E}\|\hat{\Pi}_{C_{gxy}}\big[\nabla^2_{xy}g(x_{t+1},y_{t+1})\big] - G_{t+1}\|^2
 & \leq (1-\hat{\alpha}_{t+1})\mathbb{E} \|\hat{\Pi}_{C_{gxy}}\big[\nabla^2_{xy}g(x_t,y_t)\big] - G_t\|^2 + \hat{\alpha}_{t+1}^2\sigma^2  \\
 & \quad + 2L_{gxy}^2\eta_t^2/\hat{\alpha}_{t+1}\big(\mathbb{E}\|\tilde{x}_{t+1} - x_t\|^2 + \mathbb{E}\|\tilde{y}_{t+1} - y_t\|^2 \big), \nonumber
 \end{align}
 \begin{align}
 \mathbb{E}\|\mathcal{S}_{[\mu,L_g]}\big[\nabla^2_{yy} g(x_{t+1},y_{t+1})\big] - H_{t+1}\|^2
 & \leq (1-\tilde{\alpha}_{t+1})\mathbb{E} \|\mathcal{S}_{[\mu,L_g]}\big[\nabla^2_{yy} g(x_t,y_t)\big] - H_t\|^2 + \tilde{\alpha}_{t+1}^2\sigma^2  \\
 & \quad + 2L_{gyy}^2\eta_t^2/\tilde{\alpha}_{t+1}\big(\mathbb{E}\|\tilde{x}_{t+1} - x_t\|^2 + \mathbb{E}\|\tilde{y}_{t+1} - y_t\|^2 \big). \nonumber
 \end{align}
\end{lemma}
\begin{proof}
Without loss of generality, we only consider the term $\mathbb{E}\|\mathcal{S}_{[\mu,L_g]}\big[\nabla^2_{yy} g(x_{t+1},y_{t+1})\big] - H_{t+1}\|^2$. The other terms are similar for this term.
Since $H_{t+1} = \mathcal{S}_{[\mu,L_g]}\big[ \tilde{\alpha}_{t+1}\nabla^2_{yy} g(x_{t+1},y_{t+1};\zeta_{t+1}) + (1-\tilde{\alpha}_{t+1})H_t\big]$, we have
\begin{align}
 &\mathbb{E}\|\mathcal{S}_{[\mu,L_g]}\big[\nabla^2_{yy} g(x_{t+1},y_{t+1})\big] - H_{t+1}\|^2 \nonumber \\
 & = \mathbb{E}\|\mathcal{S}_{[\mu,L_g]}\big[\nabla^2_{yy} g(x_{t+1},y_{t+1})\big] - \mathcal{S}_{[\mu,L_g]}\big[ \tilde{\alpha}_{t+1}\nabla^2_{yy} g(x_{t+1},y_{t+1};\zeta_{t+1}) + (1-\tilde{\alpha}_{t+1})H_t\big]\|^2 \nonumber \\
 & \leq \mathbb{E}\|\nabla^2_{yy} g(x_{t+1},y_{t+1}) -  \tilde{\alpha}_{t+1}\nabla^2_{yy} g(x_{t+1},y_{t+1};\zeta_{t+1}) - (1-\tilde{\alpha}_{t+1})H_t \|^2 \nonumber \\
 & = \mathbb{E}\|\tilde{\alpha}_{t+1}(\nabla^2_{yy} g(x_{t+1},y_{t+1}) - \nabla^2_{yy} g(x_{t+1},y_{t+1};\zeta_{t+1})) + (1-\tilde{\alpha}_{t+1})(\nabla^2_{yy} g(x_t,y_t) - H_t)\nonumber \\
 & \quad + (1-\tilde{\alpha}_{t+1})\big( \nabla^2_{yy} g(x_{t+1},y_{t+1}) - \nabla^2_{yy} g(x_t,y_t)\big)\|^2 \nonumber \\
 & \mathop{=}^{(i)} \mathbb{E}\| (1-\tilde{\alpha}_{t+1})(\nabla^2_{yy} g(x_t,y_t) - H_t) + (1-\tilde{\alpha}_{t+1})\big(\nabla^2_{yy} g(x_{t+1},y_{t+1}) - \nabla^2_{yy} g(x_t,y_t) \big)\|^2 \nonumber \\
 & \quad + \tilde{\alpha}_{t+1}^2\mathbb{E}\|\nabla^2_{yy} g(x_{t+1},y_{t+1}) - \nabla^2_{yy} g(x_{t+1},y_{t+1};\zeta_{t+1})\|^2 \nonumber \\
 & \leq (1-\tilde{\alpha}_{t+1})^2(1+\tilde{\alpha}_{t+1})\mathbb{E} \|\nabla^2_{yy} g(x_t,y_t) - H_t\|^2 + (1-\tilde{\alpha}_{t+1})^2(1+\frac{1}{\tilde{\alpha}_{t+1}})\mathbb{E} \|\nabla^2_{yy} g(x_{t+1},y_{t+1}) - \nabla^2_{yy} g(x_t,y_t)\|^2 \nonumber \\
 & \quad + \tilde{\alpha}_{t+1}^2\mathbb{E}\|\nabla^2_{yy} g(x_{t+1},y_{t+1}) - \nabla^2_{yy} g(x_{t+1},y_{t+1};\zeta_{t+1})\|^2 \nonumber \\
 & \mathop{\leq}^{(ii)} (1-\tilde{\alpha}_{t+1})\mathbb{E} \|\nabla^2_{yy} g(x_t,y_t) - H_t\|^2 + \frac{1}{\tilde{\alpha}_{t+1}}\mathbb{E}\|\nabla^2_{yy} g(x_{t+1},y_{t+1}) - \nabla^2_{yy} g(x_t,y_t)\|^2  + \tilde{\alpha}_{t+1}^2\sigma^2 \nonumber \\
 & \mathop{\leq}^{(iii)} (1-\tilde{\alpha}_{t+1}) \mathbb{E} \|\nabla^2_{yy} g(x_t,y_t) - H_t\|^2 + \frac{2L_{gyy}^2\eta_t^2}{\tilde{\alpha}_{t+1}}\big(\mathbb{E}\|\tilde{x}_{t+1} - x_t\|^2 + \mathbb{E}\|\tilde{y}_{t+1} - y_t\|^2 \big) + \tilde{\alpha}_{t+1}^2\sigma^2, \nonumber
\end{align}
where the equality (i) is due to $\mathbb{E}_{\zeta_{t+1}}[\nabla^2_{yy} g(x_{t+1},y_{t+1};\zeta_{t+1})]=\nabla^2_{yy} g(x_{t+1},y_{t+1})$; the inequality (ii) holds by $0\leq \tilde{\alpha}_{t+1} \leq 1$ such that  $(1-\tilde{\alpha}_{t+1})^2(1+\tilde{\alpha}_{t+1})=1-\tilde{\alpha}_{t+1}-\tilde{\alpha}_{t+1}^2+
  \tilde{\alpha}_{t+1}^3\leq 1-\tilde{\alpha}_{t+1}$ and $(1-\tilde{\alpha}_{t+1})^2(1+\frac{1}{\tilde{\alpha}_{t+1}}) \leq (1-\tilde{\alpha}_{t+1})(1+\frac{1}{\tilde{\alpha}_{t+1}}) = -\tilde{\alpha}_{t+1}+\frac{1}{\tilde{\alpha}_{t+1}}\leq \frac{1}{\tilde{\alpha}_{t+1}}$, and the inequality (iii) holds by Assumption~\ref{ass:4} and $x_{t+1}=x_t-\eta_t(\tilde{x}_{t+1}-x_t)$, $y_{t+1}=y_t-\eta_t(\tilde{y}_{t+1}-y_t)$.
\end{proof}

\begin{theorem} \label{th:A2}
(Restatement of Theorem 2)
 Under the above Assumptions (\ref{ass:1}-\ref{ass:5}, \ref{ass:8}), in the Algorithm \ref{alg:2}, let $\eta_t=\frac{k}{(m+t)^{1/2}}$ for all $t\geq 0$, $\beta_{t+1}=c_1\eta_t$, $\hat{\beta}_{t+1}=c_2\eta_t$, $\alpha_{t+1}=c_3\eta_t$, $\hat{\alpha}_{t+1}=c_4\eta_t$, $\tilde{\alpha}_{t+1}=c_5\eta_t$, $10 \leq c_1 \leq \frac{m^{1/2}}{k}$, $10\kappa^2 \leq c_2 \leq \frac{m^{1/2}}{k}$, $1 \leq c_3 \leq \frac{m^{1/2}}{k}$, $\frac{10C^2_{fy}}{\mu^2} \leq c_4 \leq \frac{m^{1/2}}{k}$, $\frac{10C^2_{fy}\kappa^2}{\mu^2} \leq c_5 \leq \frac{m^{1/2}}{k}$, $m\geq \max\big(k^2, (c_1k)^2,(c_2k)^2, (c_3k)^2, (c_4k)^2, (c_5k)^2\big)$, $k>0$, $0< \gamma \leq \min\big(\frac{m^{1/2}}{2L_Fk},\frac{\lambda\mu^2}{16\hat{L}^2},\frac{\sqrt{5}}{4\breve{L}},\frac{1}{32L^2_g\lambda},\frac{5}{8\breve{L}^2\lambda}\big)$ and $0<\lambda\leq \min\big(\frac{m^{1/2}}{2L_gk},\frac{1}{4L_g}\big)$. When $\mathcal{X}\subseteq \mathbb{R}^d$, we can get
\begin{align}
  \frac{1}{T}\sum_{t=1}^T\mathbb{E}\|\mathcal{G}(x_t,\nabla F(x_t),\gamma)\|
 & \leq \frac{1}{T} \sum_{t=1}^T \mathbb{E} \big[ \|w_t-\nabla F(x_t)\| + \frac{1}{\gamma}\|\tilde{x}_{t+1}-x_t\|  \big]   \leq \frac{\sqrt{2M}m^{1/4}}{\sqrt{T}} + \frac{\sqrt{2M}}{T^{1/4}};
\end{align}
When $\mathcal{X}= \mathbb{R}^d$, we can get
\begin{align}
  \frac{1}{T}\sum_{t=1}^T\mathbb{E}\|\nabla F(x_t)\|
 & \leq \frac{1}{T} \sum_{t=1}^T \mathbb{E} \big[ \|w_t-\nabla F(x_t)\| + \frac{1}{\gamma}\|\tilde{x}_{t+1}-x_t\|  \big]   \leq \frac{\sqrt{2M}m^{1/4}}{\sqrt{T}} + \frac{\sqrt{2M}}{T^{1/4}},
\end{align}
where $M= \frac{4(F(x_1)- F^*+g(x_1,y_1)-G(x_1))}{k\gamma} + \frac{16\sigma^2}{k} + \frac{4\lambda\sigma^2}{\gamma k} + \frac{16m\sigma^2\ln(m+T)}{k}+ \frac{4m\lambda\sigma^2\ln(m+T)}{k\gamma}$
and $\breve{L}^2=L^2_f+\frac{L^2_f}{\kappa^2}+\frac{\mu^2L^2_{gxy}}{C^2_{fy}}+\frac{\mu^2L^2_{gyy}}{C^2_{fy}\kappa^2}$.
\end{theorem}

\begin{proof}
According the line 4 of Algorithm~\ref{alg:2}, we have
\begin{align} \label{eq:F1}
\tilde{x}_{t+1} = \arg\min_{x\in \mathcal{X}}
\Big\{ \langle w_t, x-x_t\rangle + \frac{1}{2\gamma}\|x-x_t\|^2 \Big\}.
\end{align}
By using the optimal condition of the subproblem~(\ref{eq:F1}), we have
\begin{align}
 \langle w_t + \frac{1}{\gamma}(\tilde{x}_{t+1}-x_t), x_t-\tilde{x}_{t+1}\rangle \geq 0,
\end{align}
then we have
\begin{align}
 \langle w_t,\tilde{x}_{t+1} -x_t\rangle \leq -\frac{1}{\gamma}\|\tilde{x}_{t+1}-x_t\|^2.
\end{align}
According to Lemma \ref{lem:C1}, i.e., the function $F(x)$ has $L_F$-Lipschitz continuous gradient,
we have
\begin{align} \label{eq:F2}
  F(x_{t+1}) & \leq F(x_t) + \langle \nabla F(x_t), x_{t+1}-x_t\rangle + \frac{L_F}{2}\|x_{t+1}-x_t\|^2 \nonumber \\
  & = F(x_t) + \eta_t \langle \nabla F(x_t) -w_t +w_t, \tilde{x}_{t+1}-x_t\rangle + \frac{\eta^2_tL_F}{2}\|\tilde{x}_{t+1}-x_t\|^2 \nonumber \\
  & \leq F(x_t) + \eta_t \langle \nabla F(x_t) -w_t , \tilde{x}_{t+1}-x_t\rangle - \frac{\eta_t}{\gamma}\|\tilde{x}_{t+1}-x_t\|^2 + \frac{\eta^2_tL_F}{2}\|\tilde{x}_{t+1}-x_t\|^2 \nonumber \\
  & \leq F(x_t) + \eta_t\gamma\|\nabla F(x_t) -w_t\|^2 + \frac{\eta_t}{4\gamma}\|\tilde{x}_{t+1}-x_t\|^2 - \frac{\eta_t}{\gamma}\|\tilde{x}_{t+1}-x_t\|^2 + \frac{\eta^2_tL_F}{2}\|\tilde{x}_{t+1}-x_t\|^2 \nonumber \\
  & \leq F(x_t) + \eta_t\gamma\|\nabla F(x_t) -w_t\|^2 - \frac{\eta_t}{2\gamma}\|\tilde{x}_{t+1}-x_t\|^2,
\end{align}
where the last inequality holds by the above inequality~(\ref{eq:F1}), and the last inequality holds by $\gamma\leq \frac{1}{2L_F\eta_t}$.

According to the Lemma~\ref{lem:B2}, we have
\begin{align} \label{eq:F3}
  F(x_{t+1})  & \leq F(x_t) + 8\eta_t\gamma\|u_t - \nabla_xf(x_t,y_t)\|^2 + \frac{8C^2_{fy}\eta_t\gamma}{\mu^2}\|G_t-\hat{\Pi}_{C_{gxy}}\big[\nabla^2_{xy}g(x_t,y_t)\big]\|^2 \nonumber \\
  & \quad + \frac{8\kappa^2 C^2_{fy}\eta_t\gamma}{\mu^2}\|H_t-
 \mathcal{S}_{[\mu,L_g]}\big[\nabla^2_{yy} g(x_t,y_t)\big]\|^2 + 8\kappa^2\eta_t\gamma \|h_t -\Pi_{C_{fy}}\big[\nabla_yf(x_t,y_t)\big]  \|^2 \nonumber \\
 & \quad  + \frac{4\hat{L}^2\eta_t\gamma}{\mu}\big(g(x_t,y_t)-G(x_t)\big) - \frac{\eta_t}{2\gamma}\|\tilde{x}_{t+1}-x_t\|^2.
\end{align}

According to Lemma~\ref{lem:E1}, we have
\begin{align}
& g(x_{t+1},y_{t+1})-G(x_{t+1})-\big(g(x_t,y_t)-G(x_t)\big) \nonumber \\
& \leq -\frac{\eta_t\lambda\mu}{2} \big(g(x_t,y_t) -G(x_t)\big) + \frac{\eta_t}{8\gamma}\|\tilde{x}_{t+1}-x_t\|^2   -\frac{\eta_t}{4\lambda}\|\tilde{y}_{t+1}-y_t\|^2+ \eta_t\lambda\|\nabla_y g(x_t,y_t)-v_t\|^2.
\end{align}

In our Algorithm \ref{alg:2}, let $\eta_t=\frac{k}{(m+t)^{1/2}}$ for all $t\geq 0$, $\beta_{t+1}=c_1\eta_t$, $\hat{\beta}_{t+1}=c_2\eta_t$, $\alpha_{t+1}=c_3\eta_t$, $\hat{\alpha}_{t+1}=c_4\eta_t$, $\tilde{\alpha}_{t+1}=c_5\eta_t$.

Since $\eta_t=\frac{k}{(m+t)^{1/2}}$ on $t$ is decreasing and $m\geq k^2$, we have $\eta_t \leq \eta_0 = \frac{k}{m^{1/2}} \leq 1$ and $\gamma \leq \frac{m^{1/2}}{2L_Fk}\leq \frac{1}{2L_F\eta_0} \leq \frac{1}{2L_F\eta_t}$ for any $t\geq 0$. Similarly, we can get $\lambda \leq \frac{m^{1/2}}{2L_gk}\leq \frac{1}{2L_g\eta_t}$ for any $t\geq 0$.
 Due to $0 < \eta_t \leq 1$ and $m\geq (c_1k)^2$, we have $\beta_{t+1} = c_1\eta_t \leq \frac{c_1k}{m^{1/2}}\leq 1$.
 Similarly, due to $m\geq \max\big( (c_2k)^2, (c_3k)^2, (c_4k)^2, (c_5k)^2 \big)$, we have $\hat{\beta}_{t+1}\leq 1$, $\alpha_{t+1}\leq 1$,  $\hat{\alpha}_{t+1}\leq 1$, and $\tilde{\alpha}_{t+1} \leq 1$. At the same time, we have $c_i\leq \frac{m^{1/2}}{k}$ for all $i=1,2,\cdots,5$.

 According to Lemma \ref{lem:C2}, we have
 \begin{align} \label{eq:F4}
  & \mathbb{E} \|\nabla_x f(x_{t+1},y_{t+1}) - u_{t+1}\|^2 - \mathbb{E} \|\nabla_x f(x_t,y_t) - u_t\|^2  \\
  & \leq -\beta_{t+1}\mathbb{E} \|\nabla_x f(x_t,y_t) -u_t\|^2 + 2L^2_f\eta^2_t/\beta_{t+1}\big(\|\tilde{x}_{t+1}-x_t\|^2 + \|\tilde{y}_{t+1}-y_t\|^2\big) + \beta_{t+1}^2\sigma^2  \nonumber \\
  & = -c_1 \eta_t\mathbb{E} \|\nabla_x f(x_t,y_t) -u_t\|^2 + 2L^2_f\eta_t/c_1\big(\|\tilde{x}_{t+1}-x_t\|^2 + \|\tilde{y}_{t+1}-y_t\|^2\big) + c_1^2\eta_t^2\sigma^2 \nonumber \\
  & \leq -10\eta_t\mathbb{E} \|\nabla_x f(x_t,y_t) -u_t\|^2 + \frac{L^2_f\eta_t}{5}\big(\|\tilde{x}_{t+1}-x_t\|^2 + \|\tilde{y}_{t+1}-y_t\|^2\big) + \frac{m\eta^2_t\sigma^2}{k^2}, \nonumber
 \end{align}
 where the above equality holds by $\beta_{t+1}=c_1\eta_t$, and the last inequality is due to $10 \leq c_1 \leq \frac{m^{1/2}}{k}$.
 Similarly, given $ 10\kappa^2 \leq c_2 \leq \frac{m^{1/2}}{k}$, we have
 \begin{align} \label{eq:F5}
  & \mathbb{E} \|\Pi_{C_{fy}}\big[\nabla_yf(x_{t+1},y_{t+1})\big] - h_{t+1}\|^2 -  \mathbb{E} \|\Pi_{C_{fy}}\big[\nabla_yf(x_t,y_t)\big] - h_t\|^2  \\
  & \leq - 10\kappa^2\eta_t\mathbb{E} \|\Pi_{C_{fy}}\big[\nabla_yf(x_t,y_t)\big] -h_t\|^2 +
  \frac{L^2_f\eta_t}{5\kappa^2}\big(\|\tilde{x}_{t+1}-x_t\|^2 + \|\tilde{y}_{t+1}-y_t\|^2\big) + \frac{m\eta^2_t\sigma^2}{k^2}. \nonumber
 \end{align}
 Given $1 \leq c_3 \leq \frac{m^{1/2}}{k}$, we have
 \begin{align} \label{eq:F6}
  & \mathbb{E} \|\nabla_y g(x_{t+1},y_{t+1}) - v_{t+1}\|^2 -  \mathbb{E} \|\nabla_y g(x_t,y_t) - v_t\|^2  \\
  & \leq - \eta_t \mathbb{E} \|\nabla_y g(x_t,y_t) -v_t\|^2 +
   2 L^2_g\eta_t\big(\|\tilde{x}_{t+1}-x_t\|^2 + \|\tilde{y}_{t+1}-y_t\|^2\big) + \frac{m\eta^2_t\sigma^2}{k^2}. \nonumber
 \end{align}
 Given $ \frac{10C^2_{fy}}{\mu^2} \leq c_4 \leq \frac{m^{1/2}}{k}$, we have
 \begin{align}  \label{eq:F7}
  & \mathbb{E} \|\hat{\Pi}_{C_{gxy}}\big[\nabla^2_{xy}g(x_{t+1},y_{t+1})\big]-G_{t+1}\|^2 - \mathbb{E} \|\hat{\Pi}_{C_{gxy}}\big[\nabla^2_{xy}g(x_t,y_t)\big]-G_t\|^2  \\
  & \leq - \frac{10C^2_{fy}\eta_t}{\mu^2} \mathbb{E} \|\hat{\Pi}_{C_{gxy}}\big[\nabla^2_{xy}g(x_t,y_t)\big] -G_t\|^2 +
  \frac{\mu^2L^2_{gxy}\eta_t}{5C^2_{fy}} \big(\|\tilde{x}_{t+1}-x_t\|^2 + \|\tilde{y}_{t+1}-y_t\|^2\big) + \frac{m\eta^2_t\sigma^2}{k^2}. \nonumber
 \end{align}
 Given $ \frac{10\kappa^2 C^2_{fy}}{\mu^2} \leq c_5 \leq \frac{m^{1/2}}{k}$, we have
 \begin{align}  \label{eq:F8}
  & \mathbb{E} \|\mathcal{S}_{[\mu,L_g]}\big[\nabla^2_{yy} g(x_{t+1},y_{t+1})\big] - H_{t+1}\|^2  -  \mathbb{E} \|\mathcal{S}_{[\mu,L_g]}\big[\nabla^2_{yy} g(x_t,y_t)\big] - H_t\|^2  \\
  & \leq  -\frac{10\kappa^2 C^2_{fy}\eta_t}{\mu^2}  \mathbb{E} \|\mathcal{S}_{[\mu,L_g]}\big[\nabla^2_{yy} g(x_t,y_t)\big] -H_t\|^2 +
  \frac{\mu^2L^2_{gyy}\eta_t}{5\kappa^2 C^2_{fy}}\big(\|\tilde{x}_{t+1}-x_t\|^2 + \|\tilde{y}_{t+1}-y_t\|^2\big) + \frac{m\eta^2_t\sigma^2}{k^2}. \nonumber
 \end{align}

Next, we define a useful \emph{Lyapunov} function (i.e., potential function), for any $t\geq 1$
\begin{align}
 \Phi_t & = \mathbb{E}\Big [F(x_t) + g(x_t,y_t)-G(x_t) + \gamma \big(\|\nabla_x f(x_t,y_t)-u_t\|^2 + \|\Pi_{C_{fy}}\big[\nabla_yf(x_t,y_t)\big]-h_t\|^2 \nonumber \\
 & \qquad + \|\hat{\Pi}_{C_{gxy}}\big[\nabla^2_{xy}g(x_t,y_t)\big] - G_t\|^2 + \|\mathcal{S}_{[\mu,L_g]}\big[\nabla^2_{yy} g(x_t,y_t)\big] - H_t\|^2\big) + \lambda\|\nabla_y g(x_t,y_t) - v_t\|^2 \Big]. \nonumber
\end{align}
Then we have
\begin{align} \label{eq:F9}
 & \Phi_{t+1} - \Phi_t \nonumber \\
 & = F(x_{t+1}) - F(x_t) + g(x_{t+1},y_{t+1})-G(x_{t+1}) - \big( g(x_t,y_t)-G(x_t) \big)
 + \gamma\big(\mathbb{E}\|\nabla_x f(x_{t+1},y_{t+1})-u_{t+1}\|^2 \nonumber \\
 & \quad - \mathbb{E}\|\nabla_x f(x_t,y_t)-u_t\|^2 + \mathbb{E}\|\Pi_{C_{fy}}\big[\nabla_yf(x_{t+1},y_{t+1})\big]-h_{t+1}\|^2
 - \mathbb{E}\|\Pi_{C_{fy}}\big[\nabla_yf(x_t,y_t)\big]-h_t\|^2 \nonumber \\
 & \quad + \mathbb{E}\|\nabla_y g(x_{t+1},y_{t+1})-v_{t+1}\|^2
 - \mathbb{E}\|\nabla_y g(x_t,y_t)-v_t\|^2 + \mathbb{E}\|\hat{\Pi}_{C_{gxy}}\big[\nabla^2_{xy}g(x_{t+1},y_{t+1})\big]-G_{t+1}\|^2
  \nonumber \\
  & \quad - \mathbb{E}\|\hat{\Pi}_{C_{gxy}}\big[\nabla^2_{xy}g(x_t,y_t)\big]-G_t\|^2 +\mathbb{E}\|\mathcal{S}_{[\mu,L_g]}\big[\nabla^2_{yy} g(x_{t+1},y_{t+1})\big]-H_{t+1}\|^2
 - \mathbb{E}\|\mathcal{S}_{[\mu,L_g]}\big[\nabla^2_{yy} g(x_t,y_t)\big]-H_t\|^2 \big) \nonumber \\
 & \leq 8\eta_t\gamma\|u_t - \nabla_xf(x_t,y_t)\|^2 + \frac{8C^2_{fy}\eta_t\gamma}{\mu^2}\|G_t-\hat{\Pi}_{C_{gxy}}\big[\nabla^2_{xy}g(x_t,y_t)\big]\|^2 + \frac{8\kappa^2 C^2_{fy}\eta_t\gamma}{\mu^2}\|H_t-
 \mathcal{S}_{[\mu,L_g]}\big[\nabla^2_{yy} g(x_t,y_t)\big]\|^2 \nonumber \\
 & \quad + 8\kappa^2\eta_t\gamma \|h_t -\Pi_{C_{fy}}\big[\nabla_yf(x_t,y_t)\big] \|^2 + \frac{4\hat{L}^2\eta_t\gamma}{\mu}\big(g(x_t,y_t)-G(x_t)\big) - \frac{\eta_t}{2\gamma}\|\tilde{x}_{t+1}-x_t\|^2 \nonumber \\
 & \quad  -\frac{\eta_t\lambda\mu}{2} \big(g(x_t,y_t) -G(x_t)\big) + \frac{\eta_t}{8\gamma}\|\tilde{x}_{t+1}-x_t\|^2 -\frac{\eta_t}{4\lambda}\|\tilde{y}_{t+1}-y_t\|^2+ \eta_t\lambda\|\nabla_y g(x_t,y_t)-v_t\|^2 \nonumber \\
 & \quad + \gamma \bigg( -10\eta_t\mathbb{E} \|\nabla_x f(x_t,y_t) -u_t\|^2 + \frac{L^2_f\eta_t}{5}\big(\|\tilde{x}_{t+1}-x_t\|^2 + \|\tilde{y}_{t+1}-y_t\|^2\big) + \frac{m\eta^2_t\sigma^2}{k^2} \nonumber \\
 & \qquad - 10\kappa^2\eta_t\mathbb{E} \|\Pi_{C_{fy}}\big[\nabla_yf(x_t,y_t)\big] -h_t\|^2 +
  \frac{L^2_f\eta_t}{5\kappa^2}\big(\|\tilde{x}_{t+1}-x_t\|^2 + \|\tilde{y}_{t+1}-y_t\|^2\big) + \frac{m\eta^2_t\sigma^2}{k^2}  \nonumber \\
 & \qquad - \frac{10C^2_{fy}\eta_t}{\mu^2} \mathbb{E} \|\hat{\Pi}_{C_{gxy}}\big[\nabla^2_{xy}g(x_t,y_t)\big] -G_t\|^2 +
  \frac{\mu^2L^2_{gxy}\eta_t}{5C^2_{fy}} \big(\|\tilde{x}_{t+1}-x_t\|^2 + \|\tilde{y}_{t+1}-y_t\|^2\big) + \frac{m\eta^2_t\sigma^2}{k^2} \nonumber \\
 & \qquad -\frac{10\kappa^2 C^2_{fy}\eta_t}{\mu^2}  \mathbb{E} \|\mathcal{S}_{[\mu,L_g]}\big[\nabla^2_{yy} g(x_t,y_t)\big] -H_t\|^2 +
  \frac{\mu^2L^2_{gyy}\eta_t}{5\kappa^2 C^2_{fy}}\big(\|\tilde{x}_{t+1}-x_t\|^2 + \|\tilde{y}_{t+1}-y_t\|^2\big) + \frac{m\eta^2_t\sigma^2}{k^2} \bigg) \nonumber \\
 & \quad + \lambda \Big(- \eta_t \mathbb{E} \|\nabla_y g(x_t,y_t) -v_t\|^2 +
   2 L^2_g\eta_t\big(\|\tilde{x}_{t+1}-x_t\|^2 + \|\tilde{y}_{t+1}-y_t\|^2\big) + \frac{m\eta^2_t\sigma^2}{k^2} \Big)\nonumber \\
 & = - \frac{\gamma\eta_t}{4} \bigg( 8\|u_t - \nabla_xf(x_t,y_t)\|^2 + \frac{8C^2_{fy}}{\mu^2}\|G_t-\hat{\Pi}_{C_{gxy}}\big[\nabla^2_{xy}g(x_t,y_t)\big]\|^2 + \frac{8\kappa^2 C^2_{fy}}{\mu^2}\|H_t-
 \mathcal{S}_{[\mu,L_g]}\big[\nabla^2_{yy} g(x_t,y_t)\big]\|^2 \nonumber \\
 & \qquad + 8\kappa^2 \|h_t -\Pi_{C_{fy}}\big[\nabla_yf(x_t,y_t)\big] \|^2 \bigg) - \big(\frac{\lambda\mu\eta_t}{2}-\frac{4\hat{L}^2\eta_t\gamma}{\mu}\big)\big(g(x_t,y_t)-G(x_t)\big)   \nonumber \\
 & \quad - \big( \frac{1}{2\gamma} - \frac{1}{8\gamma} - \frac{\breve{L}^2\gamma}{5} - 2L^2_g\lambda\big)\eta_t\|\tilde{x}_{t+1}-x_t\|^2 - \big( \frac{1}{4\lambda} -
 \frac{\breve{L}^2\gamma}{5} - 2L^2_g\lambda\big)\eta_t\|\tilde{y}_{t+1}-y_t\|^2 + \frac{4m\gamma\sigma^2}{k^2}\eta^2_t+ \frac{m\lambda\sigma^2}{k^2}\eta^2_t  \nonumber \\
 & \leq - \frac{\gamma\eta_t}{4} \bigg( 8\|u_t - \nabla_xf(x_t,y_t)\|^2 + \frac{8C^2_{fy}}{\mu^2}\|G_t-\hat{\Pi}_{C_{gxy}}\big[\nabla^2_{xy}g(x_t,y_t)\big]\|^2 + \frac{8\kappa^2 C^2_{fy}}{\mu^2}\|H_t-
 \mathcal{S}_{[\mu,L_g]}\big[\nabla^2_{yy} g(x_t,y_t)\big]\|^2 \nonumber \\
 & \qquad + 8\kappa^2 \|h_t - \Pi_{C_{fy}}\big[\nabla_yf(x_t,y_t)\big]\|^2 \bigg) - \frac{\lambda\mu\eta_t}{4}\big(g(x_t,y_t)-G(x_t)\big)   \nonumber \\
 & \quad - \frac{\eta_t}{4\gamma}\|\tilde{x}_{t+1}-x_t\|^2 + \frac{4m\gamma\sigma^2}{k^2}\eta^2_t+ \frac{m\lambda\sigma^2}{k^2}\eta^2_t ,
\end{align}
where $\breve{L}^2=L^2_f+\frac{L^2_f}{\kappa^2}+\frac{\mu^2L^2_{gxy}}{C^2_{fy}}+\frac{\mu^2L^2_{gyy}}{C^2_{fy}\kappa^2}$, the first inequality holds by the above inequalities \eqref{eq:F5}-\eqref{eq:F8}; the last inequality is due to $0< \gamma \leq \min\big(\frac{\lambda\mu^2}{16\hat{L}^2},\frac{\sqrt{5}}{4\breve{L}},\frac{1}{32L^2_g\lambda},\frac{5}{8\breve{L}^2\lambda}\big)$ and $0<\lambda\leq \frac{1}{4L_g}$.

Since $\gamma\leq \frac{\lambda\mu^2}{16\hat{L}^2}$, according to the above inequality~(\ref{eq:F9}), we have
\begin{align} \label{eq:F10}
 &  \frac{\gamma\eta_t}{4} \bigg( 8\|u_t - \nabla_xf(x_t,y_t)\|^2 + \frac{8C^2_{fy}}{\mu^2}\|G_t-\hat{\Pi}_{C_{gxy}}\big[\nabla^2_{xy}g(x_t,y_t)\big]\|^2 + \frac{8\kappa^2 C^2_{fy}}{\mu^2}\|H_t-
 \mathcal{S}_{[\mu,L_g]}\big[\nabla^2_{yy} g(x_t,y_t)\big]\|^2 \nonumber \\
 & \qquad + 8\kappa^2 \|h_t -\Pi_{C_{fy}}\big[\nabla_yf(x_t,y_t)\big] \|^2 + \frac{16\hat{L}^2}{\mu}\big(g(x_t,y_t)-G(x_t)\big) \bigg)  + \frac{\eta_t}{4\gamma}\|\tilde{x}_{t+1}-x_t\|^2  \nonumber \\
 &  \leq \frac{\gamma\eta_t}{4} \bigg( 8\|u_t - \nabla_xf(x_t,y_t)\|^2 + \frac{8C^2_{fy}}{\mu^2}\|G_t-\hat{\Pi}_{C_{gxy}}\big[\nabla^2_{xy}g(x_t,y_t)\big]\|^2 + \frac{8\kappa^2 C^2_{fy}}{\mu^2}\|H_t-
 \mathcal{S}_{[\mu,L_g]}\big[\nabla^2_{yy} g(x_t,y_t)\big]\|^2 \nonumber \\
 & \qquad + 8\kappa^2 \|h_t - \Pi_{C_{fy}}\big[\nabla_yf(x_t,y_t)\big]\|^2 \bigg) + \frac{\lambda\mu\eta_t}{4}\big(g(x_t,y_t)-G(x_t)\big) + \frac{\eta_t}{4\gamma}\|\tilde{x}_{t+1}-x_t\|^2  \nonumber \\
 &   \leq \Phi_t -\Phi_{t+1} + \frac{4m\gamma\sigma^2}{k^2}\eta^2_t+ \frac{m\lambda\sigma^2}{k^2}\eta^2_t.
\end{align}

According to Lemma~\ref{lem:B2}, then we have
\begin{align} \label{eq:F11}
 \frac{\gamma\eta_t}{4}\|w_t-\nabla F(x_t)\|^2 + \frac{\eta_t}{4\gamma}\|\tilde{x}_{t+1}-x_t\|^2  \leq \Phi_t -\Phi_{t+1} + \frac{4m\gamma\sigma^2}{k^2}\eta^2_t+ \frac{m\lambda\sigma^2}{k^2}\eta^2_t.
\end{align}
Taking average over $t=1,2,\cdots,T$ on both sides of \eqref{eq:F11}, we have
\begin{align}
  \frac{1}{T} \sum_{t=1}^T \mathbb{E} \big[ \eta_t\|w_t-\nabla F(x_t)\|^2 + \frac{\eta_t}{\gamma^2}\|\tilde{x}_{t+1}-x_t\|^2  \big]
  \leq   \frac{4(\Phi_1 - \Phi_{T+1})}{T\gamma} + \frac{1}{T}\sum_{t=1}^T\big(\frac{16m\sigma^2}{k^2} + \frac{4m\lambda\sigma^2}{\gamma k^2}\big)\eta^2_t. \nonumber
\end{align}
Since $\eta_t=\frac{k}{(m+t)^{1/2}}$ is decreasing on $t$, i.e., $\eta_T^{-1} \geq \eta_t^{-1}$ for any $0\leq t\leq T$, we have
\begin{align}
 & \frac{1}{T} \sum_{t=1}^T \mathbb{E} \big[ \|w_t-\nabla F(x_t)\|^2 + \frac{1}{\gamma^2}\|\tilde{x}_{t+1}-x_t\|^2  \big] \nonumber \\
 & \leq \frac{4(\Phi_1 - \Phi_{T+1})}{T\eta_T\gamma} + \frac{1}{T}\sum_{t=1}^T\big(\frac{16m\sigma^2}{k^2\eta_T} + \frac{4m\lambda\sigma^2}{\gamma k^2\eta_T}\big)\eta^2_t \nonumber \\
 & \leq \frac{4(F(x_1)- F^*+g(x_1,y_1)-G(x_1) + 4\gamma\sigma^2+ \lambda\sigma^2)}{T\eta_T\gamma} + \big(\frac{16m\sigma^2}{k^2\eta_T} + \frac{4m\lambda\sigma^2}{\gamma k^2\eta_T}\big)\frac{1}{T}\sum_{t=1}^T\eta^2_t \nonumber \\
 & \leq \frac{4(F(x_1)- F^*+g(x_1,y_1)-G(x_1) + 4\gamma\sigma^2+ \lambda\sigma^2)}{T\eta_T\gamma} + \big(\frac{16m\sigma^2}{k^2\eta_T} + \frac{4m\lambda\sigma^2}{\gamma k^2\eta_T}\big)\frac{1}{T}\int_{1}^T\frac{k^2}{m+t}dt \nonumber \\
 & \leq \frac{4(F(x_1)- F^*+g(x_1,y_1)-G(x_1) + 4\gamma\sigma^2+ \lambda\sigma^2)}{T\eta_T\gamma} + \big(\frac{16m\sigma^2}{\eta_T} + \frac{4m\lambda\sigma^2}{\gamma\eta_T}\big)\frac{\ln(m+T)}{T} \nonumber\\
 & \leq \Big(\frac{4(F(x_1)- F^*+g(x_1,y_1)-G(x_1))}{k\gamma} + \frac{16\sigma^2}{k} + \frac{4\lambda\sigma^2}{\gamma k}+\big(\frac{16m\sigma^2}{\eta_T} + \frac{4m\lambda\sigma^2}{\gamma\eta_T}\big)\ln(m+T)\Big)\frac{(m+T)^{1/2}}{T}
\end{align}
where the second inequality holds by Assumptions~\ref{ass:1}-\ref{ass:3}.
Let $M= \frac{4(F(x_1)- F^*+g(x_1,y_1)-G(x_1))}{k\gamma} + \frac{16\sigma^2}{k} + \frac{4\lambda\sigma^2}{\gamma k} + \frac{16m\sigma^2\ln(m+T)}{k}+ \frac{4m\lambda\sigma^2\ln(m+T)}{k\gamma}$, we have
\begin{align}
 \frac{1}{T} \sum_{t=1}^T \mathbb{E} \big[ \|w_t-\nabla F(x_t)\|^2 + \frac{1}{\gamma^2}\|\tilde{x}_{t+1}-x_t\|^2  \big]
  \leq \frac{M(m+T)^{1/2}}{T}
\end{align}
According to the Jensen's inequality, we can get
\begin{align} \label{eq:F12}
 & \frac{1}{T} \sum_{t=1}^T \mathbb{E} \big[ \|w_t-\nabla F(x_t)\| + \frac{1}{\gamma}\|\tilde{x}_{t+1}-x_t\|  \big]  \nonumber \\
 & \leq \Big(\frac{2}{T} \sum_{t=1}^T \mathbb{E} \big[ \|w_t-\nabla F(x_t)\|^2 + \frac{1}{\gamma^2}\|\tilde{x}_{t+1}-x_t\|^2  \big] \Big)^{1/2} \nonumber \\
 & \leq \frac{\sqrt{2M}(m+T)^{1/4}}{\sqrt{T}} \leq \frac{\sqrt{2M}m^{1/4}}{\sqrt{T}} + \frac{\sqrt{2M}}{T^{1/4}}.
\end{align}

When $\mathcal{X}\subset \mathbb{R}^d$ and $\tilde{x}_{t+1}=\mathbb{P}_{\mathcal{X}}(x_t-\gamma w_t)=\arg\min_{x\in \mathcal{X}}
\big\{ \langle w_t, x-x_t\rangle + \frac{1}{2\gamma}\|x-x_t\|^2 \big\}$, we can define the gradient mapping
$\mathcal{G}(x_t,w_t,\gamma) = \frac{1}{\gamma}\big(x_t-\mathbb{P}_{\mathcal{X}}(x_t-\gamma w_t\big)$.
Meanwhile, we can also define a gradient mapping $\mathcal{G}(x_t,\nabla F(x_t),\gamma) = \frac{1}{\gamma}\big(x_t-\mathbb{P}_{\mathcal{X}}(x_t-\gamma\nabla F(x_t))\big)$.
Then we have
\begin{align} \label{eq:F13}
 \|\mathcal{G}(x_t,\nabla F(x_t),\gamma)\| & \leq \|\mathcal{G}(x_t,w_t,\gamma)\| + \|\mathcal{G}(x_t,\nabla F(x_t),\gamma)-\mathcal{G}(x_t,w_t,\gamma)\| \nonumber \\
 & = \frac{1}{\gamma} \|\tilde{x}_{t+1}-x_t\| + \frac{1}{\gamma}\|\mathbb{P}_{\mathcal{X}}(x_t-\gamma w_t-\mathbb{P}_{\mathcal{X}}(x_t-\gamma\nabla F(x_t))\| \nonumber \\
 & \leq \frac{1}{\gamma} \|\tilde{x}_{t+1}-x_t\| + \|w_t-\nabla F(x_t)\|.
\end{align}
Putting the above inequalities~(\ref{eq:F13}) into~(\ref{eq:F12}), we can obtain
 \begin{align}
 \frac{1}{T}\sum_{t=1}^T\mathbb{E}\|\mathcal{G}(x_t,\nabla F(x_t),\gamma)\|
 & \leq \frac{1}{T} \sum_{t=1}^T \mathbb{E} \big[ \|w_t-\nabla F(x_t)\| + \frac{1}{\gamma}\|\tilde{x}_{t+1}-x_t\|  \big]  \nonumber \\
 & \leq \frac{\sqrt{2M}(m+T)^{1/4}}{\sqrt{T}} \leq \frac{\sqrt{2M}m^{1/4}}{\sqrt{T}} + \frac{\sqrt{2M}}{T^{1/4}}.
\end{align}

When $\mathcal{X}= \mathbb{R}^d$ and $\tilde{x}_{t+1}=x_t-\gamma w_t=\arg\min_{x\in \mathbb{R}^d}
\big\{ \langle w_t, x-x_t\rangle + \frac{1}{2\gamma}\|x-x_t\|^2 \big\}$, we have
$w_t = \frac{1}{\gamma}(x_t - \tilde{x}_{t+1})$, then we can get
\begin{align} \label{eq:F14}
 \|\nabla F(x_t)\| \leq \|w_t\| + \|w_t-\nabla F(x_t)\| =
 \frac{1}{\gamma}\|x_t - \tilde{x}_{t+1}\| + \|w_t-\nabla F(x_t)\|.
\end{align}
Putting the above inequalities~(\ref{eq:F14}) into~(\ref{eq:F12}), we can get
 \begin{align}
 \frac{1}{T}\sum_{t=1}^T\mathbb{E}\|\nabla F(x_t)\|
  & \leq \frac{1}{T} \sum_{t=1}^T \mathbb{E} \big[ \|w_t-\nabla F(x_t)\| + \frac{1}{\gamma}\|\tilde{x}_{t+1}-x_t\|  \big]  \nonumber \\
 & \leq \frac{\sqrt{2M}(m+T)^{1/4}}{\sqrt{T}} \leq \frac{\sqrt{2M}m^{1/4}}{\sqrt{T}} + \frac{\sqrt{2M}}{T^{1/4}}.
\end{align}

\end{proof}

\subsection{Convergence Analysis of VR-MSGBiO Algorithm}
In this subsection, we detail the convergence analysis of VR-MSGBiO algorithm.

\begin{lemma} \label{lem:C3}
 Assume that the stochastic partial derivatives $u_{t+1}$, $h_{t+1}$, $v_{t+1}$, $G_{t+1}$ and $H_{t+1}$ be generated from Algorithm \ref{alg:3}, we have
 \begin{align}
 \mathbb{E}\|\nabla_x f(x_{t+1},y_{t+1}) - u_{t+1}\|^2
 & \leq (1-\beta_{t+1})^2\mathbb{E} \|\nabla_x f(x_t,y_t) - u_t\|^2 + 2\beta_{t+1}^2\sigma^2  \\
 & \quad + 4(1-\beta_{t+1})^2L_f^2\eta_t^2\big(\mathbb{E}\|\tilde{x}_{t+1} - x_t\|^2 + \mathbb{E}\|\tilde{y}_{t+1} - y_t\|^2 \big), \nonumber
 \end{align}
 \begin{align}
 \mathbb{E}\|\Pi_{C_{fy}}\big[\nabla_y f(x_{t+1},y_{t+1})\big] - h_{t+1}\|^2
 & \leq (1-\hat{\beta}_{t+1})^2\mathbb{E} \|\Pi_{C_{fy}}\big[\nabla_y f(x_t,y_t)\big] - h_t\|^2 + 2\hat{\beta}_{t+1}^2\sigma^2  \\
 & \quad + 4(1-\hat{\beta}_{t+1})^2L_f^2\eta_t^2\big(\mathbb{E}\|\tilde{x}_{t+1} - x_t\|^2 + \mathbb{E}\|\tilde{y}_{t+1} - y_t\|^2 \big), \nonumber
 \end{align}
 \begin{align}
 \mathbb{E}\|\nabla_y g(x_{t+1},y_{t+1}) - v_{t+1}\|^2
 & \leq (1-\alpha_{t+1})^2\mathbb{E} \|\nabla_y g(x_t,y_t) - v_t\|^2 + 2\alpha_{t+1}^2\sigma^2  \\
 & \quad + 4(1-\alpha_{t+1})^2L_g^2\eta_t^2\big(\mathbb{E}\|\tilde{x}_{t+1} - x_t\|^2 + \mathbb{E}\|\tilde{y}_{t+1} - y_t\|^2 \big), \nonumber
 \end{align}
 \begin{align}
 \mathbb{E}\|\hat{\Pi}_{C_{gxy}}\big[\nabla^2_{xy}g(x_{t+1},y_{t+1})\big] - G_{t+1}\|^2
 & \leq (1-\hat{\alpha}_{t+1})^2\mathbb{E} \|\hat{\Pi}_{C_{gxy}}\big[\nabla^2_{xy}g(x_t,y_t)\big] - G_t\|^2 + 2\hat{\alpha}_{t+1}^2\sigma^2  \\
 & \quad + 4(1-\hat{\alpha}_{t+1})^2L_{gxy}^2\eta_t^2\big(\mathbb{E}\|\tilde{x}_{t+1} - x_t\|^2 + \mathbb{E}\|\tilde{y}_{t+1} - y_t\|^2 \big), \nonumber
 \end{align}
 \begin{align}
 \mathbb{E}\|\mathcal{S}_{[\mu,L_g]}\big[\nabla^2_{yy} g(x_{t+1},y_{t+1})\big] - H_{t+1}\|^2
 & \leq (1-\tilde{\alpha}_{t+1})^2\mathbb{E} \|\mathcal{S}_{[\mu,L_g]}\big[\nabla^2_{yy} g(x_t,y_t)\big] - H_t\|^2 + 2\tilde{\alpha}_{t+1}\sigma^2  \\
 & \quad + 4(1-\tilde{\alpha}_{t+1})^2L_{gyy}^2\eta_t^2\big(\mathbb{E}\|\tilde{x}_{t+1} - x_t\|^2 + \mathbb{E}\|\tilde{y}_{t+1} - y_t\|^2 \big). \nonumber
 \end{align}
\end{lemma}
\begin{proof}
Without loss of generality, we only consider the term $\mathbb{E}\|\Pi_{C_{fy}}\big[\nabla_y f(x_{t+1},y_{t+1})\big]-h_{t+1}\|^2$. The other terms are similar for this term.
Since $h_{t+1} = \Pi_{C_{fy}}\big[ \nabla_y f(x_{t+1},y_{t+1};\xi_{t+1}) + (1-\hat{\beta}_{t+1})\big(h_t - \nabla_y f(x_t,y_t;\xi_{t+1})\big)\big]$, we have
\begin{align}
 &\mathbb{E}\|\Pi_{C_{fy}}\big[\nabla_y f(x_{t+1},y_{t+1})\big]-h_{t+1}\|^2 \nonumber \\
 & = \mathbb{E}\|\Pi_{C_{fy}}\big[\nabla_y f(x_{t+1},y_{t+1})\big] - \Pi_{C_{fy}}\big[ \nabla_y f(x_{t+1},y_{t+1};\xi_{t+1}) + (1-\hat{\beta}_{t+1})\big(h_t - \nabla_y f(x_t,y_t;\xi_{t+1})\big)\big]\|^2 \nonumber \\
 & \leq \mathbb{E}\|\nabla_y f(x_{t+1},y_{t+1}) -  \nabla_y f(x_{t+1},y_{t+1};\xi_{t+1}) - (1-\hat{\beta}_{t+1})\big(h_t - \nabla_y f(x_t,y_t;\xi_{t+1})\big)\|^2 \nonumber \\
 & = \mathbb{E}\|(1-\hat{\beta}_{t+1})(\nabla_y f(x_t,y_t) - h_t) + \hat{\beta}_{t+1}(\nabla_y f(x_{t+1},y_{t+1}) - \nabla_y f(x_{t+1},y_{t+1};\xi_{t+1}))\nonumber \\
 & \quad - (1-\hat{\beta}_{t+1})\big( \nabla_y f(x_{t+1}, y_{t+1}; \xi_{t+1}) - \nabla_y f(x_t,y_t;\xi_{t+1}) - (\nabla_y f(x_{t+1},y_{t+1}) - \nabla_y f(x_t,y_t))\big)\|^2 \nonumber \\
 & \mathop{=}^{(i)}(1-\hat{\beta}_{t+1})^2\mathbb{E} \|\nabla_y f(x_t,y_t) - h_t\|^2 + \mathbb{E}\big[\|\hat{\beta}_{t+1}\big(\nabla_y f(x_{t+1},y_{t+1}) - \nabla_y f(x_{t+1},y_{t+1};\xi_{t+1})\big)\nonumber \\
 & \quad - (1-\hat{\beta}_{t+1})\big( \nabla_y f(x_{t+1}, y_{t+1}; \xi_{t+1}) - \nabla_y f(x_t,y_t;\xi_{t+1}) - (\nabla_y f(x_{t+1},y_{t+1})-\nabla_y f(x_t,y_t))\big) \|^2\big] \nonumber \\
 & \leq (1-\hat{\beta}_{t+1})^2\mathbb{E} \|\nabla_y f(x_t,y_t) - h_t\|^2 + 2\hat{\beta}_{t+1}^2\mathbb{E}\|\big(\nabla_y f(x_{t+1},y_{t+1}) - \nabla f(x_{t+1},y_{t+1};\xi_{t+1})\|^2 \nonumber \\
 & \quad + 2(1-\hat{\beta}_{t+1})^2\| \nabla_y f(x_{t+1}, y_{t+1}; \xi_{t+1}) - \nabla_y f(x_t,y_t;\xi_{t+1}) - (\nabla_y f(x_{t+1},y_{t+1}) - \nabla_y f(x_t,y_t)) \|^2 \nonumber \\
 & \mathop{\leq}^{(ii)} (1-\hat{\beta}_{t+1})^2\mathbb{E} \|\nabla_y f(x_t,y_t) - h_t\|^2 + 2\hat{\beta}_{t+1}^2\sigma^2 + 2(1-\hat{\beta}_{t+1})^2\|\nabla_y f(x_{t+1},y_{t+1};\xi_{t+1}) - \nabla_y f(x_t,y_t;\xi_{t+1})\|^2 \nonumber \\
 & \mathop{\leq}^{(iii)} (1-\hat{\beta}_{t+1})^2\mathbb{E} \|\nabla_y f(x_t,y_t) - h_t\|^2 + 4(1-\hat{\beta}_{t+1})^2L_f^2\eta_t^2\big(\mathbb{E}\|\tilde{x}_{t+1} - x_t\|^2 + \mathbb{E}\|\tilde{y}_{t+1} - y_t\|^2 \big) \nonumber \\
 & \quad + 2\hat{\beta}_{t+1}^2\sigma^2, \nonumber
\end{align}
where the above equality (i) is due to $\mathbb{E}_{\xi_{t+1}}[\nabla f(x_{t+1},y_{t+1};\xi_{t+1})]=\nabla f(x_{t+1},y_{t+1})$ and  $\mathbb{E}_{\xi_{t+1}}[\nabla f(x_{t+1},y_{t+1};\xi_{t+1})-\nabla f(x_t,y_t;\xi_{t+1})]=\nabla f(x_{t+1},y_{t+1})-\nabla f(x_t,y_t)$; the inequality (ii) holds by Assumption~\ref{ass:8} and the inequality $\mathbb{E}\|\zeta - \mathbb{E}[\zeta]\|^2 = \mathbb{E}\|\zeta\|^2 - (\mathbb{E}[\zeta])^2\leq \mathbb{E}\|\zeta\|^2 $, and the inequality (iii) holds by Assumption~\ref{ass:6} and $x_{t+1}=x_t-\eta_t(\tilde{x}_{t+1}-x_t)$, $y_{t+1}=y_t-\eta_t(\tilde{y}_{t+1}-y_t)$.

\end{proof}

\begin{theorem}  \label{th:A3}
(Restatement of Theorem 3)
 Under the above Assumptions (\ref{ass:1}, \ref{ass:2}, \ref{ass:5}, \ref{ass:6}-\ref{ass:8}), in the Algorithm \ref{alg:3}, let $\eta_t=\frac{k}{(m+t)^{1/3}}$ for all $t\geq 0$, $\beta_{t+1}=c_1\eta_t^2$, $\hat{\beta}_{t+1}=c_2\eta_t^2$, $\alpha_{t+1}=c_3\eta_t^2$, $\hat{\alpha}_{t+1}=c_4\eta_t^2$, $\tilde{\alpha}_{t+1}=c_5\eta_t^2$, $m\geq \max\big(2,k^3, (c_1k)^3,(c_2k)^3, (c_3k)^3, (c_4k)^3, (c_5k)^3\big)$, $k>0$, $c_1 \geq \frac{2}{3k^3} + 10$, $c_2 \geq \frac{2}{3k^3} + 10\kappa^2$, $c_3 \geq \frac{2}{3k^3} + 1$, $c_4 \geq \frac{2}{3k^3} + \frac{10C^2_{fy}}{\mu^2}$, $c_5 \geq \frac{2}{3k^3} + +  \frac{10C^2_{fy}\kappa^2}{\mu^2}$, $0< \gamma \leq \min\big(\frac{m^{1/3}}{2L_Fk},\frac{\lambda\mu^2}{16\hat{L}^2},\frac{1}{8\check{L}},\frac{1}{64L^2_g\lambda},
 \frac{1}{32\check{L}^2\lambda},\frac{\lambda\mu}{8L_G}\big)$, $0<\lambda\leq \min\big(\frac{1}{4\sqrt{2}L_g},\frac{m^{1/3}}{2L_gk}\big)$. When $\mathcal{X}\subseteq\mathbb{R}^d$, we can get
 \begin{align}
 \frac{1}{T}\sum_{t=1}^T\mathbb{E}\|\mathcal{G}(x_t,\nabla F(x_t),\gamma)\|
 \leq \frac{1}{T} \sum_{t=1}^T \mathbb{E} \big[ \|w_t-\nabla F(x_t)\| + \frac{1}{\gamma}\|\tilde{x}_{t+1}-x_t\|  \big]  \leq \frac{\sqrt{2\breve{M}}m^{1/6}}{\sqrt{T}} + \frac{\sqrt{2\breve{M}}}{T^{1/3}};
\end{align}
When $\mathcal{X}=\mathbb{R}^d$, we can get
 \begin{align}
 \frac{1}{T}\sum_{t=1}^T\mathbb{E}\|\nabla F(x_t)\|
  \leq \frac{1}{T} \sum_{t=1}^T \mathbb{E} \big[ \|w_t-\nabla F(x_t)\| + \frac{1}{\gamma}\|\tilde{x}_{t+1}-x_t\|  \big]   \leq \frac{\sqrt{2\breve{M}}m^{1/6}}{\sqrt{T}} + \frac{\sqrt{2\breve{M}}}{T^{1/3}},
\end{align}
where $\breve{M}= \frac{4(F(x_1)- F^*+g(x_1,y_1)-G(x_1))}{k\gamma} + \frac{16\sigma^2m^{1/3}}{k^2} + \frac{4\lambda\sigma^2m^{1/3}}{\gamma k^2} + \big( 2k^2\hat{c}^2\sigma^2 + \frac{2k^2c^2_3\lambda\sigma^2}{\gamma}\big)\ln(m+T)$,
 $\check{L}^2=2L^2_f+L^2_{gxy}+L^2_{gyy}$ and $\hat{c}^2=c^2_1+c^2_2+c^2_4+c^2_5$.
\end{theorem}

\begin{proof}
According the line 4 of Algorithm~\ref{alg:3}, we have
\begin{align} \label{eq:G1}
\tilde{x}_{t+1} = \arg\min_{x\in \mathcal{X}}
\Big\{ \langle w_t, x-x_t\rangle + \frac{1}{2\gamma}\|x-x_t\|^2 \Big\}.
\end{align}
By using the optimal condition of the subproblem~(\ref{eq:G1}), we have
\begin{align}
 \langle w_t + \frac{1}{\gamma}(\tilde{x}_{t+1}-x_t), x_t-\tilde{x}_{t+1}\rangle \geq 0,
\end{align}
then we have
\begin{align}
 \langle w_t,\tilde{x}_{t+1} -x_t\rangle \leq -\frac{1}{\gamma}\|\tilde{x}_{t+1}-x_t\|^2.
\end{align}
According to Lemma \ref{lem:C1}, i.e., the function $F(x)$ has $L_F$-Lipschitz continuous gradient,
we have
\begin{align} \label{eq:G2}
  F(x_{t+1}) & \leq F(x_t) + \langle \nabla F(x_t), x_{t+1}-x_t\rangle + \frac{L_F}{2}\|x_{t+1}-x_t\|^2 \nonumber \\
  & = F(x_t) + \eta_t \langle \nabla F(x_t) -w_t +w_t, \tilde{x}_{t+1}-x_t\rangle + \frac{\eta^2_tL_F}{2}\|\tilde{x}_{t+1}-x_t\|^2 \nonumber \\
  & \leq F(x_t) + \eta_t \langle \nabla F(x_t) -w_t , \tilde{x}_{t+1}-x_t\rangle - \frac{\eta_t}{\gamma}\|\tilde{x}_{t+1}-x_t\|^2 + \frac{\eta^2_tL_F}{2}\|\tilde{x}_{t+1}-x_t\|^2 \nonumber \\
  & \leq F(x_t) + \eta_t\gamma\|\nabla F(x_t) -w_t\|^2 + \frac{\eta_t}{4\gamma}\|\tilde{x}_{t+1}-x_t\|^2 - \frac{\eta_t}{\gamma}\|\tilde{x}_{t+1}-x_t\|^2 + \frac{\eta^2_tL_F}{2}\|\tilde{x}_{t+1}-x_t\|^2 \nonumber \\
  & \leq F(x_t) + \eta_t\gamma\|\nabla F(x_t) -w_t\|^2 - \frac{\eta_t}{2\gamma}\|\tilde{x}_{t+1}-x_t\|^2,
\end{align}
where the second inequality holds by the above inequality~(\ref{eq:G1}), and the last inequality holds by $\gamma\leq \frac{1}{2L_F\eta_t}$.

According to the Lemma~\ref{lem:B2}, we have
\begin{align} \label{eq:G3}
  F(x_{t+1})  & \leq F(x_t) + 8\eta_t\gamma\|u_t - \nabla_xf(x_t,y_t)\|^2 + \frac{8C^2_{fy}\eta_t\gamma}{\mu^2}\|G_t-\hat{\Pi}_{C_{gxy}}\big[\nabla^2_{xy}g(x_t,y_t)\big]\|^2 \nonumber \\
 & \quad + \frac{8\kappa^2 C^2_{fy}\eta_t\gamma}{\mu^2}\|H_t-
 \mathcal{S}_{[\mu,L_g]}\big[\nabla^2_{yy} g(x_t,y_t)\big]\|^2 + 8\kappa^2\eta_t\gamma \|h_t -\Pi_{C_{fy}}\big[\nabla_yf(x_t,y_t)\big] \|^2 \nonumber \\
 & \quad + \frac{4\hat{L}^2\eta_t\gamma}{\mu}\big(g(x_t,y_t)-G(x_t)\big) - \frac{\eta_t}{2\gamma}\|\tilde{x}_{t+1}-x_t\|^2.
\end{align}

According to Lemma~\ref{lem:E1}, we have
\begin{align}
& g(x_{t+1},y_{t+1})-G(x_{t+1})-\big(g(x_t,y_t)-G(x_t)\big) \\
& \leq -\frac{\eta_t\lambda\mu}{2} \big(g(x_t,y_t) -G(x_t)\big) + \frac{\eta_t}{8\gamma}\|\tilde{x}_{t+1}-x_t\|^2   -\frac{\eta_t}{4\lambda}\|\tilde{y}_{t+1}-y_t\|^2+ \eta_t\lambda\|\nabla_y g(x_t,y_t)-v_t\|^2. \nonumber
\end{align}

In our Algorithm \ref{alg:3}, we set $\eta_t=\frac{k}{(m+t)^{1/3}}$ for all $t\geq 0$, $\beta_{t+1}=c_1\eta_t^2$, $\hat{\beta}_{t+1}=c_2\eta_t^2$, $\alpha_{t+1}=c_3\eta_t^2$, $\hat{\alpha}_{t+1}=c_4\eta_t^2$, $\tilde{\alpha}_{t+1}=c_5\eta_t^2$.
Since $\eta_t=\frac{k}{(m+t)^{1/3}}$ on $t$ is decreasing and $m\geq k^3$, we have $\eta_t \leq \eta_0 = \frac{k}{m^{1/3}} \leq 1$ and $\gamma \leq \frac{m^{1/3}}{2L_Fk}\leq \frac{1}{2L_F\eta_0} \leq \frac{1}{2L_F\eta_t}$ for any $t\geq 0$. Similarly, we can get $\gamma \leq \frac{m^{1/3}}{2L_gk}\leq \frac{1}{2L_g\eta_t}$ for any $t\geq 0$. Due to $0 < \eta_t \leq 1$ and $m\geq (c_1k)^3$, we have $0<\beta_{t+1} = c_1\eta_t^2 \leq c_1\eta_t \leq \frac{c_1k}{m^{1/3}}\leq 1$.
 Similarly, due to $m\geq \max\big( (c_2k)^3, (c_3k)^3, (c_4k)^3, (c_5k)^3 \big)$, we have $0<\hat{\beta}_{t+1}\leq 1$, $0<\alpha_{t+1}\leq 1$, $0<\hat{\alpha}_{t+1}\leq 1$ and $0<\tilde{\alpha}_{t+1} \leq 1$.

 According to Lemma \ref{lem:C3}, we have
 \begin{align}
  & \frac{1}{\eta_t}\mathbb{E} \|\nabla_x f(x_{t+1},y_{t+1}) - u_{t+1}\|^2 - \frac{1}{\eta_{t-1}}\mathbb{E} \|\nabla_x f(x_t,y_t) - u_t\|^2  \\
  & \leq \big(\frac{1}{\eta_t} - \frac{1}{\eta_{t-1}} - c_1\eta_t\big)\mathbb{E} \|\nabla_x f(x_t,y_t) -u_t\|^2  + 4L^2_f\eta_t\big(\|\tilde{x}_{t+1}-x_t\|^2 +\|\tilde{y}_{t+1}-y_t\|^2\big) + 2c_1^2\eta^3_t\sigma^2, \nonumber
 \end{align}
 where the above inequality is due to $0< \beta_{t+1}=c_1\eta^2_t \leq 1$.
 Similarly, according to Lemma \ref{lem:C3}, we can obtain
 \begin{align}
  & \frac{1}{\eta_t}\mathbb{E} \|\Pi_{C_{fy}}\big[\nabla_y f(x_{t+1},y_{t+1})\big] - h_{t+1}\|^2
  -  \frac{1}{\eta_{t-1}}\mathbb{E} \|\Pi_{C_{fy}}\big[\nabla_y f(x_t,y_t)\big] - h_t\|^2  \\
  & \leq \big(\frac{1}{\eta_t} - \frac{1}{\eta_{t-1}} - c_2\eta_t\big)\mathbb{E} \|\Pi_{C_{fy}}\big[\nabla_y f(x_t,y_t)\big] -h_t\|^2 +
  4L^2_f\eta_t\big(\|\tilde{x}_{t+1}-x_t\|^2 + \|\tilde{y}_{t+1}-y_t\|^2\big) + 2c_2^2\eta_t^3\sigma^2, \nonumber
 \end{align}
 \begin{align}
  & \frac{1}{\eta_t}\mathbb{E} \|\nabla_y g(x_{t+1},y_{t+1}) - v_{t+1}\|^2
  -  \frac{1}{\eta_{t-1}}\mathbb{E} \|\nabla_y g(x_t,y_t) - v_t\|^2  \\
  & \leq \big(\frac{1}{\eta_t} - \frac{1}{\eta_{t-1}} - c_3\eta_t\big)\mathbb{E} \|\nabla_y g(x_t,y_t) -v_t\|^2 +
  4L^2_g\eta_t\big(\|\tilde{x}_{t+1}-x_t\|^2 + \|\tilde{y}_{t+1}-y_t\|^2\big) + 2c_3^2\eta_t^3\sigma^2, \nonumber
 \end{align}
 \begin{align}
  & \frac{1}{\eta_t}\mathbb{E} \|\hat{\Pi}_{C_{gxy}}\big[\nabla^2_{xy} g(x_{t+1},y_{t+1})\big] - G_{t+1}\|^2
  -  \frac{1}{\eta_{t-1}}\mathbb{E} \|\hat{\Pi}_{C_{gxy}}\big[\nabla^2_{xy} g(x_t,y_t)\big] - G_t\|^2  \\
  & \leq \big(\frac{1}{\eta_t} - \frac{1}{\eta_{t-1}} - c_4\eta_t\big)\mathbb{E} \|\hat{\Pi}_{C_{gxy}}\big[\nabla^2_{xy} g(x_t,y_t)\big] -G_t\|^2 +
  4L^2_{gxy}\eta_t\big(\|\tilde{x}_{t+1}-x_t\|^2 + \|\tilde{y}_{t+1}-y_t\|^2\big) + 2c_4^2\eta_t^3\sigma^2, \nonumber
 \end{align}
 and
 \begin{align}
  & \frac{1}{\eta_t}\mathbb{E} \|\mathcal{S}_{[\mu,L_g]}\big[\nabla^2_{yy} g(x_{t+1},y_{t+1})\big] - H_{t+1}\|^2
  -  \frac{1}{\eta_{t-1}}\mathbb{E} \|\mathcal{S}_{[\mu,L_g]}\big[\nabla^2_{yy} g(x_t,y_t)\big] - H_t\|^2  \\
  & \leq \big(\frac{1}{\eta_t} - \frac{1}{\eta_{t-1}} - c_5\eta_t\big)\mathbb{E} \|\mathcal{S}_{[\mu,L_g]}\big[\nabla^2_{yy} g(x_t,y_t)\big] -H_t\|^2 +
  4L^2_{gyy}\eta_t\big(\|\tilde{x}_{t+1}-x_t\|^2 + \|\tilde{y}_{t+1}-y_t\|^2\big) + 2c_5^2\eta_t^3\sigma^2. \nonumber
 \end{align}
 By using $\eta_t = \frac{k}{(m+t)^{1/3}}$, we have
 \begin{align}
  \frac{1}{\eta_t} - \frac{1}{\eta_{t-1}} &= \frac{1}{k}\big( (m+t)^{\frac{1}{3}} - (m+t-1)^{\frac{1}{3}}\big) \leq \frac{1}{3k(m+t-1)^{2/3}} \leq \frac{1}{3k\big(m/2+t\big)^{2/3}} \nonumber \\
  & \leq \frac{2^{2/3}}{3k(m+t)^{2/3}} = \frac{2^{2/3}}{3k^3}\frac{k^2}{(m+t)^{2/3}} = \frac{2^{2/3}}{3k^3}\eta_t^2 \leq \frac{2}{3k^3}\eta_t,
 \end{align}
 where the first inequality holds by the concavity of function $f(x)=x^{1/3}$, \emph{i.e.}, $(x+y)^{1/3}\leq x^{1/3} + \frac{y}{3x^{2/3}}$; the second inequality is due to $m\geq 2$, and the last inequality is due to $0<\eta_t\leq 1$.
 Let $c_1 \geq \frac{2}{3k^3} + 10$, we have
 \begin{align} \label{eq:G4}
  & \frac{1}{\eta_t}\mathbb{E} \|\nabla_x f(x_{t+1},y_{t+1}) - u_{t+1}\|^2 - \frac{1}{\eta_{t-1}}\mathbb{E} \|\nabla_x f(x_t,y_t) - u_t\|^2  \\
  & \leq -10\eta_t\mathbb{E} \|\nabla_x f(x_t,y_t) -u_t\|^2 + 4L^2_f\eta_t\big(\|\tilde{x}_{t+1}-x_t\|^2 +\|\tilde{y}_{t+1}-y_t\|^2\big) + 2c_1^2\eta^3_t\sigma^2. \nonumber
 \end{align}
 Let $c_2 \geq \frac{2}{3k^3} + 10\kappa^2$, we have
 \begin{align} \label{eq:G5}
  & \frac{1}{\eta_t}\mathbb{E} \|\Pi_{C_{fy}}\big[\nabla_y f(x_{t+1},y_{t+1})\big] - h_{t+1}\|^2
  -  \frac{1}{\eta_{t-1}}\mathbb{E} \|\Pi_{C_{fy}}\big[\nabla_y f(x_t,y_t)\big] - h_t\|^2   \\
  & \leq -10\kappa^2\eta_t\mathbb{E} \|\Pi_{C_{fy}}\big[\nabla_y f(x_t,y_t)\big] -h_t\|^2 +
  4L^2_f\eta_t\big(\|\tilde{x}_{t+1}-x_t\|^2 +\|\tilde{y}_{t+1}-y_t\|^2\big)  + 2c_2^2\eta_t^3\sigma^2. \nonumber
 \end{align}
 Let $c_3 \geq \frac{2}{3k^3} + 1$, we have
 \begin{align} \label{eq:G6}
  & \frac{1}{\eta_t}\mathbb{E} \|\nabla_y g(x_{t+1},y_{t+1}) - v_{t+1}\|^2
  -  \frac{1}{\eta_{t-1}}\mathbb{E} \|\nabla_y g(x_t,y_t) - v_t\|^2  \\
  & \leq - \eta_t\mathbb{E} \|\nabla_y g(x_t,y_t) -v_t\|^2 +
  4L^2_g\eta_t\big(\|\tilde{x}_{t+1}-x_t\|^2 + \|\tilde{y}_{t+1}-y_t\|^2\big) + 2c_3^2\eta_t^3\sigma^2. \nonumber
 \end{align}
 Let $c_4 \geq \frac{2}{3k^3} +  \frac{10C^2_{fy}}{\mu^2}$, we have
 \begin{align} \label{eq:G7}
  &  \frac{1}{\eta_t}\mathbb{E} \|\hat{\Pi}_{C_{gxy}}\big[\nabla^2_{xy} g(x_{t+1},y_{t+1})\big] - G_{t+1}\|^2
  -  \frac{1}{\eta_{t-1}}\mathbb{E} \|\hat{\Pi}_{C_{gxy}}\big[\nabla^2_{xy} g(x_t,y_t)\big] - G_t\|^2   \\
  & \leq -\frac{10C^2_{fy}\eta_t}{\mu^2} \mathbb{E} \|\hat{\Pi}_{C_{gxy}}\big[\nabla^2_{xy} g(x_t,y_t)\big] -G_t\|^2 +
  4L^2_{gxy}\eta_t\big(\|\tilde{x}_{t+1}-x_t\|^2 + \|\tilde{y}_{t+1}-y_t\|^2\big) + 2c_4^2\eta_t^3\sigma^2. \nonumber
 \end{align}
Let $c_5 \geq \frac{2}{3k^3} + \frac{10C^2_{fy}\kappa^2}{\mu^2}$, we have
 \begin{align} \label{eq:G8}
  & \frac{1}{\eta_t}\mathbb{E} \|\mathcal{S}_{[\mu,L_g]}\big[\nabla^2_{yy} g(x_{t+1},y_{t+1})\big] - H_{t+1}\|^2
  -  \frac{1}{\eta_{t-1}}\mathbb{E} \|\mathcal{S}_{[\mu,L_g]}\big[\nabla^2_{yy} g(x_t,y_t)\big] - H_t\|^2  \\
  & \leq -\frac{10C^2_{fy}\kappa^2\eta_t}{\mu^2}\mathbb{E} \|\mathcal{S}_{[\mu,L_g]}\big[\nabla^2_{yy} g(x_t,y_t)\big] -H_t\|^2 + 4L^2_{gyy}\eta_t\big(\|\tilde{x}_{t+1}-x_t\|^2 + \|\tilde{y}_{t+1}-y_t\|^2\big) + 2c_5^2\eta_t^3\sigma^2. \nonumber
 \end{align}

Next, we define a useful \emph{Lyapunov} function (i.e., potential function), for any $t\geq 1$
 \begin{align}
 \Gamma_t & = \mathbb{E}\Big [F(x_t) + g(x_t,y_t)-G(x_t) + \frac{\gamma}{\eta_{t-1}} \big(\|\nabla_x f(x_t,y_t)-u_t\|^2 + \|\Pi_{C_{fy}}\big[\nabla_y f(x_t,y_t)\big] -h_t\|^2 \nonumber \\
 & \quad  + \|\hat{\Pi}_{C_{gxy}}\big[\nabla^2_{xy} g(x_t,y_t)\big] - G_t\|^2 + \|\mathcal{S}_{[\mu,L_g]}\big[\nabla^2_{yy} g(x_t,y_t)\big] - H_t\|^2\big) + \frac{\lambda}{\eta_{t-1}}\|\nabla_y g(x_t,y_t) - v_t\|^2 \Big]. \nonumber
 \end{align}
 Then we have
 \begin{align} \label{eq:G9}
 & \Gamma_{t+1} - \Gamma_t \nonumber \\
 & = F(x_{t+1}) - F(x_t) + g(x_{t+1},y_{t+1})-G(x_{t+1}) - \big(g(x_t,y_t)-G(x_t)\big)
 + \gamma \Big(\frac{1}{\eta_t}\mathbb{E}\|\nabla_x f(x_{t+1},y_{t+1})-u_{t+1}\|^2 \nonumber \\
 & \quad - \frac{1}{\eta_{t-1}}\mathbb{E}\|\nabla_x f(x_t,y_t)-u_t\|^2 + \frac{1}{\eta_t}\mathbb{E}\|\Pi_{C_{fy}}\big[\nabla_y f(x_{t+1},y_{t+1})\big] -h_{t+1}\|^2
 - \frac{1}{\eta_{t-1}}\mathbb{E}\|\Pi_{C_{fy}}\big[\nabla_y f(x_t,y_t)\big] -h_t\|^2 \nonumber \\
 & \quad  +  \frac{1}{\eta_t}\mathbb{E}\|\hat{\Pi}_{C_{gxy}}\big[\nabla^2_{xy} g(x_{t+1},y_{t+1})\big] -G_{t+1}\|^2
  - \frac{1}{\eta_{t-1}}\mathbb{E}\|\hat{\Pi}_{C_{gxy}}\big[\nabla^2_{xy} g(x_t,y_t)\big] -G_t\|^2 \nonumber \\
 & \quad + \frac{1}{\eta_t}\mathbb{E}\|\mathcal{S}_{[\mu,L_g]}\big[\nabla^2_{yy} g(x_{t+1},y_{t+1})\big]-H_{t+1}\|^2
 - \frac{1}{\eta_{t-1}}\mathbb{E}\|\mathcal{S}_{[\mu,L_g]}\big[\nabla^2_{yy} g(x_t,y_t)\big]-H_t\|^2 \Big) \nonumber \\
 & \quad
  +  \lambda\big(\frac{1}{\eta_t}\mathbb{E}\|\nabla_y g(x_{t+1},y_{t+1})-v_{t+1}\|^2
 - \frac{1}{\eta_{t-1}}\mathbb{E}\|\nabla_y g(x_t,y_t)-v_t\|^2\big) \nonumber \\
 & \leq 8\eta_t\gamma\|u_t - \nabla_xf(x_t,y_t)\|^2 + \frac{8C^2_{fy}\eta_t\gamma}{\mu^2}\|G_t-\hat{\Pi}_{C_{gxy}}\big[\nabla^2_{xy}g(x_t,y_t)\big]\|^2 + \frac{8\kappa^2 C^2_{fy}\eta_t\gamma}{\mu^2}\|H_t-
 \mathcal{S}_{[\mu,L_g]}\big[\nabla^2_{yy} g(x_t,y_t)\big]\|^2 \nonumber \\
 & \quad + 8\kappa^2\eta_t\gamma \|h_t -\Pi_{C_{fy}}\big[\nabla_yf(x_t,y_t)\big] \|^2 + \frac{4\hat{L}^2\eta_t\gamma}{\mu}\big(g(x_t,y_t)-G(x_t)\big) - \frac{\eta_t}{2\gamma}\|\tilde{x}_{t+1}-x_t\|^2 \nonumber \\
 & \quad -\frac{\eta_t\lambda\mu}{2} \big(g(x_t,y_t) -G(x_t)\big) + \frac{\eta_t}{8\gamma}\|\tilde{x}_{t+1}-x_t\|^2   -\frac{\eta_t}{4\lambda}\|\tilde{y}_{t+1}-y_t\|^2+ \eta_t\lambda\|\nabla_y g(x_t,y_t)-v_t\|^2  \nonumber \\
 & \quad + \gamma \bigg( -10\eta_t\mathbb{E} \|\nabla_x f(x_t,y_t) -u_t\|^2 + 4L^2_f\eta_t\big(\|\tilde{x}_{t+1}-x_t\|^2 +\|\tilde{y}_{t+1}-y_t\|^2\big) + 2c_1^2\eta^3_t\sigma^2 \nonumber \\
 & \qquad -10\kappa^2\eta_t\mathbb{E} \|\Pi_{C_{fy}}\big[\nabla_y f(x_t,y_t)\big] -h_t\|^2 +
  4L^2_f\eta_t\big(\|\tilde{x}_{t+1}-x_t\|^2 +\|\tilde{y}_{t+1}-y_t\|^2\big)  + 2c_2^2\eta_t^3\sigma^2  \nonumber \\
 & \qquad -\frac{10C^2_{fy}\eta_t}{\mu^2} \mathbb{E} \|\hat{\Pi}_{C_{gxy}}\big[\nabla^2_{xy} g(x_t,y_t)\big] -G_t\|^2 +
  4L^2_{gxy}\eta_t\big(\|\tilde{x}_{t+1}-x_t\|^2 + \|\tilde{y}_{t+1}-y_t\|^2\big) + 2c_4^2\eta_t^3\sigma^2 \nonumber \\
 & \qquad -\frac{10C^2_{fy}\kappa^2\eta_t}{\mu^2}\mathbb{E} \|\mathcal{S}_{[\mu,L_g]}\big[\nabla^2_{yy} g(x_t,y_t)\big] -H_t\|^2 + 4L^2_{gyy}\eta_t\big(\|\tilde{x}_{t+1}-x_t\|^2 + \|\tilde{y}_{t+1}-y_t\|^2\big) + 2c_5^2\eta_t^3\sigma^2 \bigg)\nonumber \\
  & \quad +\lambda\Big(- \eta_t \mathbb{E} \|\nabla_y g(x_t,y_t) -v_t\|^2 +
  4L^2_g\eta_t\big(\|\tilde{x}_{t+1}-x_t\|^2 + \|\tilde{y}_{t+1}-y_t\|^2\big) + 2c_3^2\eta_t^3\sigma^2 \Big) \nonumber \\
 & = - \frac{\gamma\eta_t}{4} \bigg( 8\eta_t\gamma\|u_t - \nabla_xf(x_t,y_t)\|^2 + \frac{8C^2_{fy}\eta_t\gamma}{\mu^2}\|G_t-\hat{\Pi}_{C_{gxy}}\big[\nabla^2_{xy}g(x_t,y_t)\big]\|^2 + \frac{8\kappa^2 C^2_{fy}\eta_t\gamma}{\mu^2}\|H_t-
 \mathcal{S}_{[\mu,L_g]}\big[\nabla^2_{yy} g(x_t,y_t)\big]\|^2 \nonumber \\
 & \qquad + 8\kappa^2\eta_t\gamma \|h_t -\Pi_{C_{fy}}\big[\nabla_yf(x_t,y_t)\big] \|^2 \bigg) -\big(\frac{\lambda\mu\eta_t}{2}-\frac{4\hat{L}^2\gamma\eta_t}{\mu}\big)\big(g(x_t,y_t)-G(x_t)\big)   \nonumber \\
 & \quad - \big( \frac{1}{2\gamma} -\frac{1}{8\gamma}- 4\check{L}^2\gamma - 4L^2_g\lambda\big)\eta_t\|\tilde{x}_{t+1}-x_t\|^2 - \big( \frac{1}{4\lambda} - 4\check{L}^2\gamma - 4L^2_g\lambda\big)\eta_t\|\tilde{y}_{t+1}-y_t\|^2 \nonumber \\
 & \quad +  2\hat{c}^2\gamma\sigma^2\eta^3_t + 2c^2_3\lambda\sigma^2\eta^3_t \nonumber \\
 & \leq - \frac{\gamma\eta_t}{4} \bigg( 8\eta_t\gamma\|u_t - \nabla_xf(x_t,y_t)\|^2 + \frac{8C^2_{fy}\eta_t\gamma}{\mu^2}\|G_t-\hat{\Pi}_{C_{gxy}}\big[\nabla^2_{xy}g(x_t,y_t)\big]\|^2  \nonumber \\
 & \qquad + \frac{8\kappa^2 C^2_{fy}\eta_t\gamma}{\mu^2}\|H_t-
 \mathcal{S}_{[\mu,L_g]}\big[\nabla^2_{yy} g(x_t,y_t)\big]\|^2 + 8\kappa^2\eta_t\gamma \|h_t -\Pi_{C_{fy}}\big[\nabla_yf(x_t,y_t)\big] \|^2 \bigg) - \frac{\eta_t}{4\gamma}\|\tilde{x}_{t+1}-x_t\|^2 \nonumber \\
 & \quad -\frac{\lambda\mu\eta_t}{4}\big(g(x_t,y_t)-G(x_t)\big) + 2\hat{c}^2\gamma\sigma^2\eta^3_t + 2c^2_3\lambda\sigma^2\eta^3_t,
 \end{align}
 where $\check{L}^2=2L^2_f+L^2_{gxy}+L^2_{gyy}$ and $\hat{c}^2=c^2_1+c^2_2+c^2_4+c^2_5$, the first inequality holds by the above inequalities \eqref{eq:G4}-\eqref{eq:G8};
 the last inequality is due to $0< \gamma \leq \min\big(\frac{\lambda\mu^2}{16\hat{L}^2},\frac{1}{8\check{L}},\frac{1}{64L^2_g\lambda},\frac{1}{32\check{L}^2\lambda}\big)$, $0<\lambda\leq \frac{1}{4\sqrt{2}L_g}$.

Since $\gamma\leq \frac{\lambda\mu^2}{16\hat{L}^2}$, according to the above inequality~(\ref{eq:G9}), we have
\begin{align} \label{eq:G10}
 &  \frac{\gamma\eta_t}{4} \bigg( 8\|u_t - \nabla_xf(x_t,y_t)\|^2 + \frac{8C^2_{fy}}{\mu^2}\|G_t-\hat{\Pi}_{C_{gxy}}\big[\nabla^2_{xy}g(x_t,y_t)\big]\|^2 + \frac{8\kappa^2 C^2_{fy}}{\mu^2}\|H_t-
 \mathcal{S}_{[\mu,L_g]}\big[\nabla^2_{yy} g(x_t,y_t)\big]\|^2 \nonumber \\
 & \qquad + 8\kappa^2 \|h_t -\Pi_{C_{fy}}\big[\nabla_yf(x_t,y_t)\big] \|^2 + \frac{16\hat{L}^2}{\mu}\big(g(x_t,y_t)-G(x_t)\big) \bigg)  + \frac{\eta_t}{4\gamma}\|\tilde{x}_{t+1}-x_t\|^2  \nonumber \\
 &  \leq \frac{\gamma\eta_t}{4} \bigg( 8\|u_t - \nabla_xf(x_t,y_t)\|^2 + \frac{8C^2_{fy}}{\mu^2}\|G_t-\hat{\Pi}_{C_{gxy}}\big[\nabla^2_{xy}g(x_t,y_t)\big]\|^2 + \frac{8\kappa^2 C^2_{fy}}{\mu^2}\|H_t-
 \mathcal{S}_{[\mu,L_g]}\big[\nabla^2_{yy} g(x_t,y_t)\big]\|^2 \nonumber \\
 & \qquad + 8\kappa^2 \|h_t - \Pi_{C_{fy}}\big[\nabla_yf(x_t,y_t)\big] \|^2 \bigg) + \frac{\lambda\mu\eta_t}{4}\big(g(x_t,y_t)-G(x_t)\big) + \frac{\eta_t}{4\gamma}\|\tilde{x}_{t+1}-x_t\|^2  \nonumber \\
 &   \leq \Gamma_t -\Gamma_{t+1} + 2\hat{c}^2\gamma\sigma^2\eta^3_t + 2c^2_3\lambda\sigma^2\eta^3_t.
\end{align}

According to Lemma~\ref{lem:B2}, then we have
\begin{align} \label{eq:G11}
 \frac{\gamma\eta_t}{4}\|w_t-\nabla F(x_t)\|^2 + \frac{\eta_t}{4\gamma}\|\tilde{x}_{t+1}-x_t\|^2  \leq \Gamma_t -\Gamma_{t+1} + 2\hat{c}^2\gamma\sigma^2\eta^3_t + 2c^2_3\lambda\sigma^2\eta^3_t.
\end{align}
Taking average over $t=1,2,\cdots,T$ on both sides of \eqref{eq:G11}, we have
\begin{align}
  \frac{1}{T} \sum_{t=1}^T \mathbb{E} \big[ \eta_t\|w_t-\nabla F(x_t)\|^2 + \frac{\eta_t}{\gamma^2}\|\tilde{x}_{t+1}-x_t\|^2  \big]
  \leq   \frac{4(\Gamma_1 - \Gamma_{T+1})}{T\gamma} + \frac{1}{T}\sum_{t=1}^T\big( 2\hat{c}^2\sigma^2 + \frac{2c^2_3\lambda\sigma^2}{\gamma}\big)\eta^3_t. \nonumber
\end{align}

Since $\eta_t=\frac{k}{(m+t)^{1/3}}$ is decreasing on $t$, i.e., $\eta_T^{-1} \geq \eta_t^{-1}$ for any $0\leq t\leq T$, we have
\begin{align}
 & \frac{1}{T} \sum_{t=1}^T \mathbb{E} \big[ \|w_t-\nabla F(x_t)\|^2 + \frac{1}{\gamma^2}\|\tilde{x}_{t+1}-x_t\|^2  \big] \nonumber \\
 & \leq \frac{4(\Gamma_1 - \Gamma_{T+1})}{T\eta_T\gamma} + \frac{1}{\eta_T T}\sum_{t=1}^T\big( 2\hat{c}^2\sigma^2 + \frac{2c^2_3\lambda\sigma^2}{\gamma}\big)\eta^3_t \nonumber \\
 & \leq \frac{4(F(x_1)- F^*+g(x_1,y_1)-G(x_1) + 4\gamma\sigma^2/\eta_0+ \lambda\sigma^2/\eta_0)}{T\eta_T\gamma} +\big( 2\hat{c}^2\sigma^2 + \frac{2c^2_3\lambda\sigma^2}{\gamma}\big)\frac{1}{\eta_T T}\sum_{t=1}^T\eta^3_t \nonumber \\
 & \leq \frac{4(F(x_1)- F^*+g(x_1,y_1)-G(x_1) + 4\gamma\sigma^2/\eta_0+ \lambda\sigma^2/\eta_0)}{T\eta_T\gamma} +\big( 2\hat{c}^2\sigma^2 + \frac{2c^2_3\lambda\sigma^2}{\gamma}\big)\frac{1}{\eta_T T}\int_{1}^T\frac{k^3}{m+t}dt \nonumber \\
 & \leq \frac{4(F(x_1)- F^*+g(x_1,y_1)-G(x_1) + 4\gamma\sigma^2/\eta_0+ \lambda\sigma^2/\eta_0)}{T\eta_T\gamma} + \big( 2\hat{c}^2\sigma^2 + \frac{2c^2_3\lambda\sigma^2}{\gamma}\big)\frac{k^3\ln(m+T)}{\eta_T T} \nonumber \\
 & \leq \Big(\frac{4(F(x_1)- F^*+g(x_1,y_1)-G(x_1))}{k\gamma} + \frac{16\sigma^2m^{1/3}}{k^2} + \frac{4\lambda\sigma^2m^{1/3}}{\gamma k^2} + \big( 2k^2\hat{c}^2\sigma^2 \nonumber \\
 & \qquad + \frac{2k^2c^2_3\lambda\sigma^2}{\gamma}\big)\ln(m+T)\Big)\frac{(m+T)^{1/3}}{T}
\end{align}
where the second inequality holds by Assumptions~\ref{ass:1}-\ref{ass:3}.
Let $\breve{M}= \frac{4(F(x_1)- F^*+g(x_1,y_1)-G(x_1))}{k\gamma} + \frac{16\sigma^2m^{1/3}}{k^2} + \frac{4\lambda\sigma^2m^{1/3}}{\gamma k^2} + \big( 2k^2\hat{c}^2\sigma^2 + \frac{2k^2c^2_3\lambda\sigma^2}{\gamma}\big)\ln(m+T)$, we have
\begin{align}
 \frac{1}{T} \sum_{t=1}^T \mathbb{E} \big[ \|w_t-\nabla F(x_t)\|^2 + \frac{1}{\gamma^2}\|\tilde{x}_{t+1}-x_t\|^2  \big]
  \leq \frac{\breve{M}(m+T)^{1/3}}{T}.
\end{align}
According to the Jensen's inequality, we can get
\begin{align} \label{eq:G12}
 & \frac{1}{T} \sum_{t=1}^T \mathbb{E} \big[ \|w_t-\nabla F(x_t)\| + \frac{1}{\gamma}\|\tilde{x}_{t+1}-x_t\|  \big]  \nonumber \\
 & \leq \Big(\frac{2}{T} \sum_{t=1}^T \mathbb{E} \big[ \|w_t-\nabla F(x_t)\|^2 + \frac{1}{\gamma^2}\|\tilde{x}_{t+1}-x_t\|^2  \big] \Big)^{1/2} \nonumber \\
 & \leq \frac{\sqrt{2\breve{M}}(m+T)^{1/6}}{\sqrt{T}} \leq \frac{\sqrt{2\breve{M}}m^{1/6}}{\sqrt{T}} + \frac{\sqrt{2\breve{M}}}{T^{1/3}}.
\end{align}

When $\mathcal{X}\subset \mathbb{R}^d$ and $\tilde{x}_{t+1}=\mathbb{P}_{\mathcal{X}}(x_t-\gamma w_t)=\arg\min_{x\in \mathcal{X}}
\big\{ \langle w_t, x-x_t\rangle + \frac{1}{2\gamma}\|x-x_t\|^2 \big\}$, we can define the gradient mapping
$\mathcal{G}(x_t,w_t,\gamma) = \frac{1}{\gamma}\big(x_t-\mathbb{P}_{\mathcal{X}}(x_t-\gamma w_t\big)$.
Meanwhile, we can also define a gradient mapping $\mathcal{G}(x_t,\nabla F(x_t),\gamma) = \frac{1}{\gamma}\big(x_t-\mathbb{P}_{\mathcal{X}}(x_t-\gamma\nabla F(x_t))\big)$.
Then we have
\begin{align} \label{eq:G13}
 \|\mathcal{G}(x_t,\nabla F(x_t),\gamma)\| & \leq \|\mathcal{G}(x_t,w_t,\gamma)\| + \|\mathcal{G}(x_t,\nabla F(x_t),\gamma)-\mathcal{G}(x_t,w_t,\gamma)\| \nonumber \\
 & = \frac{1}{\gamma} \|\tilde{x}_{t+1}-x_t\| + \frac{1}{\gamma}\|\mathbb{P}_{\mathcal{X}}(x_t-\gamma w_t-\mathbb{P}_{\mathcal{X}}(x_t-\gamma\nabla F(x_t))\| \nonumber \\
 & \leq \frac{1}{\gamma} \|\tilde{x}_{t+1}-x_t\| + \|w_t-\nabla F(x_t)\|.
\end{align}
Putting the above inequalities~(\ref{eq:G13}) into~(\ref{eq:G12}), we can obtain
 \begin{align}
 \frac{1}{T}\sum_{t=1}^T\mathbb{E}\|\mathcal{G}(x_t,\nabla F(x_t),\gamma)\|
 & \leq \frac{1}{T} \sum_{t=1}^T \mathbb{E} \big[ \|w_t-\nabla F(x_t)\| + \frac{1}{\gamma}\|\tilde{x}_{t+1}-x_t\|  \big]  \nonumber \\
 & \leq  \frac{\sqrt{2\breve{M}}(m+T)^{1/6}}{\sqrt{T}} \leq \frac{\sqrt{2\breve{M}}m^{1/6}}{\sqrt{T}} + \frac{\sqrt{2\breve{M}}}{T^{1/3}}.
\end{align}

When $\mathcal{X}= \mathbb{R}^d$ and $\tilde{x}_{t+1}=x_t-\gamma w_t=\arg\min_{x\in \mathbb{R}^d}
\big\{ \langle w_t, x-x_t\rangle + \frac{1}{2\gamma}\|x-x_t\|^2 \big\}$, we have
$w_t = \frac{1}{\gamma}(x_t - \tilde{x}_{t+1})$, then we can get
\begin{align} \label{eq:G14}
 \|\nabla F(x_t)\| \leq \|w_t\| + \|w_t-\nabla F(x_t)\| =
 \frac{1}{\gamma}\|x_t - \tilde{x}_{t+1}\| + \|w_t-\nabla F(x_t)\|.
\end{align}
Putting the above inequalities~(\ref{eq:G14}) into~(\ref{eq:G12}), we can get
 \begin{align}
 \frac{1}{T}\sum_{t=1}^T\mathbb{E}\|\nabla F(x_t)\|
  & \leq \frac{1}{T} \sum_{t=1}^T \mathbb{E} \big[ \|w_t-\nabla F(x_t)\| + \frac{1}{\gamma}\|\tilde{x}_{t+1}-x_t\|  \big]  \nonumber \\
 & \leq \frac{\sqrt{2\breve{M}}(m+T)^{1/6}}{\sqrt{T}} \leq \frac{\sqrt{2\breve{M}}m^{1/6}}{\sqrt{T}} + \frac{\sqrt{2\breve{M}}}{T^{1/3}}.
\end{align}

\end{proof}

\end{document}